\documentclass[10pt,reqno]{amsart}
\usepackage{hyperref}

\allowdisplaybreaks[4]
\usepackage{amsmath}
\usepackage{amssymb}
\usepackage{mathrsfs}
\usepackage{amsthm}
\usepackage{mathtools}
\usepackage{stmaryrd}
\usepackage{bm}
\usepackage{subfigure}
\usepackage{graphicx}
\usepackage{booktabs}
\usepackage{float}
\usepackage{geometry}
\usepackage{color}
\usepackage{enumitem}
\theoremstyle{plain}
\newtheorem{lemma}{Lemma}[section]
\newtheorem{thm}[lemma]{Theorem}

\newtheorem{coro}[lemma]{Corollary}
\newtheorem{remark}[lemma]{Remark}

\newtheorem{assumption}[lemma]{Assumption}
\numberwithin{equation}{section}
\DeclarePairedDelimiter{\normtriple}{\interleave}{\interleave}

\begin{document}
\title{Spherical Perfectly Matched Layer Method for Electromagnetic Scattering Problems with Large Wave Number}
\author[Jian Shi]{Jian Shi$^{*}$}
\author[Yonglin Li]{Yonglin Li$^{\dagger}$}
\author[Haijun Wu]{Haijun Wu$^{\ddagger}$}
\thanks{The research of Jian Shi and Haijun Wu was supported
in part by the National Key R\&D Program of China
(No. 2024YFA1012600) and NSF of China
(Nos. 12171238 and 12261160361).
The research of Yonglin Li was supported in part by NSF of China
(No. 12501562), and Hubei Provincial NSF (No. 2025AFB272).}

\thanks{$^{*}$Department of Mathematics, Nanjing University,
Nanjing 210093, China (mathshi@nju.edu.cn).}

\thanks{$^{\dagger}$School of Mathematics and Statistics, Wuhan University, Wuhan 430072, China (yonglin.li@whu.edu.cn).}

\thanks{$^{\ddagger}$Department of Mathematics, Nanjing University,
Nanjing 210093, China (hjw@nju.edu.cn).}



\begin{abstract}
We study a spherical perfectly matched layer (PML) method for the time-harmonic Maxwell scattering problem with large wavenumber. Using vector spherical harmonics and a modal transformation for the coupled transverse magnetic components, we establish the sharp wavenumber-explicit stability estimate for the truncated PML problem and prove that its inf-sup constant is of order $O(\kappa^{-1})$. We then analyze the first-order N\'ed\'elec edge element method and its continuous interior penalty counterpart. Under the condition that $\kappa^3h^2$ is sufficiently small, the preasymptotic errors in the energy norm and the $\kappa$-scaled $\bm L^2$ norm are of orders $O(\kappa h+\kappa^3h^2)$ and $O(\kappa^2h^2+\kappa^3h^2)$, respectively. Numerical experiments confirm the theoretical results.
%
\end{abstract}

\maketitle

\section{Introduction}

{
In this paper, we consider the following time-harmonic Maxwell scattering
problem:
\begin{subequations}\label{Maxwell equations}
\begin{align}
&\mathbf{curl}\,\mathbf{curl}\,\bm{E} -\kappa^2 \bm{E}=\bm{f}\quad \text{in}\  \mathbb{R}^3 \label{Maxwell equation 1}\\
&\lim_{r\to \infty}\int_{\partial B_r}\left|\mathbf{curl}\,\bm{E}\times \hat{\bm{x}}-\mathbf i \kappa\bm E\right|^2=0.
\label{Maxwell equation 2}
\end{align}
\end{subequations}
Here $r=|\bm x|$, $\hat{\bm x}=\bm x/r$, and
$\bm f\in\bm L^2(\mathbb R^3)$ is supported in
$\Omega:=B_R$, the ball of radius $R$ centered at the origin. Throughout, we
assume that the wavenumber $\kappa\geq1$.

To solve the Maxwell scattering problem \eqref{Maxwell equations}
numerically, we truncate the unbounded exterior domain by a spherical
perfectly matched layer (PML) and discretize the truncated problem using the
first-order N\'ed\'elec edge element method (EEM) of the second type and its
continuous interior penalty counterpart (CIP-EEM). The aim of this paper is
to establish sharp wavenumber-explicit stability estimates for the truncated
Maxwell PML problem and to derive preasymptotic error estimates for the EEM
and CIP-EEM.

The PML method was firstly introduced by B\'erenger in \cite{MR1294924}. It
can be realized by complex coordinate stretching; see, e.g.,
\cite{collino1998optimizing,TeixeiraChew1998}. For the time-harmonic Maxwell
equations, Bao and Wu \cite{MR2192334} firstly proved well-posedness and exponential
convergence of a three-dimensional spherical PML. A key ingredient in their
analysis is the use of uniform asymptotic expansions for the associated
Bessel and Hankel functions. {Bramble and Pasciak developed
variational frameworks for truncated spherical and Cartesian Maxwell PML
problems and proved their stability and exponential convergence
\cite{MR2291829,MR2926536}. For the spherical PML, they further established
discrete stability and quasi-optimal convergence of a N\'ed\'elec edge element
discretization \cite{MR2353940}. We also refer to
\cite{MR1621305,MR1862449,MR2048399,MR2651138,MR3023309,MR3899165}
for the stability and convergence of PML truncations for the Helmholtz
equation, and to \cite{MR2684356,MR3522967} for the corresponding results
in elastic wave scattering. Together, these analyses show that the PML
truncation error decays exponentially as either the absorption or the layer
thickness increases.}

{By contrast, the existing wavenumber-explicit PML analyses mainly
concern the scalar Helmholtz equation.} {Chen and Xiang \cite{MR3085122}
established an inf-sup lower bound proportional to $O(\kappa^{-3/2})$ for a
two-dimensional uniaxial PML (UPML) problem. Li and Wu
\cite{MR3899165} subsequently improved this lower bound to the optimal order
$O(\kappa^{-1})$ for circular and spherical PML truncations using Fourier
series and spherical harmonic expansions, respectively.} More recently,
Wang and Zheng \cite{WangZheng2026} obtained optimal
$O(\kappa^{-1})$ stability for a two-dimensional obstacle problem truncated
by a piecewise linear UPML and established the regularity estimates needed
for its discretization. {An analogous estimate has not yet been
established for Maxwell PML problems. Establishing it is the first objective
of this paper.}

{For finite element approximations of high-frequency wave problems,
the FEM error is not determined solely by the best approximation and may contain
an additional pollution error. For example, Chaumont-Frelet, Galkowski, and
Spence \cite[Theorem~1.3]
{chaumontfrelet2025sharperrorboundsedgeelement} proved for higher-order
N\'ed\'elec $h$-FEM of both types that, for
$\kappa$-oscillatory data and sufficiently small
$C_{\mathrm{sol}}(\kappa h)^{2p}$,
\[
\frac{\|\bm E-\bm E_h\|_{\bm H_\kappa(\mathbf{curl},\Omega)}}
{\|\bm E\|_{\bm H_\kappa(\mathbf{curl},\Omega)}}
\leq C\bigl(1+(\kappa h)^pC_{\mathrm{sol}}\bigr)(\kappa h)^p.
\]
Here
$\|\bm v\|_{\bm H_\kappa(\mathbf{curl},\Omega)}^2
:=\kappa^{-2}\|\mathbf{curl}\,\bm v\|_{\bm L^2(\Omega)}^2
+\|\bm v\|_{\bm L^2(\Omega)}^2$, and $C_{\mathrm{sol}}$ is the
$\bm L^2(\Omega)$-operator norm of the Maxwell solution operator. The
terms $(\kappa h)^p$ and $C_{\mathrm{sol}}(\kappa h)^{2p}$ represent the
best-approximation and pollution errors, respectively. This is a
preasymptotic error estimate: under the weaker condition
$C_{\mathrm{sol}}(\kappa h)^{2p}\ll1$, the pollution term remains visible in
the error bound. The asymptotic error estimate is obtained under the stronger
condition $C_{\mathrm{sol}}(\kappa h)^p\lesssim1$, for which the pollution term
is controlled by the best-approximation term and the interpolation order is
recovered.
To formulate verifiable mesh conditions under which the mesh resolves the
wave, we require a wavenumber-explicit $\bm L^2$-stability estimate for the
continuous Maxwell problem. If
$C_{\mathrm{sol}}\lesssim\kappa$, sufficient conditions for the
preasymptotic and the asymptotic error estimates are
$\kappa^{2p+1}h^{2p}\ll1$ and $\kappa^{p+1}h^p\lesssim1$, respectively.
Galkowski and Spence \cite{GalkowskiSpence2025} established sharp
preasymptotic estimates for the standard Helmholtz $h$-FEM. For Maxwell
equations on curved domains, Lu and Wu \cite{lu2025EEMMaxwell} established
energy-norm and $\kappa$-scaled $\bm L^2$ estimates for higher-order
second-type EEM.
Chaumont-Frelet, Galkowski, and Spence
\cite{chaumontfrelet2025sharperrorboundsedgeelement} subsequently obtained
energy-norm estimates for higher-order EEM of both types. The abstract
Galerkin frameworks in \cite{GalkowskiSpence2025,
chaumontfrelet2025sharperrorboundsedgeelement} are broadly applicable and
include the standard EEM for PML-truncated Maxwell problems. Their estimates,
however, are stated in terms of $C_{\mathrm{sol}}$, whose explicit
dependence on $\kappa$ is left unspecified.
Moreover, the standard Galerkin methods analyzed above quantify the pollution error
but do not provide a mechanism for actively reducing it.

The CIP method was introduced by Douglas and
Dupont \cite{MR440955}. It retains the conforming finite element space and
adds penalties on interelement flux jumps, namely, the jumps of the normal
derivative for the Helmholtz equation and of the tangential trace of the
discrete curl for the Maxwell equations. Zhou and Wu \cite{ZhouWu2023} proved by dispersion analysis that suitable
penalty parameters reduce the leading phase error from
$O((\kappa h)^{2p})$ to $O((\kappa h)^{2p+2})$. Numerical experiments in
\cite{MR3232452,MR3324975,MR3899165,MR4921683,MR4047031,
lu2024preasymptoticerrorestimateseem} further demonstrate that tuning the
penalty parameters can substantially reduce the pollution error.

For the Helmholtz equation, Du and Wu
\cite{MR3324975} established the linear and higher-order preasymptotic
estimates, respectively, while Li and Wu \cite{MR3899165,MR4921683}
incorporated radial PML truncation. For the Maxwell equations, Lu, Wu, and Xu
\cite{MR4047031} proved absolute stability and a first-order preasymptotic
error estimate for CIP-EEM with impedance boundary conditions. Lu and Wu
\cite{lu2024preasymptoticerrorestimateseem} subsequently obtained the unified
$O(\kappa h+\kappa^3h^2)$ estimate for first-order second-type EEM and
CIP-EEM. Since the CIP-EEM form contains mesh-dependent interior-face
penalty terms, it is not directly covered by the standard Galerkin frameworks
in \cite{GalkowskiSpence2025,
chaumontfrelet2025sharperrorboundsedgeelement}.}

{Consequently, in this paper, we begin the preasymptotic analysis
of the EEM and CIP-EEM by establishing a wavenumber-explicit stability
estimate for the truncated Maxwell PML problem. Compared with the scalar
Helmholtz PML problem \cite{MR3899165,WangZheng2026}, the analysis for the
Maxwell PML problem faces two main difficulties: the
infinite-dimensional kernel of the curl operator and the coupling
between the radial and tangential coefficients of each TM mode of the
electric field. We overcome the first
difficulty by using a Helmholtz decomposition weighted by the PML
coefficient. The key tool for treating the second difficulty is a vector
spherical harmonic expansion of the stretched electric field. This expansion
separates the angular variables from the radial dependence, but the two radial
coefficients associated with each TM mode remain coupled. Inspired by Hill's
construction of another equivalent kind of vector spherical harmonics \cite{Hill1954}; see also
\cite[Section~16.4]{Arfken2013}, we introduce an orthogonal transformation
of these two coupled coefficients. The transformation recasts the coupled TM system
as scalar radial PML equations of orders $l-1$ and $l+1$, while the uncoupled TE
coefficient is governed by an equation of order $l$. The transformed
equations can then be estimated by adapting the scalar Helmholtz PML analysis
of Li and Wu \cite{MR3899165} and using the uniform asymptotic estimates for
Bessel and Hankel functions in \cite{MR2192334}.}

Under the PML parameter condition in
Assumption~\ref{assumption:PML parameters} of Section~\ref{sec:pml-problem}, we prove
\[
\kappa\|\hat{\bm E}\|_{\bm L^2(\mathcal D)}
+\|\mathbf{curl}\,\hat{\bm E}\|_{\bm L^2(\mathcal D)}
\lesssim\|\bm f\|_{\bm L^2(\mathcal D)}
\]
for the truncated Maxwell PML problem with arbitrary
$\bm f\in\bm L^2(\mathcal D)$ and show that its inf-sup constant is of the optimal order
$O(\kappa^{-1})$, where $\mathcal D$ denotes the entire computational domain, including both the physical domain $\Omega$ and the PML layer. We also derive the piecewise
regularity needed for the finite element analysis. If $\kappa^3h^2$ is
sufficiently small, then, for sufficiently regular data, the first-order
second-type EEM and CIP-EEM attain preasymptotic energy-norm and
$\kappa$-scaled $\bm L^2$-error bounds
\[
\normtriple{\hat{\bm E}-\bm E_h}_{\gamma}
=O(\kappa h+\kappa^3h^2),
\qquad
\kappa\|\hat{\bm E}-\bm E_h\|_{\bm L^2(\mathcal D)}
=O(\kappa^2h^2+\kappa^3h^2).
\]
The factors suppressed in the $O(\cdot)$ notation are the corresponding
piecewise Sobolev norms of the data; for divergence-free data, they can be
replaced by $\|\bm f\|_{\bm L^2(\mathcal D)}$. The
standard EEM is recovered by setting all penalty parameters to zero.

The remainder of the paper is organized as follows. Section~2 introduces the
notation and preliminary results. Section~3 establishes wavenumber-explicit
stability estimates for the original Maxwell scattering problem. Section~4
formulates the spherical PML problem and derives its VSH modal formulation.
Section~5 proves the stability and inf-sup estimates for the truncated PML
problem. Section~6 presents the EEM and CIP-EEM and derives their discrete
stability and preasymptotic error estimates. Auxiliary spherical-coordinate
calculations and the weighted Helmholtz estimate are collected in the
appendices.
}

\section{Preliminaries}
\subsection{Notation}
We use the standard notation for Sobolev spaces, norms, inner products, and
trace spaces; see, e.g., \cite{MR2059447}. For a domain
$G\subset\mathbb R^d$, $(\cdot,\cdot)_G$ denotes the complex $L^2(G)$ inner
product, with the convention
\[
(u,v)_G=\int_G u\overline v\,\mathrm dx,\qquad
(\bm u,\bm v)_G=\int_G \bm u\cdot\overline{\bm v}\,\mathrm dx .
\]
In particular, $(\cdot,\cdot):=(\cdot,\cdot)_{\mathcal D}$ whenever the
truncated PML domain $\mathcal D$ is under consideration.
The corresponding norm is denoted by $\|\cdot\|_{L^2(G)}$ for scalar functions
and by $\|\cdot\|_{\bm L^2(G)}$ for vector fields. We also use the space
\begin{align*}
\bm H^s(G)&:=[H^s(G)]^3,\\
\bm H(\mathbf{curl},G)
&:=\{\bm v\in \bm L^2(G):\mathbf{curl}\,\bm v\in \bm L^2(G)\},\\
\bm H(\mathbf{div},G)
&:=\{\bm v\in \bm L^2(G):\mathbf{div}\,\bm v\in L^2(G)\},\\
\bm H_0(\mathbf{curl},G)
&:=\{\bm v\in \bm H(\mathbf{curl},G):\bm v\times\bm n=\bm 0\ \text{on }\partial G\},\\
\bm H_0(\mathbf{div},G)
&:=\{\bm v\in \bm H(\mathbf{div},G):\bm v\cdot\bm n=0\ \text{on }\partial G\},\\
\bm H(\mathbf{div}^0,G)
&:=\{\bm v\in \bm H(\mathbf{div},G):\mathbf{div}\,\bm v=0\},\\
\bm H_0(\mathbf{div}^0,G)
&:=\bm H_0(\mathbf{div},G)\cap\bm H(\mathbf{div}^0,G),\\
\bm H^s(\mathbf{curl},G)
&:=\{\bm v\in\bm H^s(G):
\mathbf{curl}\,\bm v\in\bm H^s(G)\},\\
\bm H^s(\mathbf{div},G)
&:=\{\bm v\in\bm H^s(G):
\mathbf{div}\,\bm v\in H^s(G)\}.
\end{align*}
where $\bm n$ denotes the unit outward normal to $\partial G$.
For two disjoint subdomains $G_1$ and $G_2$, the piecewise Sobolev space is
defined by
\[
\begin{aligned}
\bm H^s(G_1\cup G_2)
&:=\{\bm v\in \bm L^2(G_1\cup G_2):\bm v|_{G_i}\in \bm H^s(G_i),\ i=1,2\},\\
\|\bm v\|_{\bm H^s(G_1\cup G_2)}^2
&:=\sum_{i=1}^2\|\bm v\|_{\bm H^s(G_i)}^2 .
\end{aligned}
\]

Throughout the paper, $A\lesssim B$ means that $A\leq C B$ with a generic
constant $C$ independent of the wavenumber and discretization parameters. The
notation $A\gtrsim B$ and $A\eqsim B$ are understood analogously.

\subsection{Stability of the original problem}\label{sec:original-stability}
Next, we derive the wavenumber-explicit stability estimate for the original Maxwell scattering problem, 
which will be used in the stability analysis of the truncated PML problem.
We also refer to \cite{CFMSpence2023} for wavenumber-explicit stability
estimates for full-space Maxwell equations with radially monotone
anisotropic coefficients. 
\begin{thm} \label{stability of original problems}
Let $\bm f\in \bm L^2(\Omega)$ with $\operatorname{supp}\bm f\subset\Omega$.
Then the Maxwell scattering problem \eqref{Maxwell equations}
has a unique solution
$\bm E\in \bm H_{\mathrm{loc}}(\mathbf{curl},\mathbb R^3)$ satisfying
\begin{align}\label{stability estimates of the Maxwell problems}
\kappa \|\bm{E}\|_{\bm{L}^2(\Omega)}
+\|\mathbf{curl}\, \bm E\|_{\bm{L}^2(\Omega)}
\lesssim \|\bm f\|_{\bm{L}^2(\Omega)} .
\end{align}
\end{thm}

\begin{proof}
We first introduce the outgoing Helmholtz volume potential. Set
\begin{equation}\label{Helmholtz Green potential associated with f}
\bm{u}(\bm{x}):=\int_{\Omega} g_\kappa(|\bm{x}-\bm{y}|)\bm{f}(\bm{y})\mathrm{d}\bm{y},
\qquad
g_\kappa(r)=\frac{e^{\mathbf i\kappa r}}{4\pi r}.
\end{equation}
The components of $\bm u$ satisfy the scalar Helmholtz scattering problems
\begin{subequations}\label{scalar Helmholtz problem associated with f}
\begin{align}
 -\Delta u_j-\kappa^2u_j&=f_j\quad\text{in }\mathbb R^3,\\
 \lim_{r\to\infty}\int_{\partial B_r}
 \left|\partial_r u_j-\mathbf i\kappa u_j\right|^2&=0,
 \quad j=1,2,3.
\end{align}
\end{subequations}
Applying \cite[Lemma~3.5]{MR2684350} componentwise gives $\bm u\in \bm H^2_{\mathrm{loc}}(\mathbb R^3)$ and, in particular,
 \begin{align}\label{stability estimates of Hemholtz equationsv2}
 \kappa\|\bm u\|_{\bm L^2(\Omega)}+\|\bm u\|_{\bm H^1(\Omega)}+\kappa^{-1}\|\bm u\|_{\bm H^2(\Omega)}\lesssim\|\bm f\|_{\bm L^2(\Omega)}.
 \end{align}
Now we define
\begin{equation}\label{Maxwell solution representation by Helmholtz potential}
\bm E:=\bm u+\kappa^{-2}\nabla\mathbf{div}\,\bm u .
\end{equation}
We shall prove that $\bm{E}\in \bm{H}_{\mathrm{loc}}(\mathbf{curl},\mathbb{R}^3)$ is  the unique solution of \eqref{Maxwell equations}. Indeed, we have
\[
\mathbf{curl}\,\mathbf{curl}\,\bm E-\kappa^2\bm E
=
\mathbf{curl}\,\mathbf{curl}\,\bm u-\kappa^2\bm u-\nabla\mathbf{div}\,\bm u
=
-\Delta\bm u-\kappa^2\bm u
=
\bm f
\]
Hence $\bm E$ satisfies \eqref{Maxwell equation 1}. The Silver--M\"uller radiation condition \eqref{Maxwell equation 2} follows from the radiation property of the
outgoing Maxwell fundamental solution, see, e.g., \cite[Theorems~5.2.1--5.2.2]{Nedelec2001}.

It remains to prove uniqueness. Let
$\bm W\in \bm H_{\mathrm{loc}}(\mathbf{curl},\mathbb R^3)$ be a radiating solution of
the homogeneous Maxwell problem, that is,
\[
\mathbf{curl}\,\mathbf{curl}\,\bm W-\kappa^2\bm W=0\quad \text{in}\ \mathbb{R}^3.
\]
Taking divergence gives $\mathbf{div}\,\bm W=0$ in $\mathbb R^3$.
Thus each Cartesian component of $\bm W$ satisfies the homogeneous Helmholtz equation.
The Silver--M\"uller condition implies the Sommerfeld radiation condition
for the Cartesian components; see \cite[Theorem~6.8]{MR3971246}.
{Rellich's lemma \cite[Lemma~2.12]{MR3971246} then gives
$\bm W=0$ in the exterior of a sufficiently large ball. The unique
continuation principle \cite[Lemma~8.5]{MR3971246} yields
$\bm W=0$ in $\mathbb R^3$}.

Finally, we infer from
\eqref{stability estimates of Hemholtz equationsv2} and
\eqref{Maxwell solution representation by Helmholtz potential} that  $\mathbf{curl}\,\bm E=\mathbf{curl}\,\bm u$ and
\begin{align*}
\|\mathbf{curl}\,\bm E\|_{\bm L^2(\Omega)}
&\leq \|\bm u\|_{\bm H^1(\Omega)}
\lesssim \|\bm f\|_{\bm L^2(\Omega)},\\
\kappa\|\bm E\|_{\bm L^2(\Omega)}
&\leq
\kappa\|\bm u\|_{\bm L^2(\Omega)}
+\kappa^{-1}\|\nabla\mathbf{div}\,\bm u\|_{\bm L^2(\Omega)}
\lesssim
\|\bm f\|_{\bm L^2(\Omega)} .
\end{align*}
Combining the above two estimates proves \eqref{stability estimates of the Maxwell problems}. This completes the proof.
\end{proof}

\begin{coro}\label{each components of E satisfy the scalar Helmholtz equaion with the Sommerfeld radiation condition}
Under the conditions of Theorem \ref{stability of original problems}, assume, in addition, that $\mathbf{div}\,\bm f=0$ in $\mathbb R^3$, and let
$\bm u$ be defined by \eqref{Helmholtz Green potential associated with f}.
Then $\bm u\in \bm H^2_{\mathrm{loc}}(\mathbb{R}^3)$ is the unique radiating solution of the Maxwell scattering problem
\eqref{Maxwell equations}. Moreover,
\begin{align}\label{stability estimates of Hemholtz equations}
\kappa\|\bm u\|_{\bm L^2(\Omega)}
+\|\bm u\|_{\bm H^1(\Omega)}
+\kappa^{-1}\|\bm u\|_{\bm H^2(\Omega)}
\lesssim \|\bm f\|_{\bm L^2(\Omega)} .
\end{align}
\end{coro}
\begin{proof}
The assertion is also consistent with the discussion around \cite[Theorem~6.8]{MR3971246}; we include the proof for completeness.
By the representation \eqref{Maxwell solution representation by Helmholtz potential},
it suffices to prove $\mathbf{div}\,\bm u=0$. We use a density argument. Since $\bm{f}\in \bm H_0(\mathbf{div}^0; \Omega)$ and $[C_0^\infty(\Omega)]^3$ is dense in $\bm H_0(\mathbf{div}; \Omega)$, there exists $\bm f_n\in [C_0^\infty(\Omega)]^3$ such that
\[
\|\bm f-\bm f_n\|_{\bm L^2(\Omega)}\to0,\qquad
\|\mathbf{div}\,\bm f_n\|_{L^2(\Omega)}\to0 .
\]
Set
\[
\bm u_n(\bm x):=\int_{\Omega}g_\kappa(|\bm x-\bm y|)\bm f_n(\bm y)\,\mathrm d\bm y .
\]
Since $\nabla g_\kappa\in \bm L^1_{\mathrm{loc}}(\mathbb R^3)$ and
$\bm f_n\in [C_0^\infty(\Omega)]^3$,  by using the dominated convergence theorem and the fact that
$\nabla_{\bm x}g_\kappa(|\bm x-\bm y|)
=-\nabla_{\bm y}g_\kappa(|\bm x-\bm y|)$, integration by parts yields
\[
\mathbf{div}\,\bm u_n(\bm x)
=\int_{\Omega}g_\kappa(|\bm x-\bm y|)
\mathbf{div}\,\bm f_n(\bm y)\,\mathrm d\bm y .
\]
Applying the estimate \eqref{stability estimates of Hemholtz equationsv2} to
$\bm u-\bm u_n$ and to $\mathbf{div}\,\bm u_n$, respectively, we obtain
\[
\|\mathbf{div}\,\bm u\|_{L^2(\Omega)}
\leq
\|\mathbf{div}(\bm u-\bm u_n)\|_{L^2(\Omega)}
+\|\mathbf{div}\,\bm u_n\|_{L^2(\Omega)}
\lesssim
\|\bm f-\bm f_n\|_{\bm L^2(\Omega)}
+\kappa^{-1}\|\mathbf{div}\,\bm f_n\|_{L^2(\Omega)} .
\]
Letting $n\to\infty$ yields $\mathbf{div}\,\bm u=0$ in $\Omega$, and the unique
continuation principle \cite[Lemma~8.5]{MR3971246},  yields $\mathbf{div}\,\bm u=0$ in $\mathbb R^3$.
This completes the proof.
\end{proof}

\subsection{Vector spherical harmonic expansion}
In this subsection, we introduce the vector spherical harmonics (VSH) used in the Maxwell multipole expansions below; see, e.g., \cite[Section~5.3]{Nedelec2001}. The VSH play an important role in the stability estimates for PML problem.

The VSH are defined by
\begin{align}\label{vector spherical harmonics}
\bm{Y}_{lm}=Y_l^m \hat{\bm{x}},\quad l\geq0,
\qquad
\bm{U}_{lm}=\frac1{\sqrt{l(l+1)}}\nabla_{\mathbb{S}^2}Y_l^m,\quad
\bm{V}_{lm}=\bm{U}_{lm}\times \hat{\bm{x}},\quad l\geq1,
\end{align}
where $Y_l^m$ are the normalized scalar spherical
harmonics satisfying
\begin{equation*}
\Delta_{\mathbb S^2}Y_l^m+l(l+1)Y_l^m=0,\quad l\in \mathbb{N}_0,\ m\in \mathbb{Z}\cap[-l, l],
\end{equation*}
and $\nabla_{\mathbb S^2}$ and $\Delta_{\mathbb S^2}$ denote the
surface gradient and the Laplace--Beltrami operator on the unit sphere,
respectively. By \cite[Theorem~9.11 and Lemma~9.15]{MR2059447}, the
scalar spherical harmonics $\{Y_l^m\}$ and the VSH
defined in \eqref{vector spherical harmonics} form complete orthonormal bases
of $L^2(\mathbb S^2)$ and $\bm L^2(\mathbb S^2)$, respectively.
Consequently, for any $\rho>0$ and $\bm w\in\bm L^2(B_\rho)$, there admits
\begin{align}
\bm{w}(\bm{x})=w_{00}^{(1)}(r)\bm{Y}_{00}(\hat{\bm{x}})+\sum_{l=1}^{+\infty}\sum_{m=-l}^{l}\left(w_{lm}^{(1)}(r)\bm{Y}_{lm}(\hat{\bm{x}})+w_{lm}^{(2)}(r)\bm{U}_{lm}(\hat{\bm{x}})+w_{lm}^{(3)}(r)\bm{V}_{lm}(\hat{\bm{x}})\right). \label{Fourier expansions}
\end{align} 
The orthonormality of the VSH gives
\begin{align}\label{Parseval identity}
\|\bm w\|_{\bm L^2(B_\rho)}^2
=\int_0^\rho r^2|w_{00}^{(1)}(r)|^2\,\mathrm dr
+\sum_{l=1}^{+\infty}\sum_{m=-l}^{l}\sum_{j=1}^3
\int_0^\rho r^2|w_{lm}^{(j)}(r)|^2\,\mathrm dr .
\end{align}
\begin{lemma}\label{lemma:vsh expansion for divergence-free fields}
For $\bm w\in\bm H(\mathbf{div}^0,B_\rho)$, the coefficient $w_{00}^{(1)}$ in
\eqref{Fourier expansions} vanishes.
\end{lemma}
\begin{proof}
Testing $\mathbf{div}\,\bm w=0$ with
$\varphi(\bm x)=\eta(r)Y_0^0(\hat{\bm x})$, $\eta\in C_0^\infty(0,\rho)$, gives
\[
0=(\mathbf{div}\,\bm w,\varphi)_{B_\rho}
=-(\bm w,\nabla\varphi)_{B_\rho}.
\]
Since $Y_0^0$ is constant, $\nabla\varphi=\eta'(r)Y_0^0\hat{\bm x}$. By
\eqref{Fourier expansions} and orthogonality, we obtain
\[
\int_0^\rho r^2w_{00}^{(1)}(r)\eta'(r)\,\mathrm dr=0
\quad\forall\,\eta\in C_0^\infty(0,\rho).
\]
It follows that $r^2w_{00}^{(1)}=C$ a.e. on $(0,\rho)$. Since
\eqref{Parseval identity} implies
$\int_0^\rho r^2|w_{00}^{(1)}(r)|^2\,\mathrm dr<\infty$, necessarily $C=0$, and
hence $w_{00}^{(1)}=0$.
This proves the lemma.
\end{proof}

\subsection{Modal equations for the original problem}\label{sec_modal_original_problem}
In this subsection, assume that
$\bm f\in\bm H(\mathbf{div}^0,\mathbb R^3)$ and
$\operatorname{supp}\bm f\subset\Omega$. By Corollary~\ref{each components of E satisfy the scalar Helmholtz equaion with the Sommerfeld radiation condition},
$\bm E\in\bm H^2_{\mathrm{loc}}(\mathbb R^3)$ and
$\mathbf{div}\,\bm E=0$. Hence
Lemma~\ref{lemma:vsh expansion for divergence-free fields} gives
\begin{equation}\label{VSH expansions for E and f}
\begin{aligned}
\bm E&=\sum_{l=1}^{+\infty}\sum_{m=-l}^{l}
\left(E_{lm}^{(1)}\bm Y_{lm}
+E_{lm}^{(2)}\bm U_{lm}
+E_{lm}^{(3)}\bm V_{lm}\right),\\
\bm f&=\sum_{l=1}^{+\infty}\sum_{m=-l}^{l}
\left(f_{lm}^{(1)}\bm Y_{lm}
+f_{lm}^{(2)}\bm U_{lm}
+f_{lm}^{(3)}\bm V_{lm}\right).
\end{aligned}
\end{equation}
Substitution into \eqref{Maxwell equation 1}, together with
\eqref{VSH curl curl expansion}, yields
\begin{subequations}\label{components' equations 1}
\begin{align}
-\frac{\sqrt{l(l+1)}}{r^2}\frac{\mathrm d (rE_{lm}^{(2)})}{\mathrm d r}-\left(\kappa^2-\frac{l(l+1)}{r^2} \right)E_{lm}^{(1)}&=f_{lm}^{(1)},\label{components' equations 1.1}\\
-\frac1r\frac{\mathrm d^2 (rE_{lm}^{(2)})}{\mathrm d r^2}-\kappa^2E_{lm}^{(2)}+\frac{\sqrt{l(l+1)}}r\frac{\mathrm d E_{lm}^{(1)}}{\mathrm d r}&=f_{lm}^{(2)},\label{components' equations 1.2}
\\
-\frac1{r}\frac{\mathrm d^2 (rE_{lm}^{(3)})}{\mathrm d r^2}-\left(\kappa^2-\frac{l(l+1)}{r^2} \right)E_{lm}^{(3)}&=f_{lm}^{(3)}.\label{components' equations 1.3}
\end{align}
\end{subequations}
Using \eqref{VSH tangential curl expansion} and the orthogonality of
VSH, \eqref{Maxwell equation 2} reduces to
\begin{subequations}\label{components' equations 2}
\begin{align}
&\lim_{r\to \infty} |rE_{lm}^{(1)}|=0,
\label{modal radiation E1}\\
&\lim_{r\to \infty}\left|\frac{\mathrm d (rE_{lm}^{(2)})}{\mathrm d r}-\mathbf{i}\kappa rE_{lm}^{(2)}-\sqrt{l(l+1)}E_{lm}^{(1)}\right|=0,
\label{modal radiation E2}\\
&\lim_{r\to \infty}\left|\frac{\mathrm d (rE_{lm}^{(3)})}{\mathrm d r}-\mathbf{i}\kappa rE_{lm}^{(3)}\right|=0.
\label{modal radiation E3}
\end{align}
\end{subequations}
Applying \eqref{VSH divergence expansion} to
$\mathbf{div}\,\bm E=\mathbf{div}\,\bm f=0$ gives the two divergence relations
\begin{equation}\label{components' div condition}
\begin{aligned}
\frac{\mathrm d}{\mathrm d r}(r^2E_{lm}^{(1)})
=\sqrt{l(l+1)}\,rE_{lm}^{(2)},\quad
\frac{\mathrm d}{\mathrm d r}(r^2f_{lm}^{(1)})
=\sqrt{l(l+1)}\,rf_{lm}^{(2)}.
\end{aligned}
\end{equation}
To simplify the statement, we introduce the radial Helmholtz operator of order $l$ by
\[
\mathrm L_l[v]:=-\frac1{r^2}\frac{\mathrm d}{\mathrm d r}
\left(r^2\frac{\mathrm d v}{\mathrm d r}\right)
-\left(\kappa^2-\frac{l(l+1)}{r^2}\right)v.
\]
\begin{lemma}[Modal problem for the original scattering problem]
\label{lemma:modal problem for the original scattering problem}
Let $\bm f\in\bm H(\mathbf{div}^0,\mathbb R^3)$ satisfy
$\operatorname{supp}\bm f\subset\Omega$, and let $\bm E$ be the radiating solution of \eqref{Maxwell equations}.
Then, for every $l\in \mathbb{N}$ and $m\in\mathbb Z\cap[-l,l]$, the coefficients defined by
\eqref{VSH expansions for E and f}
satisfy the divergence relations \eqref{components' div condition} and
\begin{subequations}\label{modal operator system}
\begin{align}
\mathrm L_l\left[rE_{lm}^{(1)}\right]
&=rf_{lm}^{(1)},\label{Lu=f 1}\\
\mathrm L_l\left[rE_{lm}^{(2)}\right]
&=rf_{lm}^{(2)}
+\frac{2}{\sqrt{l(l+1)}}
\left(\kappa^2rE_{lm}^{(1)}+rf_{lm}^{(1)}\right),
\label{Lu=f 2}\\
\mathrm L_l\left[E_{lm}^{(3)}\right]
&=f_{lm}^{(3)},\label{Lu=f 3}
\end{align}
\end{subequations}
with the modal radiation conditions
\begin{equation}\label{reduced modal radiation conditions}
\lim_{r\to\infty}
\left|
\frac{\mathrm d}{\mathrm d r}(rv)-\mathbf i\kappa rv
\right|=0,
\qquad
v=rE_{lm}^{(1)},\ E_{lm}^{(2)},\ E_{lm}^{(3)}.
\end{equation}
\end{lemma}
\begin{proof}
Equation \eqref{Lu=f 3} coincides with
\eqref{components' equations 1.3}. Combining
\eqref{components' div condition} and
\eqref{components' equations 1.1}, we further obtain
\begin{align*}
\mathrm L_l\left[rE_{lm}^{(1)}\right]
&=
-\frac{\sqrt{l(l+1)}}{r}
\frac{\mathrm d(rE_{lm}^{(2)})}{\mathrm d r}
-\left(\kappa^2-\frac{l(l+1)}{r^2}\right)rE_{lm}^{(1)}
 =rf_{lm}^{(1)}.
\end{align*}
To derive \eqref{Lu=f 2}, we use the identity
\[
\mathrm L_l[(rv)']
=r(\mathrm L_l[v])'+3\mathrm L_l[v]+2\kappa^2v
\]
with $v=rE_{lm}^{(1)}$. The two relations in
\eqref{components' div condition}, together with \eqref{Lu=f 1}, then give
\[
\begin{aligned}
\mathrm L_l\left[rE_{lm}^{(2)}\right]
&=\frac{1}{\sqrt{l(l+1)}}
\mathrm L_l\left[(r^2E_{lm}^{(1)})'\right]\\
&=\frac{1}{\sqrt{l(l+1)}}
\left[r(rf_{lm}^{(1)})'+3rf_{lm}^{(1)}
+2\kappa^2rE_{lm}^{(1)}\right]\\
&=rf_{lm}^{(2)}+\frac{2}{\sqrt{l(l+1)}}
\left(rf_{lm}^{(1)}+\kappa^2rE_{lm}^{(1)}\right).
\end{aligned}
\]

We finally verify \eqref{reduced modal radiation conditions}. Noting that
\eqref{modal radiation E1} implies
$\lim_{r\to\infty}E_{lm}^{(1)}=0$, the cases
$v=E_{lm}^{(2)},\,E_{lm}^{(3)}$ follow immediately from
\eqref{components' equations 2}.
It remains to consider the case $v=rE_{lm}^{(1)}$. Since
$\operatorname{supp}\bm f\subset\Omega$, equation \eqref{Lu=f 1} is
homogeneous for $r>R$. Together with the first relation in
\eqref{components' div condition}, it yields
\begin{align*}
\sqrt{l(l+1)}\left(
\frac{\mathrm d(rE_{lm}^{(2)})}{\mathrm d r}
-\mathbf i\kappa rE_{lm}^{(2)}
-\sqrt{l(l+1)}E_{lm}^{(1)}\right)=-\mathbf i\kappa\left(
\frac{\mathrm d(r^2E_{lm}^{(1)})}{\mathrm d r}
-\mathbf i\kappa r^2E_{lm}^{(1)}\right).
\end{align*}
The left-hand side tends to zero by \eqref{modal radiation E2} which implies the condition for $v=rE_{lm}^{(1)}$ in \eqref{reduced modal radiation conditions}. This completes the proof.
\end{proof}

To solve the modal problem in
Lemma~\ref{lemma:modal problem for the original scattering problem}, we introduce the variable $z=\kappa r$, the homogeneous
equation $\mathrm L_l[v]=0$  then becomes the spherical Bessel equation of order
$l$ 
$$
\mathrm -\frac1{z^2}\frac{\mathrm d}{\mathrm d z}
\left(z^2\frac{\mathrm d v}{\mathrm d z}\right)
-\left(1-\frac{l(l+1)}{z^2}\right)v=0.
$$
which possesses two linearly independent solutions, denoted as $j_l(z)$ and $h_l(z)$, respectively corresponding to the spherical Bessel and spherical Hankel functions of the first  kind of order $l$. Applying the method of variation of parameter, by noting from \cite[Section~9.3.2]{MR2059447} that 
 \[
\lim_{r\to\infty}\left|
\frac{\mathrm d}{\mathrm d r}\bigl(rh_l(\kappa r)\bigr)
-\mathbf i\kappa r h_l(\kappa r)
\right|=0
\]
and from \cite[Eq.~(10.50.1)]{MR2723248} that the Wronskian $\mathscr{W}_r[j_l(\kappa r),h_l(\kappa r)]=\frac{\mathbf i}{\kappa r^2}$, it can be shown that
\begin{align}
rE_{lm}^{(1)}(r)
&=\mathbf i\kappa j_l(\kappa r)
\int_r^{+\infty}h_l(\kappa t)\left(tf_{lm}^{(1)}(t)\right)t^2\mathrm d t+\mathbf i\kappa h_l(\kappa r)
\int_0^r j_l(\kappa t)\left(tf_{lm}^{(1)}(t)\right)t^2\mathrm d t,
\label{representation of E1}\\
E_{lm}^{(3)}(r)&=\mathbf i\kappa j_l(\kappa r)
\int_r^{+\infty}h_l(\kappa t)f_{lm}^{(3)}(t)t^2\mathrm d t+\mathbf i\kappa h_l(\kappa r)
\int_0^r j_l(\kappa t)f_{lm}^{(3)}(t)t^2\mathrm d t.\label{representation of E3}
\end{align}
And we further deduce from  \eqref{components' div condition} that
\begin{align}
&\sqrt{l(l+1)}\,rE_{lm}^{(2)}(r)\notag
\\
&\,\,\,=\mathbf i\kappa\left[rj_l(\kappa r)\right]'
\int_r^{+\infty}h_l(\kappa t)\left(tf_{lm}^{(1)}(t)\right)t^2\mathrm d t
+\mathbf i\kappa\left[rh_l(\kappa r)\right]'
\int_0^r j_l(\kappa t)\left(tf_{lm}^{(1)}(t)\right)t^2\mathrm d t.
\label{representation of E2}
\end{align}

\section{PML formulation and multipole expansions}\label{sec:pml-problem}
In this section, we first formulate the spherical PML truncation of the Maxwell
scattering problem and state its wavenumber-explicit stability and inf-sup
estimates. We then 
derive the modal problems and representation formulas for both the truncated PML formulation. These modal reductions
provide the basis for the stability analysis in Section~\ref{sec:pml-stability}.

\subsection{Spherical PML formulation}
For a fixed $\hat R>R$, let \(L:=\hat R-R\). The spherical PML method surrounds
$\Omega$ by the finite absorbing layer
\[
\hat\Omega:=\{\bm x\in\mathbb R^3:R<|\bm x|<\hat R\}.
\]
This layer is introduced via the radial complex coordinate transform in the
exterior region; see, e.g., \cite{collino1998optimizing,TeixeiraChew1998}. More
precisely, we define the transform by
$\tilde{\bm x}(\bm x)=\tilde r\,\hat{\bm x}$ with
\[
\tilde r:=\int_0^r\alpha(s)\,\mathrm ds=r\beta(r),
\quad \alpha(r)=1+\mathbf i\sigma(r),\quad \beta(r)=1+\mathbf i\delta(r).
\]
Here $\sigma(r)\geq0$ is the PML absorption function. In this paper, we take 
\begin{align}\label{sigma and delta}
\sigma(r)&=
\begin{cases}
0, & 0\leq r\leq R,\\
\sigma_0, & r>R,
\end{cases}
&
\delta(r)&=
\begin{cases}
0, & 0\leq r\leq R,\\
\dfrac{\sigma_0(r-R)}{r}, & r>R,
\end{cases}
\end{align}
where the constant $\sigma_0>0$ is called the PML parameter.
{
\begin{assumption}\label{assumption:PML parameters}
Throughout the remainder of this paper, we assume that
$$R\simeq\hat R\simeq1,\quad \kappa R>1\quad \text{and}\quad
\kappa\sigma_0L\geq
\max\{2\kappa R+\sqrt3\kappa L,10\}.$$
\end{assumption}
}
\begin{remark}
Polynomial and piecewise constant absorbing profiles are commonly used in PML
approximations; see, for example,
\cite{MR2177884,MR2291829,MR2353940,MR2926536,MR3899165}. The above choice of
$\sigma$ is made to simplify the subsequent modal analysis and to remain
consistent with the Helmholtz PML estimates in \cite{MR3899165}.
{In particular, the PML parameter condition used in the Helmholtz
analysis \cite{MR3899165}, stated in
Assumption~\ref{assumption:PML parameters}, is also sufficient for the Maxwell
stability analysis below.} However, the sharp wavenumber-explicit
stability analysis for general continuous absorbing profiles remains a
delicate issue.
\end{remark}

The first-order spherical PML formulation for Maxwell's equations can be found
in \cite[Eqs.~(2.12)--(2.14)]{MR2192334}. In this paper, we use the equivalent second-order
spherical PML formulation and its derivation is recalled in
Appendix~\ref{app:pml-derivation} for the reader's convenience. 
Overall,
the truncated PML problem for Maxwell’s equations can be written as follows
\begin{subequations}\label{stretched Maxwell equation}
\begin{align}
\mathbf{curl}(\mathbb A\mathbf{curl}\,\hat{\bm E})
-\kappa^2\mathbb A^{-1}\hat{\bm E}
&=\bm f,\quad \text{in }\mathcal D:=B_{\hat{R}},\label{stretched Maxwell equation 1}\\
\hat{\bm E}\times\hat{\bm x}
&=\bm0,\quad \text{on }\hat\Gamma:=\partial \mathcal{D}.\label{stretched Maxwell equation 2}
\end{align}
\end{subequations}
{Here the coefficient matrix $\mathbb A$ is given by
\eqref{PML coefficient A by Jacobian} and its properties are summarized in
Lemma~\ref{lemma:PML coefficient matrix properties}.}

The variational formulation of \eqref{stretched Maxwell equation} is to find
$\hat{\bm E}\in\bm H_0(\mathbf{curl},\mathcal D)$ such that
\begin{equation}\label{variational formulation of truncated PML problem}
a(\hat{\bm E},\bm v)=(\bm f,\bm v)
\quad\forall\,\bm v\in\bm H_0(\mathbf{curl},\mathcal D),
\end{equation}
where the sesquilinear  form
\[
a(\bm u,\bm v)
:=(\mathbb A\mathbf{curl}\,\bm u,\mathbf{curl}\,\bm v)
-\kappa^2(\mathbb A^{-1}\bm u,\bm v).
\]
{We equip $\bm H_0(\mathbf{curl},\mathcal D)$ with the energy norm
\begin{equation}\label{energy norm for PML problem}
\normtriple{\bm v}
:=\left(\Re\left\{a(\bm v,\bm v)
+2\kappa^2(\mathbb A^{-1}\bm v,\bm v)\right\}\right)^{1/2}
\eqsim
\left(\|\mathbf{curl}\,\bm v\|_{\bm L^2(\mathcal D)}^2
+\kappa^2\|\bm v\|_{\bm L^2(\mathcal D)}^2\right)^{1/2}.
\end{equation}}

\subsection{Stability and regularity results}
We first state the wavenumber-explicit stability estimate and the sharp inf-sup estimate for the truncated PML problem,  which constitute the two main results of this paper. Their proofs are deferred in
Section~\ref{sec:pml-stability}.

\begin{thm}\label{lemma:L2 stability of PML Maxwell equations}
Let Assumption~\ref{assumption:PML parameters} hold. Let
$\bm f\in\bm L^2(\mathcal D)$, and  $\hat{\bm E}\in \bm{H}_0(\mathbf{curl},\mathcal{D})$ be the solution of the truncated PML
problem \eqref{stretched Maxwell equation}. Then the following estimate holds:
\begin{equation}\label{stability estimate of PML Maxwell equations}
\kappa\|\hat{\bm E}\|_{\bm L^2(\mathcal D)}
+\|\mathbf{curl}\,\hat{\bm E}\|_{\bm L^2(\mathcal D)}
\lesssim
\|\bm f\|_{\bm L^2(\mathcal D)}.
\end{equation}
\end{thm}

\begin{thm}\label{inf-sup condition}
Let Assumption~\ref{assumption:PML parameters} hold. Then the variational
formulation \eqref{variational formulation of truncated PML problem} satisfies
\begin{equation}\label{inf-sup estimate of truncated PML problem}
\inf_{\bm 0\neq \bm u\in \bm H_0(\mathbf{curl},\mathcal D)}
\sup_{\bm 0\neq \bm v\in \bm H_0(\mathbf{curl},\mathcal D)}
\frac{|a(\bm u,\bm v)|}
{\normtriple{\bm u}\normtriple{\bm v}}
\eqsim\kappa^{-1}.
\end{equation}
\end{thm}
The following lemma provide the higher regularity results for the subsequent finite element analysis.
\begin{lemma}\label{lemma:regularity of PML problem}
Let  Assumption~\ref{assumption:PML parameters} hold and $m\in\mathbb N_0$.
Let
$\bm f\in\bm L^2(\mathcal D)$, and
$\bm\omega\in\bm H_0(\mathbf{curl},\mathcal D)$ be the solution of the truncated PML problem
\eqref{stretched Maxwell equation}. Then
\begin{enumerate}
\item If
$\bm f\in\bm H(\mathbf{div},\mathcal D)
\cap\bm H^m(\mathbf{div},\Omega\cup\hat\Omega)$, we have
$\bm\omega\in\bm H^{m+1}(\mathbf{curl},\Omega\cup\hat\Omega)$ and
\begin{equation}\label{Hm+1estimateofmaxwell}
\kappa\|\bm\omega\|_{\bm H^{m+1}(\Omega\cup\hat\Omega)}
+\|\mathbf{curl}\,\bm\omega\|_{\bm H^{m+1}(\Omega\cup\hat\Omega)}
\lesssim \kappa^{m+1}C_{m+1,\bm f}.
\end{equation}
\item If, in addition, $\mathbf{div}\,\bm f=0$, we have $\bm\omega\in\bm H^{m+2}(\Omega\cup\hat\Omega)$ and
\begin{equation}\label{Hm+2estimateofmaxwell}
\|\bm\omega\|_{\bm H^{m+2}(\Omega\cup\hat\Omega)}
\lesssim \kappa^{m+1}C_{m+1,\bm f},
\end{equation}
\end{enumerate}
where
\[
C_{m+1,\bm f}:=
\|\bm f\|_{\bm L^2(\mathcal D)}
+\sum_{j=0}^{m}\kappa^{-j-1}
\|\bm f\|_{\bm H^j(\Omega\cup\hat\Omega)}
+\kappa^{-m-2}|\mathbf{div}\,\bm f|_{H^m(\Omega\cup\hat\Omega)}.
\]
\end{lemma}
\begin{proof}
Since $\mathbb A$ and
$\mathbb A^{-1}$ are piecewise smooth and uniformly coercive by Lemma~\ref{lemma:PML coefficient matrix properties},  the
regularity estimate in \cite[Theorem~2.10]{melenk2025}, together with
Theorem~\ref{lemma:L2 stability of PML Maxwell equations}, gives
\eqref{Hm+1estimateofmaxwell}. 
If $\mathbf{div}\,\bm f=0$, then
$\mathbf{div}(\mathbb A^{-1}\bm\omega)=0$. Applying the shifted estimate
\cite[Theorem~2.7]{melenk2025} with
$\bm\nu=\mathbb A^{-1}$ and $\ell=m+1$, we obtain
\[
\|\bm\omega\|_{\bm H^{m+2}(\Omega\cup\hat\Omega)}
\lesssim
\|\mathbf{curl}\,\bm\omega\|_{\bm H^{m+1}(\Omega\cup\hat\Omega)}.
\]
Combining this estimate with \eqref{Hm+1estimateofmaxwell} proves
\eqref{Hm+2estimateofmaxwell}. This completes the proof.
\end{proof}

\subsection{Modal equations for the truncated PML problem}\label{Modal equations for the truncated PML problem}
Next, similar to Section \ref{sec_modal_original_problem}, we derive the modal equations for the truncated PML problem.
For this case, we further allow
$\bm f\in\bm H(\mathbf{div}^0,\mathcal D)$ without requiring
$\operatorname{supp}\bm f\subset\Omega$, thereby considering the more general case where the source may extend into the PML region. This extension is needed for the duality argument in the subsequent finite element analysis.

Let $\hat{\bm E}$ solve
\eqref{stretched Maxwell equation 1}--\eqref{stretched Maxwell equation 2}. Denote $ 
\hat{\bm E}^{\mathrm{str}}:=\mathbb B^{-T}\hat{\bm E}, 
\hat{\bm f}:=\mathbb B\bm f,
$
where $\mathbb B$ and $J=\det(\mathbb B)$ are defined in
\eqref{Jacobian matrix of complex stretching}. The transformation identities
\eqref{curl div transformation under PML stretching} yield
\begin{equation}\label{stretched Maxwell equation for transformed field}
\widetilde{\mathbf{curl}}\,
\widetilde{\mathbf{curl}}\,\hat{\bm E}^{\mathrm{str}}
-\kappa^2\hat{\bm E}^{\mathrm{str}}
=J^{-1}\hat{\bm f}
\quad\text{in }\mathcal D.
\end{equation}
Applying $\widetilde{\mathbf{div}}$ to this equation gives
\[
-\kappa^2\widetilde{\mathbf{div}}\,
\hat{\bm E}^{\mathrm{str}}
=\widetilde{\mathbf{div}}(J^{-1}\hat{\bm f})
=J^{-1}\mathbf{div}\,\bm f=0.
\]

\begin{lemma}\label{lemma:zero mode under stretched divergence}
If $\bm w\in\bm L^2(\mathcal D)$ satisfies
$\widetilde{\mathbf{div}}\,\bm w=0$, then its coefficient
$w_{00}^{(1)}$ defined as in \eqref{Fourier expansions} vanishes.
\end{lemma}

\begin{proof}
Set $\bm q:=J\mathbb B^{-1}\bm w\in \bm{L}^2(\mathcal{D})$. Using \eqref{Jacobian matrix of complex stretching}, we have  $
J\mathbb B^{-1}
=H\operatorname{diag}\{\beta^2,\alpha\beta,\alpha\beta\}H^T$, and thus  $q_{00}^{(1)}=\beta^2w_{00}^{(1)}$. Moreover, 
\eqref{curl div transformation under PML stretching} gives
$\mathbf{div}\,\bm q=J\widetilde{\mathbf{div}}\,\bm w=0$.
Lemma~\ref{lemma:vsh expansion for divergence-free fields} (with $\rho=\hat{R}$), therefore yields $q_{00}^{(1)}=0$. Since $\beta\neq 0$, the conclusion follows.
\end{proof}

Let $\hat E_{lm}^{(j)}$ denote the VSH coefficients of $\hat{\bm E}^{\mathrm{str}}$.  Lemma~\ref{lemma:zero mode under stretched divergence} gives
$
\hat E_{00}^{(1)}=0
$ and thus
\begin{equation}\label{VSH expansions for hatE and hatf}
\begin{aligned}
\hat{\bm E}=\mathbb B^T\hat{\bm E}^{\mathrm{str}}
&=\sum_{l=1}^{+\infty}\sum_{m=-l}^{l}
{\left(\alpha
\hat E_{lm}^{(1)}\bm Y_{lm}
+\beta\hat E_{lm}^{(2)}\bm U_{lm}
+\beta\hat E_{lm}^{(3)}\bm V_{lm}
\right)},\\
\hat{\bm f}
&=\sum_{l=1}^{+\infty}\sum_{m=-l}^{l}
\left(
\alpha f_{lm}^{(1)}\bm Y_{lm}
+\beta f_{lm}^{(2)}\bm U_{lm}
+\beta f_{lm}^{(3)}\bm V_{lm}
\right).
\end{aligned}
\end{equation}
And the stretched modal equations of \eqref{stretched Maxwell equation 1} are
\begin{subequations}\label{pml components' equations 1}
\begin{align}
-\frac{\sqrt{l(l+1)}}{r^2}\frac{\mathrm d (\beta r\hat{E}_{lm}^{(2)})}{\mathrm d r}-\left(\alpha \beta^2\kappa^2-\frac{\alpha l(l+1)}{r^2} \right)\hat{E}_{lm}^{(1)}&=\alpha f_{lm}^{(1)}=:\hat{f}_{lm}^{(1)},\label{pml components' equations 1.1}\\
-\frac\beta r\frac{\mathrm d}{\mathrm d r}\left(\frac1\alpha\frac{\mathrm d (\beta r\hat{E}_{lm}^{(2)})}{\mathrm d r}\right)-\alpha\beta^2\kappa^2\hat{E}_{lm}^{(2)}+\frac{\beta \sqrt{l(l+1)}}r\frac{\mathrm d \hat{E}_{lm}^{(1)}}{\mathrm d r}&=\beta f_{lm}^{(2)}=:\hat{f}_{lm}^{(2)},\label{pml components' equations 1.2}
\\
-\frac\beta r\frac{\mathrm d}{\mathrm d r}\left(\frac1\alpha\frac{\mathrm d (\beta r\hat{E}_{lm}^{(3)})}{\mathrm d r}\right)-\left(\alpha\beta^2\kappa^2-\frac{\alpha l(l+1)}{r^2} \right)\hat{E}_{lm}^{(3)}&=\beta f_{lm}^{(3)}=:\hat{f}_{lm}^{(3)}.\label{pml components' equations 1.3}
\end{align}
\end{subequations}
Define the stretched radial Helmholtz operator of order $l$ by
\begin{equation}\label{stretched radial Helmholtz operator}
\widetilde{\mathrm L}_l[v]
:=-\frac1{r^2}\frac{\mathrm d}{\mathrm d r}
\left(\frac{\beta^2}{\alpha}r^2\frac{\mathrm d v}{\mathrm d r}\right)
-\left(\alpha\beta^2\kappa^2
-\frac{\alpha l(l+1)}{r^2}\right)v.
\end{equation}

\begin{lemma}[Modal problem for the truncated PML problem]
\label{lemma:modal problem for the truncated PML equation}
Let $\bm f\in\bm H(\mathbf{div}^0,\mathcal D)$, and let $\hat{\bm E}$ solve
\eqref{stretched Maxwell equation}. Then, for every $l\in\mathbb N$ and
$m\in\mathbb Z\cap[-l,l]$, the coefficients in
\eqref{VSH expansions for hatE and hatf} satisfy
\begin{equation}\label{components' div condition pml}
\frac{\mathrm d}{\mathrm d r}
\left(\beta^2r^2\hat E_{lm}^{(1)}\right)
=\sqrt{l(l+1)}\,\alpha\beta r\hat E_{lm}^{(2)}.
\end{equation}
For $0<r<\hat R$, the corresponding modal equations are
\begin{subequations}\label{tildeLu=f}
\begin{align}
\widetilde{\mathrm L}_l\left[\tilde r\hat E_{lm}^{(1)}\right]
&=\tilde r\hat f_{lm}^{(1)},
\label{tildeLu=f 1}\\
\widetilde{\mathrm L}_l\left[\tilde r\hat E_{lm}^{(2)}\right]
&=\tilde r\hat f_{lm}^{(2)}
+\frac{2}{\sqrt{l(l+1)}}
\left(\alpha\beta^2\kappa^2\tilde r\hat E_{lm}^{(1)}
+\tilde r\hat f_{lm}^{(1)}\right)
,
\label{tildeLu=f 2}\\
\widetilde{\mathrm L}_l\left[\hat E_{lm}^{(3)}\right]
&=\hat f_{lm}^{(3)}.
\label{tildeLu=f 3}
\end{align}
\end{subequations}
with the corresponding boundary conditions
\begin{subequations}\label{boundary condition:tildeLu=f}
\begin{align}
\lim_{r\to 0} r^{\frac12}\hat{E}_{lm}^{(1)}(r)&=0, \quad\frac{\mathrm d}{\mathrm d r}(\beta^2 r^2\hat{E}_{lm}^{(1)})=0\quad \text{at}\ r=\hat{R}.\label{boundary condition:tildeLu=f 1}
\\
\lim_{r\to 0} r^{\frac12}\hat{E}_{lm}^{(2)}(r)&=0, \, \ \ \quad\quad\quad \hat{E}_{lm}^{(2)}(r)=0\quad \text{at}\ r=\hat{R} .\label{boundary condition:tildeLu=f 2}
\\
\lim_{r\to 0}\hat{E}_{lm}^{(3)}(r)&=0,\, \ \ \quad\quad\quad \hat{E}_{lm}^{(3)}(r)=0\quad \text{at}\ r=\hat{R}.\label{boundary condition:tildeLu=f 3}
\end{align}
\end{subequations}
\end{lemma}

\begin{proof}
Note that
in the variable $\tilde r$, 
the stretched equation \eqref{stretched Maxwell equation for transformed field} has the
same form in spherical coordinates as the original equation \eqref{Maxwell equation 1}. Hence, the stretched  divergence relation \eqref{components' div condition pml} and modal equations \eqref{tildeLu=f} follow from  \eqref{components' div condition} and \eqref{modal operator system}, respectively, under the substitutions
\[
r\to\tilde r,\qquad
\bm f\to J^{-1}\hat{\bm f},
\] 
and the identity $\alpha=\mathrm d\tilde r/\mathrm d r, 
\beta=\tilde r/r$. The calculation is straightforward and thus omitted.

Now we turn to derive the boundary conditions at $r=0$ and $r=\hat R$. Since $\mathbf{div}\,\bm{f}=0$ and $\alpha=\beta=1$  in $\Omega$, $\hat{\bm E}^{\text{str}}$ satisfies the vector Helmholtz equation in $\Omega$. 
From the interior estimate of elliptic PDEs and the Sobolev embedding theorem, we have $\hat{\bm E}^{\text{str}}\in\bm H^2(B_\rho)\hookrightarrow\bm C^0(\overline{B_\rho})$ for some $\rho$, $0<\rho<R$. Therefore, by the definition of vector harmonic expansion and using the Cauchy-Schwarz inequality, we deduce that for any $r$, $ 0<r<\rho$, there holds
\begin{align*}
\left|\hat{E}_{lm}^{(1)}(r)\right|+\left|\hat{E}_{lm}^{(2)}(r)\right|=\left|(\hat{\bm{E}}^{\text{str}}(r,\cdot), \bm{Y}_{lm})_{\mathbb{S}^2}\right|+\left|(\hat{\bm{E}}^{\text{str}}(r,\cdot), \bm{U}_{lm})_{\mathbb{S}^2}\right|\lesssim \|\hat{\bm{E}}^{\text{str}}\|_{\bm{L}^{\infty}(B_\rho)}
\end{align*}
Furthermore, the continuity of $\hat{\bm E}^{\text{str}}$ implies
\[
\lim_{r\to0}\hat E_{lm}^{(3)}(r)
=(\hat{\bm E}(\bm0),\bm V_{lm})_{\mathbb S^2}=\frac{1}{\sqrt{l(l+1)}}\int_{\mathbb{S}^2}\mathbf{div}_{\mathbb{S}^2}\left(\hat{\bm E}(\bm0)\times\hat{\bm{x}}\right)\overline{Y_l^m(\hat{\bm{x}})}\mathrm{d}\hat{\bm{x}}=0.
\]
where we used the identity \(\mathbf{div}_{\mathbb S^2}(\hat{\bm x}\times\bm a)=0\) which holds for any constant vector \(\bm a\).

For the boundary condition at $r=\hat{R}$, we conclude from \eqref{VSH expansions for hatE and hatf} and
\eqref{stretched Maxwell equation 2} that
$\hat E_{lm}^{(2)}(\hat R)=\hat E_{lm}^{(3)}(\hat R)=0$, and
\eqref{components' div condition pml} then gives the condition for
$\hat E_{lm}^{(1)}$. That completes the proof.
\end{proof}

Solving the stretched modal problem
\eqref{components' div condition pml}--\eqref{boundary condition:tildeLu=f}
as in the preceding subsection, we obtain
\begin{align}
\begin{split}
\tilde{r}\hat{E}_{lm}^{(1)}(r)&=\mathbf{i}\kappa j_{l}(\kappa\tilde{r})\int_r^{\hat{R}}h_{l}(\kappa\tilde{t})\left(\tilde{t}\hat{f}_{lm}^{(1)}(t)\right)t^2\mathrm d t\\
 &\,+\mathbf{i}\kappa h_{l}(\kappa\tilde{r})\int_0^{r} j_{l}(\kappa\tilde{t})\left(\tilde{t}\hat{f}_{lm}^{(1)}(t)\right)t^2\mathrm d t+\hat{C}_{lm}^{(1)} j_{l}(\kappa\tilde{r}),
\end{split}\label{expression of PML solution hatE1}\\
\sqrt{l(l+1)}\ \alpha\tilde{r}\hat{E}_{lm}^{(2)}(r)&=\mathbf{i}\kappa\left[\tilde{r}j_{l}(\kappa\tilde{r})\right]'\int_r^{\hat{R}}h_{l}(\kappa\tilde{t})\left(\tilde t \hat{f}_{lm}^{(1)}(t)\right)t^2\mathrm d t\label{expression of PML solution hatE2}\\
&\,+\mathbf{i}\kappa\left[\tilde{r}h_{l}(\kappa\tilde{r})\right]'\int_0^{r} j_{l}(\kappa\tilde{t})\left(\tilde{t}\hat{f}_{lm}^{(1)}(t)\right)t^2\mathrm d t+\hat{C}_{lm}^{(2)} \left[\tilde{r}j_{l}(\kappa\tilde{r})\right]'.\notag
\\
\begin{split}
\hat E_{lm}^{(3)}(r)&=\mathbf{i}\kappa j_{l}(\kappa\tilde{r})\int_r^{\hat{R}}h_{l}(\kappa\tilde{t})\hat{f}_{lm}^{(3)}(t)t^2\mathrm d t\\
 &\,+\mathbf{i}\kappa h_{l}(\kappa\tilde{r})\int_0^{r} j_{l}(\kappa\tilde{t})\hat{f}_{lm}^{(3)}(t)t^2\mathrm d t +\hat{C}_{lm}^{(3)} j_{l}(\kappa\tilde{r}).
\end{split}
\end{align}
Recalling $\tilde{\hat{R}}:=\int_0^{\hat{R}}\alpha(r)\mathrm{d}r=\hat{R}+\bm{i}\sigma_0L$, the symbols $\hat{C}_{lm}^{(j)}$ above are
\begin{align*}
\hat{C}_{lm}^{(1)}=\hat{C}_{lm}^{(2)}&=-\mathbf{i}\kappa\frac{h_l(\kappa\tilde{\hat R})+\kappa\tilde{\hat R}h_l'(\kappa\tilde{\hat R})}{j_l(\kappa\tilde{\hat R})+\kappa\tilde{\hat R}j_l'(\kappa\tilde{\hat R})}
\int_0^{\hat R} j_{l}(\kappa\tilde{t})\left(\tilde{t}\hat{f}_{lm}^{(1)}(t)\right)t^2\mathrm d t,\\
\hat{C}_{lm}^{(3)}&=-\mathbf{i}\kappa\frac{h_{l}(\kappa\tilde{\hat R})}{j_{l}(\kappa\tilde{\hat R})}\int_0^{\hat R} j_{l}(\kappa\tilde{t})\hat{f}_{lm}^{(3)}(t)t^2\mathrm d t.
\end{align*}

\section{Proofs of
Theorems~\ref{lemma:L2 stability of PML Maxwell equations} and
\ref{inf-sup condition}}\label{sec:pml-stability}
{In this section, we give the proofs of
Theorems~\ref{lemma:L2 stability of PML Maxwell equations} and
\ref{inf-sup condition}.}

\subsection{\texorpdfstring{Proof of Theorem~\ref{lemma:L2 stability of PML Maxwell equations}}{Proof of Theorem~\ref{lemma:L2 stability of PML Maxwell equations}}}
{We divide the proof of
Theorem~\ref{lemma:L2 stability of PML Maxwell equations} into seven steps.}
\par\smallskip
\noindent\textit{Step 1. Reduction to divergence-free source.}
For $\bm{f}\in \bm{L}^2(\mathcal{D})$,
 by the weighted Helmholtz decomposition
\cite[Theorem~4.23]{MR3288313}, applied to
$\mathbb A\bm f$ with weight $\mathbb A^{-1}$, there exists a decomposition $\bm f=\bm f_0+\mathbb A^{-1}\nabla p$ with $\bm f_0\in\bm H(\mathbf{div}^0,\mathcal D)$ and
$p\in H_0^1(\mathcal D)$ such that
\begin{align}\label{boundedness of Helmholtz decomposition}
\|\bm f_0\|_{\bm L^2(\mathcal D)}
+\|\nabla p\|_{\bm L^2(\mathcal D)}
\lesssim\|\bm f\|_{\bm{L}^2(\mathcal D)}.
\end{align}
Let $\hat{\bm{E}}_0:=\hat{\bm{E}}+\kappa^{-2}\nabla p$. Since $\nabla H_0^1(\mathcal{D})\subset \bm{H}_0(\mathbf{curl}, \mathcal{D})$, it follows that $\hat{\bm{E}}_0\in \bm{H}_0(\mathbf{curl}, \mathcal{D})$ and solves the truncated PML problems \eqref{stretched Maxwell equation} with source $\bm{f}_0$, that is,
\begin{subequations}\label{Maxwell PML problem divfree}
\begin{align}
\mathbf{curl}(\mathbb A\mathbf{curl}\,\hat{\bm E}_0)
-\kappa^2\mathbb A^{-1}\hat{\bm E}_0
&=\bm f_0,\quad \text{in }\mathcal D\\
\hat{\bm E}_0\times\hat{\bm x}
&=\bm0,\quad \text{on }\hat\Gamma.
\end{align}
\end{subequations}
Once the $L^2$-part of
\eqref{stability estimate of PML Maxwell equations} is established for
\eqref{Maxwell PML problem divfree}, the definition of $\hat{\bm{E}}_0$ and
\eqref{boundedness of Helmholtz decomposition} gives
\begin{align*}
\|\hat{\bm E}\|_{\bm L^2(\mathcal D)}\leq
\|\hat{\bm E}_0\|_{\bm L^2(\mathcal D)}+\kappa^{-2}\|\nabla p\|_{\bm L^2(\mathcal D)}\lesssim \kappa^{-1}\|\bm f\|_{\bm L^2(\mathcal D)}.
\end{align*}
Consequently, it remains to establish the above estimate under the assumption
$\bm f\in\bm H(\mathbf{div}^0,\mathcal D)$.
\par\smallskip
\noindent\textit{Step 2. Reduction to modal estimates.}
As discussed in Section \ref{Modal equations for the truncated PML problem}, if $\bm f\in\bm H(\mathbf{div}^0,\mathcal D)$, then $\hat{\bm{E}}$ can be written as \eqref{VSH expansions for hatE and hatf}. Using the Parseval identity \eqref{Parseval identity}, $\|\hat{\bm E}\|_{\bm L^2(\mathcal D)}\lesssim \kappa^{-1}\|\bm f\|_{\bm L^2(\mathcal D)}$ is equivalent to
\begin{align}\label{k estimate 123}
\sum_{j=1}^3\int_0^{\hat R}
r^2\left|\hat E_{lm}^{(j)}\right|^2\,\mathrm dr
\lesssim
\frac{1}{\kappa^2}
\sum_{j=1}^3\int_0^{\hat R}
r^2\left|\hat f_{lm}^{(j)}\right|^2\,\mathrm dr,\quad l\in\mathbb N, m\in\mathbb Z\cap[-l,l].
\end{align}
\par\smallskip
\noindent\textit{Step 3. Estimate for $\hat{E}_{lm}^{(3)}$.}
The following stability result of Helmholtz PML truncation problem established in  \cite[Theorem 3.1]{MR3899165} plays a  key role for our subsequent arguments.
\begin{lemma}\label{lemma stability of Helmholtz pml equation}
Let Assumption~\ref{assumption:PML parameters} hold. 
Then the Helmholtz truncated PML problem
\begin{subequations}
\begin{align}
-\frac{1}{r^2}\frac{\partial}{\partial r}\left(\frac{\beta^2r^2}{\alpha}\frac{\partial \hat u}{\partial r}\right)-\frac{\alpha}{r^2}\Delta_{\mathbb{S}^2}\hat u-\alpha\beta^2\kappa^2\hat u&=f,\quad \text{in} \ \mathcal{D},\label{PML Helmholtz equations 1}\\
\hat u&=0,\quad \text{on}\ \hat\Gamma,\label{PML Helmholtz equations 2}
\end{align}
\end{subequations}
in the polar coordinates, satisfies the following estimate:
\begin{align}\label{L2-estimate of HelmholtzPML}
\|\hat u\|_{L^2(\mathcal{D})}\leq \kappa^{-1}\|f\|_{L^2(\mathcal{D})}.
\end{align}
\end{lemma}
Take $\hat u=\hat E_{lm}^{(3)}Y_l^m$ in Lemma~\ref{lemma stability of Helmholtz pml equation}. 
By \eqref{tildeLu=f 3}--\eqref{boundary condition:tildeLu=f 3},
$\hat u$ solves
\eqref{PML Helmholtz equations 1}--\eqref{PML Helmholtz equations 2}
with $f=\hat f_{lm}^{(3)}Y_l^m$. Hence,
\[
\int_0^{\hat R}r^2\left|\hat E_{lm}^{(3)}(r)\right|^2\,\mathrm dr
\lesssim
\kappa^{-2}\int_0^{\hat R}r^2
\left|\hat f_{lm}^{(3)}(r)\right|^2\,\mathrm dr.
\]
\par\smallskip
\noindent\textit{Step 4. Reduction of $\hat{E}^{(1)}_{lm}$ and $\hat{E}^{(2)}_{lm}$.}  For $l\in\mathbb N_1, m\in\mathbb Z\cap[-l,l]$, we denote the shifted transform
\begin{equation}\label{eq:shifted-definitions}
\begin{aligned}
 \hat E_{lm}^-&:=
 \sqrt{\frac{l}{2l+1}}\,\hat E_{lm}^{(1)}
 +\sqrt{\frac{l+1}{2l+1}}\,\hat E_{lm}^{(2)},
 &
 \hat E_{lm}^+&:=
 \sqrt{\frac{l+1}{2l+1}}\,\hat E_{lm}^{(1)}
 -\sqrt{\frac{l}{2l+1}}\,\hat E_{lm}^{(2)},\\
 \hat f_{lm}^-&:=
 \sqrt{\frac{l}{2l+1}}\,\hat f_{lm}^{(1)}
 +\sqrt{\frac{l+1}{2l+1}}\,\hat f_{lm}^{(2)},
 &
 \hat f_{lm}^+&:=
 \sqrt{\frac{l+1}{2l+1}}\,\hat f_{lm}^{(1)}
 -\sqrt{\frac{l}{2l+1}}\,\hat f_{lm}^{(2)}.
\end{aligned}
\end{equation}
A direct calculation implies the transform is isometric, that is,
\begin{align} \label{eq:norm-identity}
 |\hat E_{lm}^{(1)}|^2+|\hat E_{lm}^{(2)}|^2
 =|\hat E_{lm}^-|^2+|\hat E_{lm}^+|^2
,\quad
 |\hat f_{lm}^{(1)}|^2+|\hat f_{lm}^{(2)}|^2
 =|\hat f_{lm}^-|^2+|\hat f_{lm}^+|^2.
\end{align}
Hence, we can turn to prove
\begin{align}\label{k estimate 12}
\int_0^{\hat R}
r^2\left(\left|\hat E_{lm}^{-}\right|^2+\left|\hat E_{lm}^{+}\right|^2\right)\,\mathrm dr
\lesssim
\frac{1}{\kappa^2}\int_0^{\hat R}
r^2\left(\left|\hat f_{lm}^{-}\right|^2+\left|\hat f_{lm}^{+}\right|^2\right)\,\mathrm dr.
\end{align}
For $q\geq 0$, we define the Green potential of order $q$ by
\begin{align}\label{eq:radial PML Green potential}
(\mathcal{N}_q g)(r):=\mathbf{i}\kappa j_{q}(\kappa\tilde{r})\int_r^{\hat{R}}h_{q}(\kappa\tilde{t})g(t)t^2\mathrm d t
 +\mathbf{i}\kappa h_{q}(\kappa\tilde{r})\int_0^{r} j_{q}(\kappa\tilde{t})g(t)t^2\mathrm d t
\end{align}
It follows from the Wronskian $\mathscr{W}_r[j_l(\kappa \tilde{r}),h_l(\kappa \tilde{r})]=\frac{\mathbf{i}\alpha}{\beta^2\kappa r^2}$ that $N_q g$ solves $\tilde{L}_q[v]=g$.
\begin{lemma}\label{lem:shifted-equations}
Denote $T_{lm}(r):=\tilde{r}^3J^{-1}\hat{f}_{lm}^{(1)}(r)$. For every $l\in\mathbb N$ and
$m\in\mathbb Z\cap[-l,l]$, we have
\begin{subequations}\label{eq:exact-shifted-representations}
\begin{align}
 \hat E_{lm}^-(r)
 =&(\mathcal{N}_{l-1}\hat{f}_{lm}^-)(r)
 +\sqrt{\frac{l+1}{2l+1}}
 \left[
 \frac{\kappa \hat{C}_{lm}^{(1)}}{\sqrt{l(l+1)}}
 -\frac{\mathbf{i}\kappa T_{lm}(\hat{R})}{\sqrt{l(l+1)}}h_{l-1}(\kappa \tilde{\hat{R}})
 \right]j_{l-1}(\kappa\tilde r),
 \label{eq:exact-shifted-representation-minus}\\
 \hat E_{lm}^+(r)
 =&(\mathcal N_{l+1}\hat f_{lm}^+)(r)
 +\sqrt{\frac{l}{2l+1}}
 \left[
 \frac{\kappa \hat{C}_{lm}^{(1)}}{\sqrt{l(l+1)}}
 +\frac{\mathbf{i}\kappa T_{lm}(\hat{R})}{\sqrt{l(l+1)}}h_{l+1}(\kappa \tilde{\hat{R}})
 \right]j_{l+1}(\kappa\tilde r).
 \label{eq:exact-shifted-representation-plus}
\end{align}
\end{subequations}
Moreover,  $\hat E_{lm}^-$ and $\hat E_{lm}^+$ satisfy, respectively,
\begin{equation}\label{eq:shifted-equations}
 \widetilde{\mathrm L}_{l-1}[\hat E_{lm}^-]=\hat f_{lm}^- ,
 \qquad
 \widetilde{\mathrm L}_{l+1}[\hat E_{lm}^+]=\hat f_{lm}^+,
\end{equation}
with the boundary conditions at $r=0$:
\begin{align}
\lim_{r\to 0^+} r^{\frac12} \hat{E}_{lm}^{-}(r)=\lim_{r\to 0^+} r^{\frac12} \hat{E}_{lm}^{+}(r)=0,
\end{align}
and the coupled boundary conditions at $r=\hat{R}$:
\begin{align}
 \tilde r\left(
 \sqrt l\,\frac{\mathrm d\hat E_{lm}^-}{\mathrm d\tilde r}
 +\sqrt{l+1}\,\frac{\mathrm d\hat E_{lm}^+}{\mathrm d\tilde r}
 \right)
 &=(l-1)\sqrt l\,\hat E_{lm}^-
 -(l+2)\sqrt{l+1}\,\hat E_{lm}^+,
 \label{eq:boundary-one}\\
  \sqrt{l+1}\,\hat E_{lm}^--\sqrt l\,\hat E_{lm}^+&=0.
 \label{eq:boundary-two}
\end{align}
\end{lemma}
\begin{proof}
We first deduce from the divergence relation 
\begin{equation}\label{eq:stretched source divergence relation}
\frac{\mathrm d}{\mathrm d r}
\left(\beta^2r^2J^{-1}\hat f_{lm}^{(1)}\right)
=\sqrt{l(l+1)}\,\alpha\beta rJ^{-1}\hat f_{lm}^{(2)}.
\end{equation}
that $A_{lm}:=\tilde r^{-1}T_{lm}\in H^1(0,\hat R)
\hookrightarrow C([0,\hat R])$. Since $A_{lm}/r\in L^2(0,\hat R)$,
$A_{lm}(0)=0$ and, consequently, $T_{lm}(0)=0$.
Set
\[
a_l:=\sqrt{\frac{l}{2l+1}},\qquad
b_l:=\sqrt{\frac{l+1}{2l+1}}.
\]
For $s_q=j_q$ or $h_q$, we infer from
\cite[Eq.~(10.51.1)]{MR2723248} that
\begin{align}
\left[\tilde r s_q(\kappa\tilde r)\right]'
+\alpha q s_q(\kappa\tilde r)
&=\alpha\kappa\tilde r s_{q-1}(\kappa\tilde r),
\label{eq:spherical function lowering relation}\\
\left[\tilde r s_q(\kappa\tilde r)\right]'
-\alpha(q+1)s_q(\kappa\tilde r)
&=-\alpha\kappa\tilde r s_{q+1}(\kappa\tilde r).
\label{eq:spherical function raising relation}
\end{align}
By combining \eqref{expression of PML solution hatE1}--\eqref{expression of PML solution hatE2} according to
\eqref{eq:shifted-definitions} and using
\eqref{eq:spherical function lowering relation}, we obtain
\begin{equation}\label{eq:shifted field minus before integration}
\begin{aligned}
\hat E_{lm}^-(r)
&=\frac{\mathbf i\kappa^2b_l}{\sqrt{l(l+1)}}
\left[
j_{l-1}(\kappa\tilde r)
\int_r^{\hat R}h_l(\kappa\tilde t)
\left(\tilde t\hat f_{lm}^{(1)}(t)\right)t^2\,\mathrm dt
\right.\\
&\hspace{25mm}\left.
+h_{l-1}(\kappa\tilde r)
\int_0^rj_l(\kappa\tilde t)
\left(\tilde t\hat f_{lm}^{(1)}(t)\right)t^2\,\mathrm dt
\right]
+\frac{\kappa b_l\hat C_{lm}^{(1)}}{\sqrt{l(l+1)}}
j_{l-1}(\kappa\tilde r).
\end{aligned}
\end{equation}
Next, using \eqref{eq:spherical function raising relation} and integration by parts, and  \eqref{eq:shifted-definitions} and
\eqref{eq:stretched source divergence relation} 
 yield, for
$0\leq r_1<r_2\leq\hat R$,
\begin{equation}\label{eq:integration identity minus}
\begin{aligned}
\int_{r_1}^{r_2}s_l(\kappa\tilde t)
\left(\tilde t\hat f_{lm}^{(1)}(t)\right)t^2\,\mathrm dt
&=\frac{\sqrt{l(l+1)}}{\kappa b_l}\int_{r_1}^{r_2}s_{l-1}(\kappa\tilde t)
\hat f_{lm}^-(t)t^2\,\mathrm dt-
\left[T_{lm}(t)s_{l-1}(\kappa\tilde t)\right]\big|_{r_1}^{r_2}.
\end{aligned}
\end{equation}
Hence, by applying \eqref{eq:integration identity minus} to the two integrals in
\eqref{eq:shifted field minus before integration}, with
$s_l=h_l$ on $(r,\hat R)$ and $s_l=j_l$ on $(0,r)$, we get \eqref{eq:exact-shifted-representation-minus} while  \eqref{eq:exact-shifted-representation-plus} follows
analogously and thus we omit the detail here.

Moreover, since $\widetilde{\mathrm L}_q[\mathcal{N}_q g]=g$ and $\widetilde{\mathrm L}_q[j_q(\kappa\tilde r)]=0$,  the representations in
\eqref{eq:exact-shifted-representations} imply
\eqref{eq:shifted-equations}. 
Finally, the conditions \eqref{eq:boundary-one}--\eqref{eq:boundary-two}
follow directly from \eqref{boundary condition:tildeLu=f 1}--\eqref{boundary condition:tildeLu=f 2} and the transformation
\eqref{eq:shifted-definitions}.
This completes the proof.
\end{proof}

\par\smallskip
\noindent{\textit{Step 5. Estimates for $\hat E_{lm}^-$ and $\hat E_{lm}^+$: $l\geq2$.}}
{
For $l\geq2$, set $\nu_-:=l-\frac12$ and $\nu_+:=l+\frac32$.
It follows from \eqref{eq:shifted-equations} that
$u_{lm}^{\pm}:=\tilde r^{1/2}\hat E_{lm}^{\pm}$ satisfy
\begin{equation}\label{eq:transformed shifted equations}
-\frac1r\frac{\mathrm d}{\mathrm d r}
\left(\frac{\beta}{\alpha}r
\frac{\mathrm d u_{lm}^{\pm}}{\mathrm d r}\right)
+\left(\frac{\alpha}{\beta}\frac{\nu_{\pm}^2}{r^2}
-\alpha\beta\kappa^2\right)u_{lm}^{\pm}
=\beta^{-1/2}r^{1/2}\hat f_{lm}^{\pm}.
\end{equation}
Let $
V_0:=\{v\in H^1(0,\hat R):v(0)=0\}$ and 
define the coupled space
\begin{equation}\label{eq:coupled modal space}
\mathcal V_l:=\left\{(v^-,v^+)\in V_0\times V_0:
\sqrt{l+1}\,v^-(\hat R)-\sqrt l\,v^+(\hat R)=0\right\}.
\end{equation}
For $\nu\geq 1$, define the sesquilinear form
\begin{equation}\label{eq:modal sesquilinear form}
a_\nu (u,v):=\int_0^{\hat R} \frac{\beta }{\alpha}r\frac{\mathrm d u}{\mathrm d r}\overline{\frac{\mathrm d v}{\mathrm d r}}+\left(\frac{\alpha}{\beta}\frac{\nu^2}{r}-\alpha\beta\kappa^2r\right) u\overline{v}\ \mathrm d r+\frac32u(\hat R)\overline{v(\hat R)}.
\end{equation}
The variational formulation of \eqref{eq:shifted-equations} with
\eqref{eq:boundary-one}--\eqref{eq:boundary-two} is to find
$(u_{lm}^-,u_{lm}^+)\in\mathcal V_l$ such that
\begin{align}
a_{\nu_-}(u_{lm}^-,v^-)+a_{\nu_+}(u_{lm}^+,v^+)=\int_0^{\hat R}\beta^{-1/2}r^{3/2}
\left(\hat f_{lm}^-\overline{v^-}
+\hat f_{lm}^+\overline{v^+}\right)\,\mathrm dr,
\qquad \forall\,(v^-,v^+)\in\mathcal V_l.
\label{eq:coupled variational problem}
\end{align}
}
The following lemma extends the coercivity estimates of $a_\nu$ in
\cite[Lemma~3.3]{MR3899165}.
\begin{lemma}\label{Coercivity of a}
Assume that $\sigma_0\geq\sqrt{3}$. There exists a constant $\lambda$
independent of $\kappa$ such that
\begin{enumerate}[label=(\roman*)]
\item If $\nu^2\geq \kappa^2R^2+\lambda\kappa\hat R^2$, then  for any $u\in V_0$, we have
\begin{align}\label{estimate of Reav}
\Re a_\nu(u,u)
\gtrsim
\int_0^{\hat R}\left(r\left|\frac{\mathrm d u}{\mathrm d r}\right|^2
+\lambda\kappa r|u|^2\right)\,\mathrm d r.
\end{align}
\item If $1\leq\nu^2<\kappa^2R^2+\lambda\kappa\hat R^2$, then for any $u\in V_0$, we have
\begin{align}\label{estimate of Imav}
-\Im a_\nu(u,u)
\gtrsim
\int_R^{\hat R}\left(r\left|\frac{\mathrm d u}{\mathrm d r}\right|^2
+\lambda\kappa r|u|^2\right)\,\mathrm d r.
\end{align}
\item If $1\leq\nu_-^2<\kappa^2R^2+\lambda\kappa\hat R^2\leq\nu_+^2$, 
then, for any $(u^-,u^+)\in\mathcal V_l$,
\begin{align}\label{eq:Garding inequality in the mixed range}
\kappa\int_0^{\hat R}r
\left(|u^-|^2+|u^+|^2\right)\,\mathrm dr\lesssim
\Re\left\{\left(1+\frac{\mathbf i}{\sigma_0}\right)
\left[a_{\nu_-}(u^-,u^-)+a_{\nu_+}(u^+,u^+)\right]\right\}
+\kappa\int_0^Rr|u^-|^2\,\mathrm dr,
\end{align}
\end{enumerate}
\end{lemma}
\begin{proof} Since $\alpha=1+\mathbf i\sigma$ and
$\beta=1+\mathbf i\delta$, a direct calculation implies
\begin{align}
\Re a_\nu(u,u)&=\int_0^{\hat R} \frac{1+\sigma\delta}{1+\sigma^2} r\left| \frac{\mathrm d u}{\mathrm d r}\right|^2 +\left(\frac{1+\sigma\delta}{1+\delta^2}\frac{\nu^2}{r}-\kappa^2(1-\sigma\delta)r\right) |u|^2\ \mathrm d r+\frac32|u(\hat R)|^2,\label{Re a =}\\
-\Im a_\nu(u,u)&=\int_R^{\hat R} \frac{\sigma-\delta}{1+\sigma^2} r\left| \frac{\mathrm d u}{\mathrm d r}\right|^2 +\left(\frac{\delta-\sigma}{1+\delta^2}\frac{\nu^2}{r}\ +\kappa^2(\sigma+\delta)r\right) |u|^2\ \mathrm d r.\label{-Im a=}
\end{align}
From \cite[(3.8)]{MR3899165}, we have
\[
\frac{1+\sigma\delta}{1+\delta^2}\frac{\nu^2}{r}
-\kappa^2(1-\sigma\delta)r
\geq \lambda\kappa r,\qquad r\in[0,\hat R],
\]
and hence, (i) holds. Next we prove (ii).  From \eqref{sigma and delta} and \eqref{-Im a=}, we get
$$
-\Im a_\nu(u,u)=\int_R^{\hat R} \frac{\sigma R}{1+\sigma^2} \left| \frac{\mathrm d u}{\mathrm d r}\right|^2 +\left(-\frac{\sigma R}{1+\delta^2}\frac{\nu^2}{r^2}\ +\kappa^2(\sigma+\delta)r\right) |u|^2\ \mathrm d r.
$$
By using Newton-Leibniz formula and H\"{o}lder  inequality, there holds 
\begin{equation}\label{eq:Newton-Leibniz shift estimate}
|u(r)|^2\leq 2|u(r+c)|^2+2c\int_r^{r+c}\left|\frac{\mathrm d u}{\mathrm ds}\right|^2\mathrm d s,\quad \ R\leq r\leq \hat R-c.
\end{equation}
Taking $c=\hat R-R-\tau/\kappa$ in
\eqref{eq:Newton-Leibniz shift estimate} and integrating with respect to $r$
over $(R,\hat R-c)$, we obtain
\begin{align*}
\int_R^{R+\tau/\kappa}|u|^2\mathrm d r&\leq 2\int_{\hat R-\tau/\kappa}^{\hat R}|u|^2\mathrm d r+\frac{2L\tau}{\kappa}\int_R^{\hat R}\left|\frac{\mathrm d u}{\mathrm dr}\right|^2\mathrm d r\\
&\leq 2\int_{ R+\tau/\kappa}^{\hat R}|u|^2\mathrm d r+\frac{2L\tau}{\kappa}\int_R^{\hat R}\left|\frac{\mathrm d u}{\mathrm dr}\right|^2\mathrm d r.
\end{align*}
Hence,
\begin{align*}
-\Im a_\nu(u,u)=&\frac12\int_R^{\hat R} \frac{\sigma R}{1+\sigma^2} \left| \frac{\mathrm d u}{\mathrm d r}\right|^2 +\frac{\kappa\sigma_0 R}{4L\tau(1+\sigma_0^2)}\left(\frac{2L\tau}{\kappa}\int_R^{\hat R}\left|\frac{\mathrm d u}{\mathrm dr}\right|^2\mathrm d r\right)\\
&+\int_R^{\hat R}\left(-\frac{\sigma R}{1+\delta^2}\frac{\nu^2}{r^2}\ +\kappa^2(\sigma+\delta)r\right) |u|^2\ \mathrm d r\\
\geq& \frac12\int_R^{\hat R} \frac{\sigma R}{1+\sigma^2} \left| \frac{\mathrm d u}{\mathrm d r}\right|^2\mathrm d r+\int_R^{\hat R} \lambda \kappa r|u|^2\mathrm d r\\
&+\int_R^{R+\tau/\kappa}\left(\kappa^2(\sigma+\delta)r-\frac{\sigma R}{1+\delta^2}\frac{\nu^2}{r^2}-\lambda\kappa r +\frac{\kappa\sigma R}{4L\tau(1+\sigma^2)}\right)|u|^2\mathrm d r\\
&+\int_{R+\tau/\kappa}^{\hat R}\left(\kappa^2(\sigma+\delta)r-\frac{\sigma R}{1+\delta^2}\frac{\nu^2}{r^2}-\lambda\kappa r- \frac{\kappa\sigma R}{2L\tau(1+\sigma^2)}\right)|u|^2\mathrm d r
\end{align*}
 Denote by
$$
g(r):=\kappa^2(\sigma+\delta)r-\frac{\sigma R}{1+\delta^2}\frac{\nu^2}{r^2}-\lambda\kappa r,\quad R\leq r\leq \hat R.
$$
If $\lambda<2\kappa\sigma_0$, then $g(r)$ is monotonically increasing. To prove (ii), it is sufficient to ask
$$
g(R)\geq -\frac{\kappa\sigma R}{4L\tau(1+\sigma^2)},\quad g(R+\tau/\kappa)\geq \frac{\kappa\sigma R}{2L\tau(1+\sigma^2)}.
$$
Note that $\frac{\sigma R}{1+\delta^2}\leq \sigma r $ holds for $r\in (R,\hat R)$,   one sufficient condition for the above inequalities is
\begin{align}\label{sufficient condition 1}
\nu^2\leq \kappa^2R^2+\min\left\{ -\frac{\lambda\kappa R^2}{\sigma}+\frac{\kappa R^2}{4L\tau(1+\sigma^2)},\ 3\kappa R\tau-\frac{\lambda\kappa\hat{R}^2}{\sigma}-\frac{\kappa R\hat R}{2L\tau(1+\sigma^2)}\right\}.
\end{align}
Taking $\tau^2=\frac{\hat{R}}{\sigma_0L}$. Note that $\sigma_0^2\geq 3$ and choose 
$\lambda=\frac{3R^2}{16\hat{R}^2L\tau(1+\sigma_0^2)}$, a direct computation shows
$$
3\kappa R\tau-\frac{\lambda\kappa\hat{R}^2}{\sigma}-\frac{\kappa R\hat R}{2L\tau(1+\sigma^2)}\geq  -\frac{\lambda\kappa R^2}{\sigma}+\frac{\kappa R^2}{4L\tau(1+\sigma^2)}\geq \lambda \kappa \hat{R}^2,
$$
and indeed $\lambda<2\kappa\sigma_0$, which indicates \eqref{sufficient condition 1}, and hence, (ii) holds. 
We finally prove (iii). Since
\[
\Re\left\{\left(1+\frac{\mathbf i}{\sigma_0}\right)a_\nu(u,u)\right\}
=\Re a_\nu(u,u)-\sigma_0^{-1}\Im a_\nu(u,u),
\]
it follows from \eqref{Re a =}--\eqref{-Im a=} and
$\delta=\sigma_0(r-R)/r$ that, for every $\nu\geq1$ and $u\in V_0$,
\begin{align}\label{eq:rotated PML lower bound}
&\Re\left\{\left(1+\frac{\mathbf i}{\sigma_0}\right)
a_\nu(u,u)\right\}
-\int_0^Rr
\left|\frac{\mathrm d u}{\mathrm d r}\right|^2
+\left(\frac{\nu^2}{r}-\kappa^2r\right)|u|^2\,\mathrm dr
\notag\\
&\qquad\gtrsim
\int_R^{\hat R}r\left|\frac{\mathrm d u}{\mathrm d r}\right|^2\,\mathrm dr
+\kappa^2\int_R^{\hat R}(r-R)|u|^2\,\mathrm dr
+|u(\hat R)|^2.
\end{align}
We claim that
\begin{align}\label{eq:one-dimensional PML estimate}
\kappa\int_R^{\hat R}r|u|^2\,\mathrm dr
\lesssim\int_R^{\hat R}r\left|\frac{\mathrm d u}{\mathrm d r}\right|^2\,\mathrm dr
+\kappa^2\int_R^{\hat R}(r-R)|u|^2\,\mathrm dr
+|u(\hat R)|^2.
\end{align}
Indeed, if $\kappa L<2$, then \eqref{eq:one-dimensional PML estimate} directly follows from the one-dimensional Friedrichs inequality 
\[
\int_R^{\hat R}|u|^2\,\mathrm dr
\lesssim
L^2\int_R^{\hat R}\left|\frac{\mathrm d u}{\mathrm d r}\right|^2\,\mathrm dr
+L|u(\hat R)|^2.
\]
Suppose now that $\kappa L\geq2$. Then
$R+2\kappa^{-1}\leq\hat R$. Integrating
\eqref{eq:Newton-Leibniz shift estimate} over
$(R,R+\kappa^{-1})$ with $c=\kappa^{-1}$ gives
\begin{align*}
\kappa\int_R^{\hat R}r|u|^2\,\mathrm dr
&=\kappa\int_R^{R+\kappa^{-1}}r|u|^2\,\mathrm dr
+\kappa\int_{R+\kappa^{-1}}^{\hat R}r|u|^2\,\mathrm dr\\
&\lesssim
\int_R^{R+2\kappa^{-1}}r
\left|\frac{\mathrm d u}{\mathrm d r}\right|^2\,\mathrm dr
+\kappa^2\int_{R+\kappa^{-1}}^{\hat R}(r-R)|u|^2\,\mathrm dr,
\end{align*}
where we have used $\kappa^2(r-R)\geq\kappa$ for
$r\geq R+\kappa^{-1}$.
Thus, \eqref{eq:one-dimensional PML estimate} also holds in this case.

Using \eqref{eq:one-dimensional PML estimate} in
\eqref{eq:rotated PML lower bound},  we obtain
\begin{align}\label{eq:rotated PML L2 control}
&\Re\left\{\left(1+\frac{\mathbf i}{\sigma_0}\right)
a_\nu(u,u)\right\}-\int_0^R\left(\frac{\nu^2}{r}-\kappa^2r\right)|u|^2\,\mathrm dr
\gtrsim
\kappa\int_R^{\hat R}r|u|^2\,\mathrm dr.
\end{align}
Let $r\in (0,R)$. For $\nu=\nu_+$, the lower bound of $\nu_+$ in (iii) indicates
\begin{equation}\label{eq:positive shifted physical potential}
\frac{\nu_+^2}{r}-\kappa^2r
\geq\frac{\nu_+^2-\kappa^2R^2}{R}
\geq\frac{\lambda\kappa\hat R^2}{R}
\geq\lambda\kappa r.
\end{equation}
For $\nu=\nu_-$. If $\nu_-\geq\kappa r$, then
$\nu_-^2/r-\kappa^2r\geq0$. If $\nu_-<\kappa r$, then
the identity
$\nu_+^2=\nu_-^2+4\nu_-+4$,  yields
\[
\frac{\nu_-^2}{r}-\kappa^2r
\geq\frac{\nu_+^2-\kappa^2R^2}{R}-\frac{4\nu_-+4}{R}
\geq-\frac{4\nu_-+4}{R}.
\]
Noting that
$\nu_-\geq1$ and $R\simeq1$, we have $(4\nu_-+4)/R\lesssim \nu_-<\kappa r$. Hence, in both cases,
\begin{equation}\label{eq:negative shifted physical potential}
\frac{\nu_-^2}{r}-\kappa^2r\gtrsim-\kappa r.
\end{equation}
Combining \eqref{eq:rotated PML L2 control}-\eqref{eq:negative shifted physical potential}, we conclude (iii) holds. This completes the proof.
\end{proof}
We now prove \eqref{k estimate 12} for $l\geq 2$. The argument is analogous to Cases~1 and~2 in the proof of
\cite[Theorem~3.1]{MR3899165}. For
$\nu_-^2\geq\kappa^2R^2+\lambda\kappa\hat R^2$,
\eqref{eq:coupled variational problem}, \eqref{estimate of Reav}, and the
Cauchy--Schwarz inequality yield
\begin{align*}
\lambda\kappa&\int_0^{\hat{R}} r\left(|u_{lm}^-|^2+|u_{lm}^+|^2\right)\mathrm{d}r\lesssim \Re[a_{\nu_-}(u_{lm}^-,u_{lm}^-)+a_{\nu_+}(u_{lm}^+,u_{lm}^+)]
\\
&
\lesssim \left(\frac1{\lambda\kappa}\int_0^{\hat{R}} r^2\left(|\hat{f}_{lm}^-|^2+|\hat{f}_{lm}^+|^2\right)\mathrm{d}r\right)^{\frac12}\left(\lambda\kappa\int_0^{\hat{R}} r\left(|u_{lm}^-|^2+|u_{lm}^+|^2\right)\mathrm{d}r\right)^{\frac12}
\end{align*}
Therefore, \eqref{k estimate 12} follows from $u_{lm}^{\pm}=\tilde r^{1/2}\hat E_{lm}^{\pm}$. 
For $
\nu_-^2<\kappa^2R^2+\lambda\kappa\hat R^2\leq\nu_+^2$, 
\eqref{eq:coupled variational problem} and
\eqref{eq:Garding inequality in the mixed range}, together with the
Cauchy--Schwarz and Young inequalities, yield
\begin{align*}
&\kappa\int_0^{\hat R}r
\left(|u_{lm}^-|^2+|u_{lm}^+|^2\right)\,\mathrm dr\\
&\quad\lesssim
\left(\int_0^{\hat R}r^2
\left(|\hat f_{lm}^-|^2+|\hat f_{lm}^+|^2\right)\,\mathrm dr\right)^{1/2}
\left(\int_0^{\hat R}r
\left(|u_{lm}^-|^2+|u_{lm}^+|^2\right)\,\mathrm dr\right)^{1/2}
+\kappa\int_0^Rr|u_{lm}^-|^2\,\mathrm dr\\
&\quad\leq
\frac{\kappa}{2}\int_0^{\hat R}r
\left(|u_{lm}^-|^2+|u_{lm}^+|^2\right)\,\mathrm dr
+\frac{C}{\kappa}\int_0^{\hat R}r^2
\left(|\hat f_{lm}^-|^2+|\hat f_{lm}^+|^2\right)\,\mathrm dr
+C\kappa\int_0^Rr|u_{lm}^-|^2\,\mathrm dr.
\end{align*}
Thus,
\begin{align}\label{eq:mixed-mode reduction to physical estimate}
\int_0^{\hat R}r
\left(|u_{lm}^-|^2+|u_{lm}^+|^2\right)\,\mathrm dr\lesssim
\kappa^{-2}\int_0^{\hat R}r^2
\left(|\hat f_{lm}^-|^2+|\hat f_{lm}^+|^2\right)\,\mathrm dr
+\int_0^Rr|u_{lm}^-|^2\,\mathrm dr.
\end{align}
For
$\nu_+^2<\kappa^2R^2+\lambda\kappa\hat R^2$,
by the argument as in the preceding case but using \eqref{estimate of Imav} , we have
\begin{align}\label{eq:mixed-mode reduction to physical estimate 2}
\int_0^{\hat{R}} r\left(|u_{lm}^-|^2+|u_{lm}^+|^2\right)\mathrm{d}r\lesssim \frac1{\lambda^2\kappa^2}\int_0^{\hat{R}} r^2\left(|\hat{f}_{lm}^-|^2+|\hat{f}_{lm}^+|^2\right)\mathrm{d}r+\int_0^{R} r\left(|u_{lm}^-|^2+|u_{lm}^+|^2\right)\mathrm{d}r
\end{align}
Consequently, \eqref{eq:mixed-mode reduction to physical estimate} and
\eqref{eq:mixed-mode reduction to physical estimate 2} reduce the proof to establishing that for $\nu_\pm^2<\kappa^2R^2+\lambda\kappa\hat R^2$,
\begin{align}\label{eq:low-mode-physical-estimate}
\int_0^{R} r|u_{lm}^\pm|^2\mathrm{d}r\lesssim \frac1{\kappa^2}\int_0^{\hat{R}} r^2\left(|\hat{f}_{lm}^-|^2+|\hat{f}_{lm}^+|^2\right)\mathrm{d}r.
\end{align}
{We prove the estimate for $u_{lm}^-$; the other case is analogous and omitted. Using \eqref{eq:exact-shifted-representation-minus}, we obtain
\begin{align}
\int_0^R r|u_{lm}^-|^2\,\mathrm dr
&=\int_0^R r^2|\hat E_{lm}^-|^2\,\mathrm dr
\lesssim I_1+I_2+I_3,
\label{eq:three-term-decomposition-minus}
\end{align}
where
\begin{align*}
I_1
&:=\int_0^R r^2
\left|(\mathcal N_{l-1}\hat f_{lm}^-)(r)\right|^2\,\mathrm dr,\quad I_2
:=\frac{l+1}{2l+1}\int_0^R r^2
\left|
\frac{\kappa\hat C_{lm}^{(1)}}{\sqrt{l(l+1)}}
j_{l-1}(\kappa r)
\right|^2\,\mathrm dr,\\
I_3
&:=\frac{l+1}{2l+1}\int_0^R r^2
\left|
\frac{\kappa T_{lm}(\hat R)}{\sqrt{l(l+1)}}
h_{l-1}(\kappa\tilde{\hat R})j_{l-1}(\kappa r)
\right|^2\,\mathrm dr.
\end{align*}

\noindent\textit{(i) Estimate of $I_1$.}
For $0<r<R$, write
\begin{align*}
(\mathcal N_{l-1}\hat f_{lm}^-)(r)
={}&\underbrace{\mathbf i\kappa j_{l-1}(\kappa r)
\int_r^R h_{l-1}(\kappa t)\hat f_{lm}^-(t)t^2\,\mathrm dt
+\mathbf i\kappa h_{l-1}(\kappa r)
\int_0^r j_{l-1}(\kappa t)\hat f_{lm}^-(t)t^2\,\mathrm dt}
_{=:p_{lm,1}^-(r)}\\
&+\underbrace{\mathbf i\kappa j_{l-1}(\kappa r)
\int_R^{\hat R}h_{l-1}(\kappa\tilde t)
\hat f_{lm}^-(t)t^2\,\mathrm dt}_{=:p_{lm,2}^-(r)}.
\end{align*}
We infer from \cite[Eqs.~(10.47.3) and (10.47.5)]{MR2723248} that
\begin{equation}\label{j-J h-H}
j_q(z)=\left(\frac{\pi}{2z}\right)^{1/2}J_{q+1/2}(z),
\qquad
h_q(z)=\left(\frac{\pi}{2z}\right)^{1/2}H_{q+1/2}^{(1)}(z);
\end{equation}
Therefore,
$r^{1/2}p_{lm,1}^-$ and $r^{1/2}p_{lm,2}^-$ coincide with the first two
terms in the three-dimensional representation
\cite[Eq.~(2.30)]{MR3899165}. The argument in
\cite[Eqs.~(3.17)--(3.18)]{MR3899165}, together with their three-dimensional
counterparts in \cite[Section~3.1.2]{MR3899165}, yield
\begin{equation}\label{eq:green-potential-physical-estimate}
I_1
\lesssim\kappa^{-2}
\int_0^{\hat R}r^2|\hat f_{lm}^-(r)|^2\,\mathrm dr.
\end{equation}
\noindent\textit{(ii) Estimate of $I_2$.} Let $\nu=l+\frac12\geq \frac52$. Denote
$z=\kappa \tilde{\hat{R}}$ and $\hat{z}_\nu=-\bm{i}z/\nu=\frac{\kappa\sigma_0L}{\nu}-\bm{i}\frac{\kappa\hat{R}}{\nu}$, 
and 
$$
\mathcal{A}(\nu,z):=\frac{h_l(z)+zh_l'(z)}{j_l(z)+zj_l'(z)}\overset{\eqref{j-J h-H}}{=}\frac{H_\nu(z)+2zH_\nu'(z)}{J_\nu(z)+2zJ_\nu'(z)}=\frac{H_\nu(z)}{J_{\nu}(z)}\frac{1+2z\frac{H_\nu'(z)}{H_\nu(z)}}{1+2z\frac{J_\nu'(z)}{J_\nu(z)}}:=\frac{H_\nu(z)}{J_{\nu}(z)}\mathcal{B}(\nu,z).
$$
We first refer from the proof of \cite[Theorem 3.1]{MR3899165} that
\begin{align}
\left|\frac{H_q(z)}{J_{q}(z)}\right|\leq e^{-\kappa\sigma_0L},\quad q\geq 1.
\end{align}
To estimate $\mathcal{B}(\nu,z)$, we introduce the modified Bessel functions $I_\nu, K_n$, which satisfy 
$$
J_\nu(z)=e^{\frac{\nu\pi\bm{i}}{2}}I_\nu(-\bm{i}z),\quad H_\nu(z)=-\frac{2\bm{i}}{\pi}e^{-\frac{\nu\pi\bm{i}}{2}}K_\nu(-\bm{i}z),\quad -\frac{\pi}2<arg\, z\leq \pi,
$$
see \cite[Eqs.~(10.27.6)--(10.27.8)]{MR2723248};  The asymptotic expansions of $I_\nu$ and $K_\nu$ in \cite[Lemmas~5.1]{MR2192334} gives
\begin{align}
1+2z\frac{H_\nu'(z)}{H_\nu(z)}&=1-2\bm{i}z\frac{K_\nu'(\nu\hat{z}_\nu)}{K_{\nu}(\nu\hat{z}_\nu)}=1-2\nu(1+\hat{z}_\nu^2)^{\frac12}\frac{1+\eta_4(\nu,\hat{z}_\nu)}{1+\eta_2(\nu,\hat{z}_\nu)},\\
1+2z\frac{J_\nu'(z)}{J_\nu(z)}&=1-2\bm{i}z\frac{I_\nu'(\nu\hat{z}_\nu)}{I_{\nu}(\nu\hat{z}_\nu)}=1+2\nu(1+\hat{z}_\nu^2)^{\frac12}\frac{1+\eta_3(\nu,\hat{z}_\nu)}{1+\eta_1(\nu,\hat{z}_\nu)},
\end{align}
The remainder estimates of $\eta_j, j=1,2,3,4$ in \cite[Lemmas~5.2]{MR2192334} yields
\begin{align*}
|\mathcal{B}(\nu,z)|\leq \frac{\frac{1}{2\nu|1+\hat{z}_\nu^2|^{1/2}}+\frac{1+\widehat{N}_1(\nu,\hat{z}_\nu)}{1-\widehat{M}_1(\nu,\hat{z}_\nu)}}{\frac{1-\widehat{N}_1(\nu,\hat{z}_\nu)}{1+\widehat{M}_1(\nu,\hat{z}_\nu)}-\frac{1}{2\nu|1+\hat{z}_\nu^2|^{1/2}}}:=\mathcal{C}(\nu,z),\quad \text{if}\ \mathcal{C}(\nu,z)>0.
\end{align*}
where the notations $\widehat{M}_1$ and $\widehat{N}_1$ are defined as follows:
\[
\widehat M_1(\nu,\hat{z}_\nu)
:=\exp\left(\frac{2M_1(\hat{z}_\nu)}{\nu}\right)\frac{2M_1(\hat{z}_\nu)}{\nu},
\qquad
\widehat N_1(\nu,\hat{z}_\nu)
:=\frac{N_1(\hat{z}_\nu)}{2\nu}
+\left(1+\frac{N_1(\hat{z}_\nu)}{2\nu}\right)\widehat M_1(\nu,\hat{z}_\nu).
\]
Assumption~\ref{assumption:PML parameters} gives
$\kappa\sigma_0L\geq10$, $\kappa\hat R>1$, and
$\sigma_0L\geq\hat R$. Hence \cite[Lemma~5.2]{MR2192334} yields
\begin{align*}
&M_1(\hat{z}_\nu)\leq \frac1{12}+\frac1{6\sqrt{5}}+\frac{|\Im(\hat{z}_\nu)|}{\Re(\hat{z}_\nu)}\min\left\{\left(\frac{4}{27}\right)^{\frac14},\frac1{\Re(\hat{z}_\nu)}\right\}\leq 0.158+\min\left\{0.621,\frac{\nu}{\kappa\sigma_0L}\right\},\\
&N_1(\hat{z}_\nu)=\min\left\{\left(\frac{4}{27}\right)^{\frac14},\frac1{\Re(\hat{z}_\nu)}\right\}\leq \min\left\{0.621,\frac{\nu}{\kappa(\sigma_0^2L^2+\hat{R}^2)^{1/2}}\right\}.
\end{align*}
If $\nu\geq9/2$, then $M_1(\hat{z}_\nu)\leq0.779$ and hence
$\widehat M_1(\hat{z}_\nu)<0.5$. If $5/2\leq\nu<9/2$, then
\[
\frac{2M_1(\hat{z}_\nu)}{\nu}
\leq\frac{0.316}{\nu}+\frac2{\kappa\sigma_0L}
\leq0.3264,
\qquad
\widehat M_1(\nu,\hat{z}_\nu)\leq0.3264e^{0.3264}<0.5.
\]
Moreover, since $\kappa\sigma_0L\geq 10$, we obtain
$$
N_1(\hat{z}_\nu)<\frac{\nu}{\kappa\sigma_0L}<\frac{\nu}{10},\quad \widehat{N}_1(\nu,\hat{z}_\nu)\leq\frac1{20}+\frac{21}{40}=0.575.
$$
Noting from $\kappa\sigma_0L\geq 10, \kappa\hat{R}>1$ that 
$$
\frac{1}{2\nu|1+\hat{z}_\nu^2|^{1/2}}\leq\frac{1}{2\sqrt{2\kappa\sigma_0L}}<\frac1{2\sqrt{20}}<0.112<\frac{1-\widehat{N}_1(\nu,\hat{z}_\nu)}{1+\widehat{M}_1(\nu,\hat{z}_\nu)}.
$$
Therefore, $\mathcal{C}(\nu,z)>0$ and we obtain
\begin{align*}
|\mathcal{B}(\nu,z)|\leq \mathcal{C}(\nu,z)<\frac{0.112+\frac{1+0.575}{1-0.5}}{\frac{1-0.575}{1+0.5}-0.112}<20,\quad |\mathcal{A}(\nu,z)|\lesssim e^{-\kappa\sigma_0L}.
\end{align*}
Combining this bound with the three-dimensional version of
\cite[Eqs.~(3.21)--(3.24)]{MR3899165} yields
\begin{equation}\label{eq:estimate of correction coefficient}
|\hat C_{lm}^{(1)}|
\lesssim e^{-\kappa\sigma_0L/2}
\left(\int_0^{\hat R}r^2
|\hat f_{lm}^{(1)}(r)|^2\,\mathrm dr\right)^{1/2}.
\end{equation}
Since $|j_q(x)|\leq1$ for $x\geq0$ and $q\geq0$, the definition of $I_2$
and \eqref{eq:estimate of correction coefficient} give
\begin{align}
I_2
&\lesssim \kappa^2e^{-\kappa\sigma_0L}
\int_0^{\hat R}r^2|\hat f_{lm}^{(1)}(r)|^2\,\mathrm dr
\lesssim \kappa^{-2}
\int_0^{\hat R}r^2|\hat f_{lm}^{(1)}(r)|^2\,\mathrm dr,
\label{eq:estimate of I2}
\end{align}
where the last inequality follows from
$\kappa\sigma_0L\geq2\kappa R$ and $R\simeq1$.

\noindent\textit{(iii) Estimate of $I_3$.}
As shown after \eqref{eq:stretched source divergence relation},
$A_{lm}=\tilde r^2J^{-1}\hat f_{lm}^{(1)}(r)
=\tilde r^{-1}T_{lm}(r)\in H^1(0,\hat R)$ and $A_{lm}(0)=0$. 
Using the Newton--Leibniz formula and 
\eqref{eq:stretched source divergence relation} yield
\begin{equation}\label{eq:outer source trace estimate}
\begin{aligned}
\frac{|T_{lm}(\hat R)|}{\sqrt{l(l+1)}}
=\frac{|\tilde{\hat R}|}{\sqrt{l(l+1)}}
\left|\int_0^{\hat R}A_{lm}'(r)\,\mathrm dr\right|
\lesssim
\left(\int_0^{\hat R}r^2
|\hat f_{lm}^{(2)}(r)|^2\,\mathrm dr\right)^{1/2}.
\end{aligned}
\end{equation}
We then recall the following uniform  estimate for the
Hankel functions \cite[Lemma~2.2]{MR2177884}.
\begin{lemma}\label{lemma:uniform exponential decay of Hankel functions}
For any $\nu\in\mathbb R$, $z\in\mathbb C_{++}:=
\{z\in\mathbb C:\Im z\geq0,\ \Re z\geq0\}$, and $\Theta\in\mathbb R$
such that $0<\Theta\leq|z|$, we have
\begin{equation}\label{eq:uniform exponential decay of Hankel functions}
|H_\nu^{(1)}(z)|
\leq
e^{-\Im z\left(1-\Theta^2/|z|^2\right)^{1/2}}
|H_\nu^{(1)}(\Theta)|.
\end{equation}
\end{lemma}
By Assumption~\ref{assumption:PML parameters} and
the choose of $\lambda$ in Lemma~\ref{Coercivity of a}, we have
$\kappa^2R^2+\lambda\kappa\hat R^2
\leq\kappa^2\sigma_0^2L^2$. Hence
\[
0<\nu_-<\kappa\sigma_0L<\kappa|\tilde{\hat R}|.
\]
Taking $\nu=\nu_-$, $z=\kappa \tilde{\hat R}=\kappa\hat R+\mathbf i\kappa\sigma_0L$, and $\Theta=\nu_-$ in
Lemma~\ref{lemma:uniform exponential decay of Hankel functions}, we obtain
\[
|H_{\nu_-}^{(1)}(z)|
\leq
e^{-\kappa\sigma_0L
\left(1-\nu_-^2/(|z|^2)\right)^{1/2}}
|H_{\nu_-}^{(1)}(\nu_-)|.
\]
Since $\nu_-<\kappa\sigma_0L$ and $\sigma_0L\geq\hat R$, we have
\[
\kappa\sigma_0L
\left(1-\frac{\nu_-^2}{|z|^2}\right)^{1/2}
\geq
\frac{\kappa\sigma_0L\hat R}{\sqrt{\hat{R}^2+\sigma_0^2L^2}}
\geq\frac{\kappa\hat R}{\sqrt2}.
\]
Moreover, 
\cite[Eq.~(10.19.8)]{MR2723248} gives
$|H_{\nu_-}^{(1)}(\nu_-)|\lesssim\nu_-^{-1/3}\lesssim1$.
Hence,
we deduce that
\begin{align}\label{ex-decay of h}
|h_{l-1}(\kappa \tilde{\hat R})|=|\sqrt{\pi/(2z)}
H_{\nu_-}^{(1)}(z)|
\lesssim
\kappa^{-1/2}e^{-\kappa\hat R/\sqrt2}
\lesssim\kappa^{-2}.
\end{align}
Here we used $R\simeq\hat R\simeq1$ and $\kappa R>1$. Therefore, it follows from
\eqref{eq:outer source trace estimate} and \eqref{ex-decay of h} that
\begin{equation}\label{eq:estimate of I3}
I_3\lesssim \kappa^{-2}
\int_0^{\hat R}r^2|\hat f_{lm}^{(2)}(r)|^2\,\mathrm dr.
\end{equation}
Hence, \eqref{eq:green-potential-physical-estimate},
\eqref{eq:estimate of I2}, and \eqref{eq:estimate of I3}, together with
\eqref{eq:norm-identity}, yield
\eqref{eq:low-mode-physical-estimate}.  This completes the proof of \eqref{k estimate 12} for $l\geq2$.
}

\par\smallskip
\noindent{\textit{Step 6. Estimates for $\hat E_{lm}^-$ and $\hat E_{lm}^+$: $l=1$.}}
{
Set $z:=\kappa\tilde{\hat R}$.
For $l=1$, the representations
\eqref{eq:exact-shifted-representations} become
\begin{align}
\hat E_{1m}^-
={}&\mathcal N_0\hat f_{1m}^-
+\sqrt{\frac23}\left[
\frac{\kappa\hat C_{1m}^{(1)}}{\sqrt2}
-\frac{\mathbf i\kappa T_{1m}(\hat R)}{\sqrt2}h_0(z)
\right]j_0(\kappa\tilde r),
\label{eq:exact representation l1 minus}\\
\hat E_{1m}^+
={}&\mathcal N_2\hat f_{1m}^+
+\frac1{\sqrt3}\left[
\frac{\kappa\hat C_{1m}^{(1)}}{\sqrt2}
+\frac{\mathbf i\kappa T_{1m}(\hat R)}{\sqrt2}h_2(z)
\right]j_2(\kappa\tilde r).
\label{eq:exact representation l1 plus}
\end{align}
We estimate the three terms in these formulas separately. Recall that
\begin{equation}\label{eq:fixed order spherical Bessel Hankel formulas}
\begin{alignedat}{3}
j_0(w)&=\frac{\sin w}{w},\quad &
j_1(w)&=\frac{\sin w}{w^2}-\frac{\cos w}{w},\quad &
j_2(w)&=\left(\frac3{w^3}-\frac1w\right)\sin w-\frac3{w^2}\cos w,\\
h_0(w)&=-\mathbf i\frac{e^{\mathbf i w}}w,\quad &
h_1(w)&=-e^{\mathbf i w}\left(\frac1w+\frac{\mathbf i}{w^2}\right),\quad &
h_2(w)&=\mathbf i e^{\mathbf i w}
\left(\frac1w+\frac{3\mathbf i}{w^2}-\frac3{w^3}\right).
\end{alignedat}
\end{equation}
Since $\kappa R>1$, we have $\tilde r=r$ whenever $\kappa r\leq1$;
moreover, $|\tilde r|\geq r$. Hence, for $q=0,1,2$,
\begin{align}\label{est of low mode j-h}
|j_q(\kappa\tilde r)|
\lesssim
\begin{cases}
(\kappa r)^q, & \kappa r\leq1,\\
(\kappa r)^{-1}e^{\kappa\Im\tilde r}, & \kappa r\geq1,
\end{cases}
\qquad
|h_q(\kappa\tilde r)|
\lesssim
\begin{cases}
(\kappa r)^{-q-1}, & \kappa r\leq1,\\
(\kappa r)^{-1}e^{-\kappa\Im\tilde r}, & \kappa r\geq1.
\end{cases}
\end{align}
Hence, for $q=0,2$ and $0<s\leq t\leq\hat R$,
\[
\left|\kappa j_q(\kappa\tilde s)h_q(\kappa\tilde t)\right|
\lesssim
\begin{cases}
s^qt^{-q-1},&\kappa t\leq1,\\
\kappa^qs^qt^{-1}e^{-\kappa\Im\tilde t},
&\kappa s\leq1\leq\kappa t,\\
(\kappa st)^{-1}
e^{-\kappa(\Im\tilde t-\Im\tilde s)},&1\leq\kappa s.
\end{cases}
\]
Integration over the three regions above gives
\[
\int_0^{\hat R}\!\int_0^{\hat R}r^2t^2
\left|\kappa j_q(\kappa\tilde r_<)
h_q(\kappa\tilde r_>)\right|^2\,\mathrm dt\,\mathrm dr
\lesssim\kappa^{-2},
\qquad q=0,2,
\]
where $r_<:=\min\{r,t\}$ and $r_>:=\max\{r,t\}$.
The Hilbert--Schmidt inequality therefore yields
\begin{equation}\label{eq:l1 Green potential estimate}
\int_0^{\hat R}r^2|\mathcal N_qg|^2\,\mathrm dr=\int_0^{\hat R}r^2\left|\int_0^{\hat{R}} \kappa j_q(\kappa\tilde r_<)
h_q(\kappa\tilde r_>)g(t)t^2\mathrm{d}t\right|^2\,\mathrm dr
\lesssim\kappa^{-2}\int_0^{\hat R}r^2|g|^2\,\mathrm dr.
\end{equation}
Using \eqref{est of low mode j-h} and $\hat{R}\simeq 1$,  a direct computation implies
\begin{align}
\int_0^{\hat R}r^2|j_q(\kappa\tilde r)|^2\,\mathrm dr
&\lesssim\kappa^{-2}+\kappa^{-3}e^{2\Im z}
\lesssim\kappa^{-3}e^{2\Im z},
\qquad q=0,1,2,
\label{eq:l1 fixed Bessel estimates}\\
|h_q(z)|&\lesssim\kappa^{-1}e^{-\Im z},\qquad q=0,2.
\label{eq:l1 fixed Hankel estimates}
\end{align}
Together with \eqref{eq:outer source trace estimate}, yield
\begin{equation}\label{eq:l1 endpoint estimate}
\kappa^2|T_{1m}(\hat R)|^2
|h_q(z)|^2\int_0^{\hat R}r^2
|j_q(\kappa\tilde r)|^2\,\mathrm dr
\lesssim\kappa^{-2}\int_0^{\hat R}r^2
|\hat f_{1m}^{(2)}|^2\,\mathrm dr,\quad q=0,2.
\end{equation}
It remains to estimate $\hat C_{1m}^{(1)}$. Direct calculation gives
\begin{align*}
h_1(z)+zh_1'(z)
&=e^{\mathbf iz}\left(-\mathbf i+\frac1z+\frac{\mathbf i}{z^2}\right),\\
j_1(z)+zj_1'(z)
&=\frac12e^{-\mathbf iz}
\left(\mathbf i+\frac1z-\frac{\mathbf i}{z^2}\right)
+\frac12e^{\mathbf iz}
\left(-\mathbf i+\frac1z+\frac{\mathbf i}{z^2}\right).
\end{align*}
Since $\Im z=\kappa\sigma_0L\geq10$ and $|z|\geq\Im z$, we have
\begin{align*}
|h_1(z)+zh_1'(z)|
&\leq e^{-\Im z}(1+|z|^{-1}+|z|^{-2}),\\
|j_1(z)+zj_1'(z)|
&\geq\frac12e^{\Im z}(1-|z|^{-1}-|z|^{-2})
-\frac12e^{-\Im z}(1+|z|^{-1}+|z|^{-2})
\gtrsim e^{\Im z}.
\end{align*}
Consequently,
\begin{equation}\label{eq:l1 derivative reflection estimate}
\left|
\frac{h_1(z)+zh_1'(z)}{j_1(z)+zj_1'(z)}
\right|\lesssim e^{-2\Im z}.
\end{equation}
The formula for $\hat C_{1m}^{(1)}$, the Cauchy--Schwarz inequality,
\eqref{eq:l1 fixed Bessel estimates}, and
$|\tilde r|\leq|\tilde{\hat R}|$ give
\begin{align*}
|\hat C_{1m}^{(1)}|^2
\lesssim
\kappa^2 e^{-4\Im z}
\left(\int_0^{\hat R}r^2
|\tilde rj_1(\kappa\tilde r)|^2\,\mathrm dr\right)
\left(\int_0^{\hat R}r^2|\hat f_{1m}^{(1)}|^2\,\mathrm dr\right)
\lesssim\kappa^{-1}e^{-2\Im z}
\int_0^{\hat R}r^2|\hat f_{1m}^{(1)}|^2\,\mathrm dr.
\end{align*}
Consequently, \eqref{eq:l1 fixed Bessel estimates} yields
\begin{equation}\label{eq:l1 correction estimate}
\kappa^2|\hat C_{1m}^{(1)}|^2\int_0^{\hat R}r^2
|j_q(\kappa\tilde r)|^2\,\mathrm dr
\lesssim\kappa^{-2}\int_0^{\hat R}r^2
|\hat f_{1m}^{(1)}|^2\,\mathrm dr,\quad q=0,2.
\end{equation}
Finally, \eqref{eq:exact representation l1 minus}--\eqref{eq:exact representation l1 plus}, together with
\eqref{eq:l1 Green potential estimate}, \eqref{eq:l1 endpoint estimate} and
\eqref{eq:l1 correction estimate}
, give
\[
\int_0^{\hat R}r^2
\left(|\hat E_{1m}^-|^2+|\hat E_{1m}^+|^2\right)\,\mathrm dr
\lesssim\kappa^{-2}
\int_0^{\hat R}r^2
\left(|\hat f_{1m}^-|^2+|\hat f_{1m}^+|^2
+|\hat f_{1m}^{(1)}|^2+|\hat f_{1m}^{(2)}|^2\right)\,\mathrm dr.
\]
By \eqref{eq:norm-identity}, this is precisely \eqref{k estimate 12} for
$l=1$.
}

\par\smallskip
\noindent{\textit{Step 7. Completion of the proof.}}
{
The estimates in Steps 3, 5, and 6, together with
\eqref{eq:norm-identity}, prove \eqref{k estimate 123}. The reductions in
Steps 1 and 2 then give
\begin{equation*}
\|\hat{\bm E}\|_{\bm L^2(\mathcal D)}
\lesssim\kappa^{-1}\|\bm f\|_{\bm L^2(\mathcal D)}.
\end{equation*}
Taking $\bm v=\hat{\bm E}$ in
\eqref{variational formulation of truncated PML problem} and using
Lemma~\ref{lemma:PML coefficient matrix properties}, we obtain
\begin{align*}
\|\mathbf{curl}\,\hat{\bm E}\|_{\bm L^2(\mathcal D)}^2
\lesssim
\kappa^2\|\hat{\bm E}\|_{\bm L^2(\mathcal D)}^2
+\|\bm f\|_{\bm L^2(\mathcal D)}
\|\hat{\bm E}\|_{\bm L^2(\mathcal D)}
\lesssim \|\bm f\|_{\bm L^2(\mathcal D)}^2.
\end{align*}
Combining the two estimates proves
\eqref{stability estimate of PML Maxwell equations} and completes the proof of
Theorem~\ref{lemma:L2 stability of PML Maxwell equations}.
}

\subsection{\texorpdfstring{Proof of Theorem~\ref{inf-sup condition}}{Proof of Theorem~\ref{inf-sup condition}}}
Following the argument of \cite[Theorem~3.5]{MR3899165}, we first prove the
lower bound. For any $\bm0\neq\bm u\in
\bm H_0(\mathbf{curl},\mathcal D)$, let $\bm z\in
\bm H_0(\mathbf{curl},\mathcal D)$ solve the adjoint problem
\begin{equation}\label{adjoint problem for inf-sup condition}
a(\bm v,\bm z)
=2\kappa^2(\mathbb A^{-1}\bm v,\bm u)
\qquad
\forall\,\bm v\in\bm H_0(\mathbf{curl},\mathcal D).
\end{equation}
Since \(\mathbb A^T=\mathbb A\), \eqref{adjoint problem for inf-sup condition} implies that \(\overline{\bm z}\) satisfies \eqref{variational formulation of truncated PML problem} with source \(2\kappa^2\mathbb A^{-1}\overline{\bm u}\). Applying
Theorem~\ref{lemma:L2 stability of PML Maxwell equations}, we have
\begin{equation}\label{adjoint estimate for inf-sup condition}
\normtriple{\bm z}
\lesssim
\kappa^2\|\bm u\|_{\bm L^2(\mathcal D)}
\lesssim
\kappa\normtriple{\bm u}.
\end{equation}
Taking $\bm v=\bm u$ in
\eqref{adjoint problem for inf-sup condition}, we obtain
\[
\Re\,a(\bm u,\bm u+\bm z)
=\Re\left\{a(\bm u,\bm u)
+2\kappa^2(\mathbb A^{-1}\bm u,\bm u)\right\}
=\normtriple{\bm u}^2.
\]
Together with \eqref{adjoint estimate for inf-sup condition}, this yields
\begin{equation}\label{lower bound for inf-sup condition}
\sup_{\bm0\neq\bm v\in\bm H_0(\mathbf{curl},\mathcal D)}
\frac{|a(\bm u,\bm v)|}
{\normtriple{\bm u}\normtriple{\bm v}}
\gtrsim\kappa^{-1}.
\end{equation}
For the upper bound,  let $\bm e_1=(1,0,0)^T, \bm e_2=(0,1,0)^T$ and choose a fixed nonzero
$\chi\in C_0^\infty(\Omega)$. Set
\[
\bm u_\kappa(\bm x):=\chi(\bm x)e^{\mathbf i\kappa x_1}\bm e_2.
\]
Then $\bm u_\kappa\in\bm H_0(\mathbf{curl},\mathcal D)$ and, since
$\mathbb A=\mathbb I$ in $\Omega$,
\begin{equation}\label{plane-wave test for inf-sup upper bound}
a(\bm u_\kappa,\bm v)
=\big(\mathbf{curl}\,\mathbf{curl}\,\bm u_\kappa
-\kappa^2\bm u_\kappa,\bm v\big)_\Omega
\qquad
\forall\,\bm v\in\bm H_0(\mathbf{curl},\mathcal D).
\end{equation}
A direct calculation gives
\[
\mathbf{curl}\,\mathbf{curl}\,\bm u_\kappa
-\kappa^2\bm u_\kappa
=e^{\mathbf i\kappa x_1}
\left[
\nabla\partial_2\chi-(\Delta\chi)\bm e_2
+\mathbf i\kappa
\big((\partial_2\chi)\bm e_1-2(\partial_1\chi)\bm e_2\big)
\right],
\]
and therefore
\[
\|\mathbf{curl}\,\mathbf{curl}\,\bm u_\kappa
-\kappa^2\bm u_\kappa\|_{\bm L^2(\Omega)}
\lesssim 1+\kappa\lesssim\kappa,
\qquad
\normtriple{\bm u_\kappa}
\geq\kappa\|\chi\|_{L^2(\Omega)}.
\]
Since $\chi$ is fixed and nonzero, the Cauchy--Schwarz inequality,
\eqref{plane-wave test for inf-sup upper bound}, and
\eqref{energy norm for PML problem} give
\begin{equation}\label{upper bound for inf-sup condition}
\sup_{\bm0\neq\bm v\in\bm H_0(\mathbf{curl},\mathcal D)}
\frac{|a(\bm u_\kappa,\bm v)|}
{\normtriple{\bm u_\kappa}\normtriple{\bm v}}
\lesssim
\sup_{\bm0\neq\bm v\in\bm H_0(\mathbf{curl},\mathcal D)}
\frac{\kappa\|\bm v\|_{\bm L^2(\Omega)}}
{\kappa\|\chi\|_{L^2(\Omega)}\normtriple{\bm v}}
\lesssim
\frac{1}{\kappa\|\chi\|_{L^2(\Omega)}}
\lesssim\kappa^{-1}.
\end{equation}
Combining \eqref{lower bound for inf-sup condition} and
\eqref{upper bound for inf-sup condition} proves
\eqref{inf-sup estimate of truncated PML problem}. This completes the proof.

\section{Edge element methods and preasymptotic error analysis}\label{sec:discrete-error}
In this section, we analyze the CIP-EEM based on first-order
N\'ed\'elec elements of the second type for the truncated PML problem; the
standard EEM corresponds to vanishing penalty parameters. The key ingredients
are weighted Helmholtz decompositions and elliptic projection estimates, which
lead to the discrete stability and preasymptotic error estimates under the mesh
condition that $\kappa^3h^2$ is sufficiently small.
\subsection{Mesh setting and CIP-EEM}
Let $\{\mathcal T_h\}_{h>0}$ be a family of conforming, shape-regular,
and quasi-uniform curvilinear tetrahedral meshes of $\mathcal D$ in the sense of
\cite[Section~8.1]{MelenkSauter2024Maxwell}; see also
\cite[Section~2]{lu2024preasymptoticerrorestimateseem}. In particular,
$\overline{\mathcal D}
=\bigcup_{K\in\mathcal T_h}\overline K$ and the element map $F_K: \overline{\hat K}\to \overline{K}$ is  $C^1$-diffeomorphic. Here, $\hat K$ is the reference tetrahedron.
 The meshes are further asked to  resolve the interface $\Gamma=\partial\Omega$, namely,
$K^\circ\cap\Gamma=\emptyset$ for every $K\in\mathcal T_h$.
For $K\in\mathcal T_h$ and a face $F\subset\partial K$, set
$h_K:=\operatorname{diam}(K)$, $h_F:=\operatorname{diam}(F)$, and
$h:=\max_{K\in\mathcal T_h}h_K$. The shape regularity and
quasi-uniformity give $h_K\eqsim h_F\eqsim h$.

We additionally require
$\mathcal T_h$ to be a normalizable regular triangulation as specified in
\cite[Assumption~8.1]{MelenkSauter2024Maxwell}. This is the geometric
assumption required for the uniform curvilinear interpolation estimates
used below.
We also set
$
\mathcal T_h^\Omega
:=\{K\in\mathcal T_h:K\subset\overline\Omega\},  
\mathcal F_h^I
:=\{F=\partial K^+\cap\partial K^-:
K^+,K^-\in\mathcal T_h^\Omega\}.
$
And for $F=\partial K^+\cap\partial K^-\in\mathcal F_h^I$, we define the tangential jump
\[
\llbracket\bm v\rrbracket_F
:=\bm v|_{K^+}\times\bm\nu^+
+\bm v|_{K^-}\times\bm\nu^-,
\]
where  $\bm\nu^\pm$ are the unit outward normals to $K^\pm$, respectively. 

Let $V_h$ denote the first-order N\'ed\'elec edge element space of the
second type:
\[
V_h:=\left\{\bm v_h\in\bm H_0(\mathbf{curl},\mathcal D):
(DF_K)^T(\bm v_h|_K\circ F_K)\in\mathcal P_1(\hat K)^3
\quad\forall\,K\in\mathcal T_h\right\}.
\]
To formulate the CIP-EEM, we introduce the space
\[
V:=\left\{\bm v\in\bm H_0(\mathbf{curl},\mathcal D):
\bm v|_K\in\bm H^1(\mathbf{curl},K)
\quad\forall\,K\in\mathcal T_h\right\}.
\]
Then $V_h\subset V$. And define the penalty term  on $V\times V$ by
$$
J(\bm u,\bm v):=\sum_{F\in \mathcal{F}^{I}_h}\gamma_Fh_F\langle \llbracket\mathbf{curl}\,\bm u\rrbracket,\llbracket\mathbf{curl}\,\bm v\rrbracket\rangle_F.
$$
where $\gamma_F, F\in \mathcal{F}_h^I$ are the penalty parameters. Clearly, if $\bm{v}\in \bm{H}^1(\mathbf{curl},\Omega)$, then we have
$$
J(\bm{v},\bm{v}_h)=0,\quad \forall\ \bm{v}_h\in V_h.
$$ 
Define the sesquilinear form
\[
a_\gamma(\bm u,\bm v):=a(\bm u,\bm v)+J(\bm u,\bm v).
\]
Consequently, the CIP-EEM for \eqref{stretched Maxwell equation} reads as follows: find $\bm E_h\in V_h$ such that
\begin{align}\label{agamma}
a_\gamma({\bm E}_h,\bm v_h)=(\bm f, \bm v_h),
\qquad \forall\,\bm v_h\in V_h.
\end{align}
\begin{remark}
\begin{enumerate}[label=(\roman*),leftmargin=*]
\item If $\gamma_F=0$ for all $F\in\mathcal F_h^I$, then the CIP-EEM
\eqref{agamma} reduces to the standard EEM.

\item If $\bm f\in\bm H(\mathbf{div},\mathcal D)$, then the CIP-EEM satisfies
the Galerkin orthogonality
\begin{equation}\label{CIP Galerkin orthogonality}
a_\gamma(\hat{\bm E}-\bm E_h,\bm v_h)=0,
\qquad \forall\,\bm v_h\in V_h.
\end{equation}
Indeed, Lemma~\ref{lemma:regularity of PML problem} with $m=0$ implies
$J(\hat{\bm E},\bm v_h)=0$ for every $\bm v_h\in V_h$, and
\eqref{CIP Galerkin orthogonality} follows by subtracting \eqref{agamma}
from \eqref{variational formulation of truncated PML problem}.
\end{enumerate}
\end{remark}

To measure the numerical error of CIP-EEM, we also introduce the penalized energy norm by
\begin{equation*}
\normtriple{\bm v}_\gamma
:=\Big(
\normtriple{\bm v}^2
+\sum_{F\in\mathcal F_h^I}|\gamma_F|h_F
\|\llbracket\mathbf{curl}\,\bm v\rrbracket\|_{\bm L^2(F)}^2
\Big)^{1/2}.
\end{equation*}
\subsection{Weighted Helmholtz decompositions}
\label{subsection:weighted Helmholtz decompositions}
A basic difficulty in the analysis of Maxwell problems is that the curl
operator has an infinite-dimensional kernel.
 For the PML problem considered
here, the weighted Helmholtz decomposition is a key tool for overcoming this
difficulty. We therefore introduce its continuous and discrete forms for the
subsequent analysis.
\begin{lemma}[$\mathbb M$-weighted Helmholtz decomposition]
\label{lemma:weighted Helmholtz decomposition} Let  $\mathbb M\in L^\infty(\mathcal{D};\mathbb C^{3\times3})$ such that $\mathbb{M}(\bm{x})$ is symmetric for almost all $\bm x\in  \mathcal D$. Furthermore, assume that there exists a constant $\hat{c}$ such that $
\Re\bigl(\mathbb M(\bm x)\bm\xi\cdot\overline{\bm\xi}\bigr)
\geq \hat{c}|\bm\xi|^2
$ for all $\bm\xi\in\mathbb C^3$ and almost all $\bm x\in \mathcal D$.
Then for every
$\bm w\in\bm H_0(\mathbf{curl},\mathcal D)$, there admits a unique decomposition
\begin{equation}\label{weighted Helmholtz decomposition}
\bm w=\Pi^\nabla_{\mathbb{M}}\bm w+\Pi^{\mathbf{curl}}_{\mathbb{M}}\bm w,
\end{equation}
where $\Pi^\nabla_{\mathbb{M}}\bm w\in \nabla H_0^1(\mathcal D)$ satisfies
$
(\mathbb{M}\,\nabla\xi,\Pi^\nabla_{\mathbb{M}}\bm w)=(\mathbb{M}\,\nabla\xi, \bm w)
$ for all $\xi\in H_0^1(D)$.
Moreover,
\begin{equation}\label{properties of weighted curl component}
\mathbf{curl}\,\Pi^{\mathbf{curl}}_{\mathbb M}\bm w
=\mathbf{curl}\,\bm w,
\qquad
\mathbf{div}\bigl(\mathbb M^H
\Pi^{\mathbf{curl}}_{\mathbb M}\bm w\bigr)=0.
\end{equation}
\end{lemma}
\begin{proof}
For the proof of the lemma, we refer to 
\cite[Theorem~4.23]{MR3288313}, applied with weight $A=\mathbb M^{H}$.
\end{proof}
For the Galerkin counterpart of
Lemma~\ref{lemma:weighted Helmholtz decomposition},  we define the quadratic Lagrange finite element space
\begin{equation}\label{quadratic scalar finite element space}
S_h:=\left\{q_h\in H_0^1(\mathcal D):
q_h|_K\circ F_K\in\mathcal P_2(\hat K)
\quad\forall\,K\in\mathcal T_h\right\}.
\end{equation}
Then $\nabla S_h\subset V_h$. For any
$\bm w_h\in V_h$,  there exists a splitting
\begin{equation}\label{discrete weighted gradient projection}
\bm w_h=\Pi^\nabla_{\mathbb{M},h}\bm{w}_h+\Pi^{\mathbf{curl}}_{\mathbb{M},h}\bm{w}_h
\end{equation}
such that $\Pi^\nabla_{\mathbb{M},h}\bm{w}_h\in \nabla S_h$ is determined by
\[
(\mathbb{M}\,\nabla\xi_h,\Pi^\nabla_{\mathbb{M},h}\bm{w}_h)
=(\mathbb{M}\,\nabla\xi_h,\bm w_h),
\qquad \forall\,\xi_h\in S_h.
\]
The well-posedness of the projection operator $\Pi^\nabla_{\mathbb{M},h}$ follows from Lax-Milgram lemma. 
The following lemma gives the continuous and discrete weighted Helmholtz
decompositions of $\bm w_h\in V_h$.
\begin{lemma}\label{lemma:weighted continuous-discrete Helmholtz estimate}
Under the assumption of lemma \ref{lemma:weighted Helmholtz decomposition}, we further suppose that $\mathbb{M}$ is piecewise $C^1$ smooth, precisely, its components are in $C^1(\Omega\cup\hat{\Omega})$. Then for every $\bm w_h\in V_h$, we have
\begin{align}
&\bm w_h=\Pi^\nabla_{\mathbb{M}}\bm{w}_h+\Pi^{\mathbf{curl}}_{\mathbb{M}}\bm{w}_h=\Pi^\nabla_{\mathbb{M},h}\bm{w}_h+\Pi^{\mathbf{curl}}_{\mathbb{M},h}\bm{w}_h,\\
\label{curl-curlh}
&\|\Pi^{\mathbf{curl}}_{\mathbb{M}}\bm w_h
-\Pi^{\mathbf{curl}}_{\mathbb{M},h}\bm w_h\|_{\bm L^2(\mathcal D)}
\lesssim
h\|\mathbf{curl}\,\bm w_h\|_{\bm L^2(\mathcal D)}.
\end{align}
\end{lemma}
\begin{proof}
For the readers' convenience, the proof is given in
Appendix~\ref{appendix:weighted Helmholtz estimate}; see also
\cite[Lemma~7.6]{MR2059447} for the corresponding unweighted argument.
\end{proof}

\subsection{Approximation estimates and elliptic projections}
For interpolation estimates on curved meshes, we refer to
\cite[Theorem~10.2 and Appendix~B]
{chaumontfrelet2025sharperrorboundsedgeelement} and
\cite[Lemma~2.1]{lu2025EEMMaxwell} for N\'ed\'elec elements of the first and
second types, respectively. For the affine case, see also
\cite[Theorem~8.15]{MR2059447}.

Following \cite{lu2025EEMMaxwell}, we denote the mapped second-type N\'ed\'elec linear interpolant by
\[
\pi_h:\bm H(\mathbf{curl},\mathcal D)
\cap\bm H^1(\mathbf{curl},\Omega\cup\hat\Omega)\to\bm V_h.
\]
\begin{lemma}\label{interpolation theorem}
For any
$\bm u\in\bm H(\mathbf{curl},\mathcal D)
\cap\bm H^1(\mathbf{curl},\Omega\cup\hat\Omega)$, there holds
\begin{align}
\|\bm u-\pi_h\bm u\|_{\bm L^2(\mathcal D)}
&\lesssim
h\|\bm u\|_{\bm H^1(\Omega\cup\hat\Omega)}
+h^2\|\mathbf{curl}\,\bm u\|_{\bm H^1(\Omega\cup\hat\Omega)},
\label{low regularity interpolation estimate}
\end{align}
If, in addition, $\bm u\in\bm H^2(\Omega\cup\hat\Omega)$, then
\begin{align}\label{high regularity interpolation estimates}
\|\bm{u}-\pi_h\bm{u}\|_{\bm{L}^2(\mathcal{D})}\lesssim h^2\|\bm{u}\|_{\bm{H}^2(\Omega\cup\hat{\Omega})},\quad 
\|\mathbf{curl}(\bm{u}-\pi_h\bm{u})\|_{\bm{L}^2(\mathcal{D})}\lesssim h\|\mathbf{curl}\, \bm{u}\|_{\bm{H}^1(\Omega\cup\hat{\Omega})}
\end{align}
\end{lemma}
\begin{proof}
The estimate \eqref{low regularity interpolation estimate}
follows from \cite[Theorem~10.2 and Appendix~B]
{chaumontfrelet2025sharperrorboundsedgeelement}, with the first-type
reference interpolant replaced by its second-type counterpart as in the
proof of \cite[Lemma~2.1]{lu2025EEMMaxwell}. The remaining estimates follow
from \cite[Lemma~2.1]{lu2025EEMMaxwell}, applied separately on $\Omega$ and
$\hat\Omega$.
\end{proof}
\begin{lemma}\label{lemma:penalized interpolation estimates}
For any $\bm u\in\bm H(\mathbf{curl},\mathcal D)
\cap\bm H^2(\Omega\cup\hat\Omega)$, there holds
\begin{equation}\label{penalized interpolation estimate}
\normtriple{\bm u-\pi_h\bm u}_\gamma
\lesssim
\kappa h^2\|\bm u\|_{\bm H^2(\Omega\cup\hat\Omega)}
+h\|\mathbf{curl}\,\bm u\|_{\bm H^1(\Omega\cup\hat\Omega)}.
\end{equation}
\end{lemma}
\begin{proof}
The local trace inequality and Lemma~\ref{interpolation theorem} give
\begin{align*}
&\sum_{F\in\mathcal F_h^I}h_F
\|\llbracket\mathbf{curl}(\bm u-\pi_h\bm u)\rrbracket\|_{\bm L^2(F)}^2\lesssim
\sum_{K\in\mathcal T_h^\Omega}
\left(
\|\mathbf{curl}(\bm u-\pi_h\bm u)\|_{\bm L^2(K)}^2
+h_K^2|\mathbf{curl}(\bm u-\pi_h\bm u)|_{\bm H^1(K)}^2
\right)\\
&\qquad\lesssim
h^2\|\mathbf{curl}\,\bm u\|_{\bm H^1(\Omega)}^2.
\end{align*}
The result now follows from Lemma~\ref{interpolation theorem}, the definition
of $\normtriple{\cdot}_\gamma$, and $|\gamma_F|\lesssim1$.
\end{proof}
We next define the sesquilinear forms used for the elliptic projections:
\begin{align*}
&b(\bm u,\bm v)= a(\bm u, \bm v)+2\kappa^2(\mathbb{A}^{-1}\bm u,\bm v),\\
&b_\gamma(\bm u,\bm v)=b(\bm u, \bm v)+J(\bm u,\bm v).
\end{align*}

\begin{lemma}\label{lemma:the continuity and coercive of bgamma}
There exists a constant $\alpha_0>0$  such that if $\Re(\gamma_F)\geq -\alpha_0$ and $|\gamma_F|\lesssim 1$,  then
\begin{align}
|b_\gamma(\bm{u},\bm{v})|&\lesssim \normtriple{\bm{u}}_\gamma\normtriple{\bm{v}}_\gamma,\quad \forall \ \bm{u},\bm{v}\in V,
\label{the continuity of bgamma}\\
\Re\, b_\gamma(\bm{v}_h,\bm{v}_h)&\gtrsim \normtriple{\bm{v}_h}^2_\gamma,\quad \forall \ \bm{v}_h\in V_h.
\label{the coercive of bgamma}
\end{align}
\end{lemma}
\begin{proof}
By using the Cauchy-Schwarz inequality,   the continuity \eqref{the continuity of bgamma} of $b_\gamma(\cdot,\cdot)$ follows from 
\begin{align*}
&|b(\bm{u},\bm{v})|\lesssim\normtriple{\bm{u}} \normtriple{\bm{v}},\quad\forall \bm{u},\bm{v}\in V\\
&|J(\bm{u},\bm{v})|\leq \sum_{F\in \mathcal{F}^{I}_h}|\gamma_F|h_F\|\llbracket\mathbf{curl}\,\bm u\rrbracket\|_{\bm{L}^2(F)} \|\llbracket\mathbf{curl}\,\bm v\rrbracket\|_{\bm{L}^2(F)}.
\end{align*}
The local trace inequality and the inverse inequality give a
constant $\beta_0>0$ such that
$$
 \sum_{F\in \mathcal{F}^{I}_h}h_F\|\llbracket\mathbf{curl}\,\bm{v}_h\rrbracket\|_{\bm{L}^2(F)}^2\leq \beta_0\|\mathbf{curl}\, \bm{v}_h\|_{\bm{L}^2(\Omega)}^2.
$$
That indicates for any $\bm{v}_h\in V_h$
\begin{align}\label{eqeeg}
\normtriple{\bm{v}_h}\leq \normtriple{\bm{v}_h}_\gamma\leq (1+\beta_0^{1/2}) \normtriple{\bm{v}_h}.
\end{align}
Denote by $\alpha_0=\frac{1}{2\beta_0}$ and let $\Re(\gamma_F)\geq -\alpha_0$, it is clear that
\begin{align*}
\Re\, b_\gamma(\bm{v}_h,\bm{v}_h)&=\normtriple{\bm{v}_h}^2+ \sum_{F\in \mathcal{F}^{I}_h}\Re(\gamma_F)h_F\|\llbracket\mathbf{curl}\,\bm{v}_h\rrbracket\|_{\bm{L}^2(F)}^2\geq \frac12\normtriple{\bm{v}_h}^2.
\end{align*}
which together with \eqref{eqeeg} yields \eqref{the coercive of bgamma} holds. That completes the proof.
\end{proof}

For any $\bm w\in V$, define its elliptic projections $P_h^{\pm}\bm w\in V_h$ as
\begin{align}\label{def:elliptic op}
b_\gamma(\bm w-P_h^+\bm w,\bm v_h)=0,\quad  b_\gamma(\bm v_h,\bm w-P_h^-\bm w)=0,\quad \forall\ \bm{v}_h\in V_h
\end{align}
\begin{lemma}\label{elliptic projection error estimate}
Suppose that  $\bm w\in \bm{H}_0(\mathbf{curl},\mathcal{D})\cap\bm H^2(\Omega\cup \hat{\Omega})$.
Then
\begin{align}\label{elliptic projection error estimate1}
\normtriple{\bm w-P_h^{\pm}\bm w}_\gamma
&\lesssim \kappa h^2\|\bm{w}\|_{\bm{H}^2(\Omega\cup \hat{\Omega})}+h \|\mathbf{curl}\,\bm{w}\|_{\bm{H}^1(\Omega\cup\hat\Omega)},\\
\label{elliptic projection error estimate2}
\|\bm w-P_h^{\pm}\bm w\|_{\bm{L}^2(\mathcal{D})}
&\lesssim C_{\kappa h} h^2\|\bm w\|_{\bm{H}^2(\Omega\cup\hat\Omega)}.
\end{align}
where the constant $C_{\kappa h}=1+\kappa h+\kappa^2h^2$. 
\end{lemma}
\begin{proof}
We present the error estimates for the elliptic operator $P_h^{+}$  while the proof for $P_h^{-}$ is omitted as it follows from an analogous argument, with the slight modification that the Helmholtz decomposition associated with $\mathbb{A}^{-H}$ is employed instead of $\mathbb{A}^{-1}$.

We first prove \eqref{elliptic projection error estimate1}. Recall that  for any $\bm{w}\in V$ and any $\bm{v}_h\in V_h$, there holds
\begin{align}\label{definition of P_h^+}
b_\gamma(\bm w-P_h^+\bm w,\bm v_h)=0.
\end{align}
From the lemma \ref{interpolation theorem} and triangle inequality, there holds
\begin{align}\label{triangle inequality: elliptic projection error estimate}
\normtriple{\bm{w}-P_h^{+}\bm w}_\gamma\leq \normtriple{\bm{w}-\pi_h\bm{w}}_\gamma+\normtriple{\pi_h\bm{w}-P_h^{+}\bm w}_\gamma.
\end{align}
Set $\bm\eta_h:=\pi_h\bm w-P_h^+\bm w$. By
\eqref{definition of P_h^+}, Lemma~\ref{lemma:the continuity and coercive of bgamma},
and \eqref{penalized interpolation estimate},
\begin{align*}
\normtriple{\bm\eta_h}_\gamma^2
&\lesssim \Re\,b_\gamma(\bm\eta_h,\bm\eta_h)
=\Re\,b_\gamma(\pi_h\bm w-\bm w,\bm\eta_h)\\
&\lesssim
\normtriple{\bm w-\pi_h\bm w}_\gamma
\normtriple{\bm\eta_h}_\gamma.
\end{align*}
Together with \eqref{triangle inequality: elliptic projection error estimate}
and \eqref{penalized interpolation estimate}, this proves
\eqref{elliptic projection error estimate1}.

Next, we prove \eqref{elliptic projection error estimate2}.
For brevity, write $\bm{\Psi}=\bm{w}-P_h^+\bm{w}\in V$, $\bm{\Phi}_h=P_h^+\bm{w}-\pi_h\bm{w}\in V_h$.
By the definition of $P_h^+$, we know
\begin{align}\label{orth:w-Pw}
(\mathbb{A}^{-1}\bm{\Psi},\nabla \xi_h)=\kappa^{-2}b_\gamma(\bm{\Psi},\nabla \xi_h)=0\quad \forall\ \xi_h\in S_h.
\end{align}
Applying  Lemma \ref{lemma:weighted continuous-discrete Helmholtz estimate} with $\mathbb{M}=\mathbb{A}^{-1}$, we have (the subscript $\mathbb{A}^{-1}$ is omitted here for simplification) 
\begin{align*}
\bm{\Phi}_h=\Pi^\nabla \bm{\Phi}_h+\Pi^{\mathbf{curl}}\bm{\Phi}_h=\Pi^\nabla_h \bm{\Phi}_h+\Pi^{\mathbf{curl}}_h\bm{\Phi}_h.
\end{align*}
Set
$\bm s:=\mathbb A^{-H}\Pi^{\mathbf{curl}}\bm\Phi_h$.
By \eqref{properties of weighted curl component},
$\mathbf{div}\,\bm s=0$.
Therefore, from \eqref{orth:w-Pw}  and $\Pi^\nabla_h \bm{\Phi}_h\in \nabla S_h$, we have
\begin{align}\label{decom-Psi}
(\mathbb{A}^{-1}\bm{\Psi},\bm{\Psi})
&=(\mathbb{A}^{-1}\bm{\Psi},\bm{w}-\pi_h\bm{w})-(\mathbb{A}^{-1}\bm{\Psi},\bm{\Phi}_h)\notag\\
&=(\mathbb{A}^{-1}\bm{\Psi},\bm{w}-\pi_h\bm{w})-(\mathbb{A}^{-1}\bm{\Psi},\Pi^{\mathbf{curl}}_h\bm{\Phi}_h-\Pi^{\mathbf{curl}}\bm{\Phi}_h)-(\mathbb{A}^{-1}\bm{\Psi},\Pi^{\mathbf{curl}}\bm{\Phi}_h)
\end{align}
To analyze the last term above. We introduce the auxiliary problem
\begin{align*}
\mathbf{curl}\,[\mathbb{A}\mathbf{curl}\, \overline{\bm{z}} ]+\kappa^2 \mathbb{A}^{-1} \overline{\bm{z}}& =\overline{\bm{s}},\quad \text{in}\ \mathcal{D}\\
\overline{\bm{z}}\times \hat{\bm x}&=\bm{0},\quad \text{on}\ \hat\Gamma.
\end{align*}
or in the variational form 
$$
b(\bm{v},\bm{z})=(\bm{v},\bm{s}),\quad \forall\ \bm{v}\in \bm{H}_0(\mathbf{curl},\mathcal{D}).
$$
Taking $\bm{v}=\bm{z}$. By the Young inequality, there holds that
$$
\kappa\|\bm{z}\|_{\bm{L}^2(\mathcal{D})}+\|\mathbf{curl}\,\bm{z}\|_{\bm{L}^2(\mathcal{D})}
\lesssim \kappa^{-1}\|\bm s\|_{\bm L^2(\mathcal{D})}.
$$
As in the proof of
Lemma~\ref{lemma:regularity of PML problem},
\cite[Theorems~2.7 and~2.10]{melenk2025} implies that if
$\mathbf{div}\,\bm s=0$, then
\begin{align}\label{+kstab}
\normtriple{\bm{z}}+\|\bm{z}\|_{\bm{H}^1(\Omega\cup\hat{\Omega})}+\kappa^{-1} \|\bm{z}\|_{\bm{H}^{2}(\Omega \cup \hat{\Omega})}& \lesssim \kappa^{-1}\|\bm{s}\|_{\bm{L}^2(\mathcal{D})}.
\end{align}
	Taking $\bm v=\bm\Psi$. Since
	$\bm z\in\bm H^2(\Omega\cup\hat\Omega)$, we have
	$J(\bm\Psi,\bm z)=0$. Thus \eqref{def:elliptic op} gives
	\begin{align*}
	(\mathbb{A}^{-1}\bm{\Psi},\Pi^{\mathbf{curl}}\bm{\Phi}_h)=b_\gamma(\bm{\Psi},\bm{z})=b_\gamma(\bm{\Psi},\bm{z}-\pi_h \bm{z}).
	\end{align*}
	By \eqref{the continuity of bgamma},
	\eqref{penalized interpolation estimate}, and \eqref{+kstab},
	\begin{align}\label{estimate of duality term for elliptic projection}
	|(\mathbb{A}^{-1}\bm{\Psi},\Pi^{\mathbf{curl}}\bm{\Phi}_h)|
&\lesssim (1+\kappa h) h\normtriple{\bm{\Psi}}_\gamma\|\Pi^{\mathbf{curl}}\bm{\Phi}_h\|_{\bm{L}^2(\mathcal{D})}\notag\\
&\lesssim (1+\kappa h) h\normtriple{\bm{\Psi}}_\gamma(\|\bm{\Psi}\|_{\bm{L}^2(\mathcal{D})}+\|\bm{w}-\pi_h\bm{w}\|_{\bm{L}^2(\mathcal{D})})
\end{align}
	Using \eqref{decom-Psi}, \eqref{curl-curlh}, and
	\eqref{estimate of duality term for elliptic projection}, together with
	Lemma~\ref{lemma:PML coefficient matrix properties} and Young's inequality,
	we obtain
	\begin{align}\label{intermediate L2 estimate for elliptic projection}
	\|\bm{\Psi}\|_{\bm{L}^2(\mathcal{D})}^2\lesssim |(\mathbb{A}^{-1}\bm{\Psi},\bm{\Psi})|\lesssim  \|\bm{w}-\pi_h\bm{w}\|_{\bm{L}^2(\mathcal{D})}^2+h^2\|\mathbf{curl}\, \bm{\Phi}_h\|_{\bm{L}^2(\mathcal{D})}^2+(1+\kappa h)^2 h^2\normtriple{\bm{\Psi}}_\gamma^2.
	\end{align}
	Moreover, the identity
	$\bm\Phi_h=-\bm\Psi+(\bm w-\pi_h\bm w)$,
	\eqref{elliptic projection error estimate1}, and
	\eqref{high regularity interpolation estimates} imply
	\begin{align}\label{curl and energy bounds for elliptic projection error}
	\|\mathbf{curl}\,\bm\Phi_h\|_{\bm L^2(\mathcal D)}
	+\normtriple{\bm\Psi}_\gamma
	\lesssim
	\kappa h^2\|\bm w\|_{\bm H^2(\Omega\cup\hat\Omega)}
	+h\|\mathbf{curl}\,\bm w\|_{\bm H^1(\Omega\cup\hat\Omega)}.
	\end{align}
	Combining \eqref{intermediate L2 estimate for elliptic projection},
	\eqref{curl and energy bounds for elliptic projection error}, and
	\eqref{high regularity interpolation estimates}, we obtain
	\begin{align*}
	\|\bm\Psi\|_{\bm L^2(\mathcal D)}
	&\lesssim h^2\|\bm w\|_{\bm H^2(\Omega\cup\hat\Omega)}
	+(1+\kappa h)h^2
	\left(
	\kappa h\|\bm w\|_{\bm H^2(\Omega\cup\hat\Omega)}
	+\|\mathbf{curl}\,\bm w\|_{\bm H^1(\Omega\cup\hat\Omega)}
	\right)\\
	&\lesssim C_{\kappa h}h^2
	\|\bm w\|_{\bm H^2(\Omega\cup\hat\Omega)}.
	\end{align*}
	which completes the proof.
\end{proof}

\subsection{Preasymptotic error estimates for CIP-EEM}
\begin{thm}
{Let Assumption~\ref{assumption:PML parameters} hold, and suppose that
the conditions of Lemma~\ref{lemma:the continuity and coercive of bgamma} are
satisfied.} Let $\bm E_h$ be the solution to \eqref{agamma}. There exists a constant $C_0$ independent of $\kappa$ and $h$ such that if 
$\kappa^3h^2\leq C_0,$
then the following stability estimation holds
\begin{align}\label{stab-CIP-EEM}
\|\mathbf{curl}\,\bm E_h\|_{\bm{L}^2(\mathcal{D})}+\kappa\|\bm E_h\|_{\bm{L}^2(\mathcal{D})}&\lesssim \|\bm f\|_{\bm{L}^2(\mathcal{D})},
\end{align}
\end{thm}
\begin{proof}
Let  $\bm E_h$ solves CIP-EEM problem \eqref{agamma}, equivalently,
\begin{align*}
b_\gamma(\bm{E}_h,\bm{E}_h)=(\bm{f},\bm{E}_h)+2\kappa^2(\mathbb{A}^{-1}\bm{E}_h,\bm{E}_h),\quad \forall\ \bm v\in  V_h.
\end{align*}
By the Young inequality, we deduce from  Lemma \ref{lemma:the continuity and coercive of bgamma} and \eqref{eqeeg} that
\begin{align}\label{curlEleqkE+f}
\|\mathbf{curl}\, \bm{E}_h\|^2_{\bm{L}^2(\mathcal{D})}&\leq \normtriple{\bm{E}_h}^2\lesssim |b_\gamma(\bm{E}_h,\bm{E}_h)|
\lesssim \kappa^2\|\bm{E}_h\|^2_{\bm{L}^2(\mathcal{D})}+\|\bm{f}\|^2_{\bm{L}^2(\mathcal{D})}
\end{align}
Now, we focus on the stability analysis for $\|\bm{E}_h\|^2_{\bm{L}^2(\mathcal{D})}$. Clearly, there holds
\begin{align}\label{agammaeq}
-\kappa^2(\mathbb{A}^{-1}\bm{E}_h,\nabla\bm{\xi}_h)=a_\gamma(\bm{E}_h,\nabla\bm{\xi}_h)=(\bm f,\nabla\bm{\xi}_h),\quad \forall \bm{\xi}_h\in S_h.
\end{align}
Applying  Lemma \ref{lemma:weighted continuous-discrete Helmholtz estimate} with $\mathbb{M}=\mathbb{A}^{-1}$, we have (the subscript $\mathbb{A}^{-1}$ is omitted here for simplification) \begin{align*}
\bm E_h=\Pi^\mathbf{\nabla}\bm E_h+\Pi^\mathbf{curl}\bm E_h=\Pi^\mathbf{\nabla}_h\bm E_h+\Pi^\mathbf{curl}_h\bm E_h.
\end{align*}
Since $\Pi^\mathbf{\nabla}_h\bm E_h\in \nabla S_h$, we deduce from \eqref{agammaeq} that
\begin{align}\label{decomEh}
&(\mathbb{A}^{-1}\bm E_h,\bm E_h)
=(\mathbb{A}^{-1}\bm E_h,\Pi^\mathbf{\nabla}_h\bm E_h)+(\mathbb{A}^{-1}\bm E_h,\Pi_h^\mathbf{curl}\bm E_h)\notag\\
=&-\kappa^{-2}(\bm f,\Pi^\mathbf{\nabla}_h\bm E_h)
+(\mathbb{A}^{-1}\bm{E}_h,\Pi_h^\mathbf{curl}\bm E_h-\Pi^\mathbf{curl}\bm E_h)+(\mathbb{A}^{-1}\bm E_h,\Pi^\mathbf{curl}\bm E_h)
\end{align}
For the second term of \eqref{decomEh}, from \eqref{curl-curlh} and  \eqref{curlEleqkE+f}, there holds
\begin{align}\label{estimate of curl -curlh Eh}
\|\Pi_h^\mathbf{curl}\bm E_h-\Pi^\mathbf{curl}\bm E_h\|_{\bm{L}^2(\mathcal{D})}\lesssim h\|\mathbf{curl}\,\bm E_h\|_{\bm{L}^2(\mathcal{D})}
\leq \kappa h\left(\|\bm E_h\|_{\bm{L}^2(\mathcal{D})}+\kappa^{-1}\|\bm f\|_{\bm{L}^2(\mathcal{D})}\right).
\end{align}
To analyze the last term above, we 
 consider the following dual problem
\begin{subequations}\label{dual Maxwell equation}
\begin{align}\label{dual Maxwell equation 1}
\mathbf{curl}\,[\mathbb{A}\mathbf{curl}\, \overline{\bm z} ]-\kappa^2 \mathbb{A}^{-1} \overline{\bm z}& =\overline{\bm s},\quad \text{in}\ \mathcal{D},\\
\overline{\bm z}\times \hat{\bm x}&=\bm{0},\quad \text{on}\ \hat\Gamma.\label{dual Maxwell equation 2}
\end{align}
\end{subequations}
where $\mathbf{div} \overline{\bm s}=\overline{\mathbf{div} \bm s}=0.$ And equivalently, in the  variational form 
$$
a(\bm{v},\bm{z})=(\bm{v},\bm{s}),\quad \forall\ \bm{v}\in \bm{H}_0(\mathbf{curl},\mathcal{D}).
$$
Theorem~\ref{lemma:L2 stability of PML Maxwell equations} and
Lemma~\ref{lemma:regularity of PML problem} with $m=0$ give the following
estimate for \eqref{dual Maxwell equation}:
\begin{align}\label{-kstab}
\|\bm{z}\|_{\bm{H}^2(\Omega\cup \hat{\Omega})}+\|\mathbf{curl}\,\bm{z}\|_{\bm{H}^1(\Omega\cup \hat{\Omega})}+\kappa\|\mathbf{curl}\, \bm{z}\|+\kappa\|\bm{z}\|_{\bm{H}^1(\Omega\cup \hat{\Omega})}+\kappa^2\|\bm{z}\|_{\bm{L}^2(\mathcal{D})}\lesssim \kappa \|\bm s\|_{\bm{L}^2(\mathcal{D})}
\end{align}
Take $\bm s=\mathbb{A}^{-H}\Pi^\mathbf{curl}\bm E_h$ and
$\bm v=\bm E_h$. Since $\bm z\in\bm H^2(\Omega\cup\hat\Omega)$, we have  $J(\bm E_h,\bm z)=0$ and then \eqref{def:elliptic op} give
\begin{align*}
(\mathbb{A}^{-1}\bm E_h,\Pi^\mathbf{curl}\bm E_h)&=a(\bm E_h, \bm z)=a_\gamma(\bm E_h, \bm z)\\
&=a_\gamma(\bm E_h,P_h^{-}\bm{z})+a_\gamma(\bm E_h, \bm{z}-P_h^{-}\bm{z})\\
&=(\bm{f},P_h^-\bm{z})-2\kappa^2(\mathbb{A}^{-1}\bm{E}_h,\bm{z}-P_h^{-}\bm{z})\\
&=(\bm{f},\bm{z})-(\bm{f},\bm{z}-P_h^-\bm{z})-2\kappa^2(\mathbb{A}^{-1}\bm{E}_h,\bm{z}-P_h^{-}\bm{z}).
\end{align*}
By the Cauchy-Schwarz inequality and Lemma \ref{elliptic projection error estimate} and \eqref{-kstab}, we obtain
\begin{align}\label{prod_EPcurlE}
|(\mathbb{A}^{-1}\bm E_h,\Pi^\mathbf{curl}\bm E_h)|\lesssim C_{\kappa h}
\left(
\kappa^3h^2\|\bm{E}_h\|_{\bm{L}^2(\mathcal{D})}+\kappa^{-1}\|\bm{f}\|_{\bm{L}^2(\mathcal{D})}
\right)
\|\Pi^\mathbf{curl}\bm{E}_h\|_{\bm{L}^2(\mathcal{D})}.
\end{align}
From \eqref{prod_EPcurlE} and \eqref{estimate of curl -curlh Eh}, by applying the Young inequality to  \eqref{decomEh} , we conclude that
\begin{align}\label{Eleqk3h2E+f}
\|\bm{E}_h\|_{\bm{L}^2(\mathcal{D})}\leq C_{\kappa h}
\left(
 \kappa^{-1}\|\bm f\|_{\bm{L}^2(\mathcal{D})}+\kappa h\|\bm E_h\|_{\bm{L}^2(\mathcal{D})}+\kappa^3h^2\|\bm{E}_h\|_{\bm{L}^2(\mathcal{D})}
\right).
\end{align}
Since $\kappa\gtrsim1$, the mesh condition
$\kappa^3h^2\leq C_0$ implies $C_{\kappa h}\lesssim1$. Therefore, for
sufficiently small $C_0$, the resulting
$L^2$ estimate \eqref{Eleqk3h2E+f} together with \eqref{curlEleqkE+f}, proves
\eqref{stab-CIP-EEM}.
\end{proof}

\begin{thm}
Let Assumption~\ref{assumption:PML parameters} and the assumptions of
Lemma~\ref{lemma:the continuity and coercive of bgamma} hold, and let
$\bm E_h$ solve \eqref{agamma}. There exists a constant $C_0$ independent of $\kappa$ and $h$ such that if $$\kappa^3h^2\leq C_0,$$ then
the preasymptotic error estimation holds: 
\begin{enumerate}
\item For $\bm{f}\in  \bm{H}(\mathbf{div}^0,\mathcal{D})$, we have
\begin{align}
\interleave\hat{\bm E}-\bm E_h\interleave_\gamma&\lesssim (\kappa h+\kappa^3h^2)\|\bm f\|_{\bm{L}^2(\mathcal{D})},\\
{\kappa}\|\hat{\bm E}-\bm E_h\|_{\bm L^2(\mathcal{D})}&\lesssim (\kappa^2h^2+\kappa^3h^2)\|\bm f\|_{\bm{L}^2(\mathcal{D})}.
\end{align}
\item For $\bm{f}\in  \bm{H}(\mathbf{div},\mathcal{D})\cap\bm{H}^1(\mathbf{div},\Omega\cup\hat{\Omega})$, we have
\begin{align}
\interleave\hat{\bm E}-\bm E_h\interleave_\gamma&\lesssim (\kappa h+\kappa^3h^2)C_{2,\bm{f}},\\
{\kappa}\|\hat{\bm E}-\bm E_h\|_{\bm L^2(\mathcal{D})}&\lesssim (\kappa^2h^2+\kappa^3h^2)C_{2,\bm{f}}.
\end{align}
\end{enumerate}
\end{thm}
\begin{proof}
Set $\hat{\bm E}-\bm E_h=\bm\eta-\bm\xi_h$ with $\bm\eta:=\hat{\bm E}-P_h^+\hat{\bm E}$ and $\bm\xi_h:=\bm E_h-P_h^+\hat{\bm E}$.  By
\eqref{CIP Galerkin orthogonality} and the definition of $P_h^+$, we have
\begin{align*}
a_\gamma(\bm\xi_h,\bm v_h)
&=a_\gamma(\bm\eta,\bm v_h)
=-2\kappa^2(\mathbb A^{-1}\bm\eta,\bm v_h),
\qquad \forall\,\bm v_h\in V_h.
\end{align*}
Thus $\bm\xi_h$ is the CIP-EEM solution with source
$-2\kappa^2\mathbb A^{-1}\bm\eta$. The discrete stability estimate
\eqref{stab-CIP-EEM}, the norm equivalence \eqref{eqeeg}, and the uniform
boundedness of $\mathbb A^{-1}$ give
\begin{equation}\label{stability estimate of discrete projection remainder}
\normtriple{\bm\xi_h}_\gamma
+\kappa\|\bm\xi_h\|_{\bm L^2(\mathcal D)}
\lesssim
\kappa^2\|\bm\eta\|_{\bm L^2(\mathcal D)}.
\end{equation}
Since
$\|\mathbf{curl}\,\hat{\bm E}\|_{\bm H^1(\Omega\cup\hat\Omega)}
\lesssim\|\hat{\bm E}\|_{\bm H^2(\Omega\cup\hat\Omega)}$,
\eqref{elliptic projection error estimate1}--
\eqref{elliptic projection error estimate2} yield
\begin{align}
\normtriple{\hat{\bm E}-\bm E_h}_\gamma
&\lesssim
(h+\kappa^2h^2)
\|\hat{\bm E}\|_{\bm H^2(\Omega\cup\hat\Omega)},
\label{abstract preasymptotic energy error estimate}\\
\kappa\|\hat{\bm E}-\bm E_h\|_{\bm L^2(\mathcal D)}
&\lesssim
(\kappa h^2+\kappa^2h^2)
\|\hat{\bm E}\|_{\bm H^2(\Omega\cup\hat\Omega)}.
\label{abstract preasymptotic L2 error estimate}
\end{align}
Here we used $C_{\kappa h}\lesssim1$ and $\kappa\gtrsim1$, which follow from
the mesh condition and Assumption~\ref{assumption:PML parameters}.
Finally, \eqref{Hm+2estimateofmaxwell} with $m=0$ gives
$\|\hat{\bm E}\|_{\bm H^2}\lesssim\kappa\|\bm f\|_{\bm L^2}$ in the
first case, while \eqref{Hm+1estimateofmaxwell} with $m=1$ gives
$\|\hat{\bm E}\|_{\bm H^2}\lesssim\kappa C_{2,\bm f}$ in the second.
Substitution into the above estimates completes the proof.
\end{proof}

\section{Numerical experiments}\label{sec:numerical-experiments}
In this section, numerical experiments are presented for the spherical PML problem to investigate the
convergence rates and pollution errors of the EEM and CIP-EEM. The conforming
tetrahedral meshes are generated by Gmsh \cite{GeuzaineRemacle2009}.

\subsection{A single-mode VSH exact solution}
{
\setlength{\abovedisplayskip}{6pt}
\setlength{\belowdisplayskip}{6pt}
\setlength{\abovedisplayshortskip}{4pt}
\setlength{\belowdisplayshortskip}{4pt}
To simulate the PML truncated problem \eqref{stretched Maxwell equation},
we take $\Omega=B_{0.5}(\bm 0)$ and $\mathcal D=B_1(\bm 0)$. Set the PML parameter
$\sigma_0=2$. And the exact solution $\hat{\bm{E}}$ and source $\bm{f}$ is chosen to be the single VSH mode corresponding
to $(l,m)=(2,0)$ as follows:
\begin{equation}\label{numerical phase matched exact field}
  \hat{\bm E}
  =\alpha\hat E_{20}^{(1)}\bm Y_{20}
   +\beta\hat E_{20}^{(2)}\bm U_{20},\quad \bm{f}=f_{20}^{(1)}\bm Y_{20}
    +f_{20}^{(2)}\bm U_{20}.
\end{equation}
Here $\hat E_{20}^{(1)}$ and $\hat E_{20}^{(2)}$ are given in
\eqref{expression of PML solution hatE1}--\eqref{expression of PML solution hatE2}. The corresponding source coefficients are
\begin{equation}\label{numerical phase matched source coefficients}
\bigl(f_{20}^{(1)}(r),f_{20}^{(2)}(r)\bigr)
=\left(\frac{\sqrt6}{r}\chi(r)j_2(\kappa r),\,
\chi'(r)j_2(\kappa r)+\chi(r)\left(\frac{j_2(\kappa r)}{r}
+\kappa j_2'(\kappa r)\right)\right),
\end{equation}
where $\chi(r)$ is the cutoff function supported in $\Omega$ defined by
\begin{equation}\label{numerical phase matched cutoff}
  \chi(r)=
  \begin{cases}
    \left(1-\dfrac{r^2}{\rho^2}\right)^4,&0\leq r<\rho,\\[1mm]
    0,&r\geq\rho,
  \end{cases}
  \qquad \rho=0.68R<R.
\end{equation}
It follows from \eqref{components' div condition}
that
$
\bm{f}\in \bm{H}(\mathbf{div}^0,\mathcal{D})
$ and then from Lemma~\ref{lemma:regularity of PML problem} that
$\hat{\bm E}\in \bm{H}^2(\Omega\cup\hat{\Omega})$.
}
\begingroup
\setlength{\intextsep}{8pt}
\subsection{Convergence and pollution tests}
{
We discretize the truncated PML problem by the first-order N\'ed\'elec edge
element of the second type.  The relative energy error is measured by
\begin{equation}\label{numerical relative energy error}
  \mathcal E_h
  :=\frac{\normtriple{\hat{\bm E}-\bm E_h}}
  {\normtriple{\hat{\bm E}}}.
\end{equation}
To quantify the mesh resolution, we introduce the global number of degrees
of freedom per wavelength
\begin{equation}\label{numerical global nlambda}
  N_\lambda
  :=\frac{2\pi}{\kappa}
  \left(\frac{\operatorname{ndof}}{|\mathcal D|}\right)^{1/3}.
\end{equation}
The penalty parameters for the CIP-EEM, for simplicity, are chosen as
$
\gamma_F=-\frac{\sqrt{2}}{24}-0.01\mathrm{i}, F\in\mathcal F_h^{I}.
$
The real part follows from the dispersion analysis of the Helmholtz
CIP-EEM \cite{ZhouWu2023} for the uniform mesh, whereas the small negative imaginary part is
added, as in \cite{lu2024preasymptoticerrorestimateseem}, to
improve stability. Regarding the choice of penalty parameters, we also refer the reader to the
optimization algorithm in \cite{LiLiYang2024}, developed for the Helmholtz
CIP-EEM, which may be useful in further
reducing the pollution error of the CIP-EEM on general unstructured meshes.
\begin{figure}[htbp]
  \centering
  \begin{minipage}[t]{0.49\textwidth}
    \centering
    \includegraphics[width=\linewidth]
    {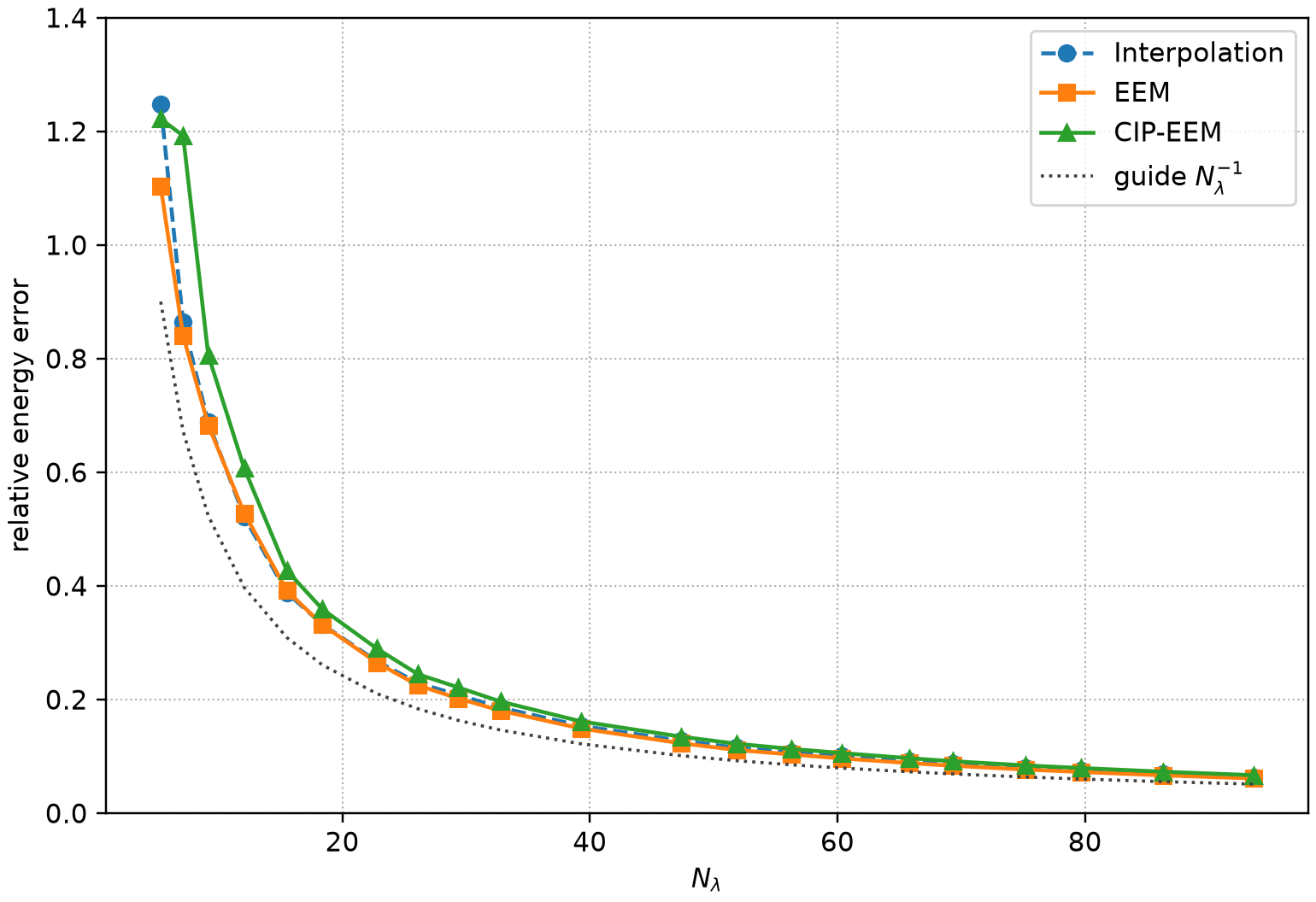}
    \par\smallskip
    (a) $\kappa=10$.
  \end{minipage}\hfill
  \begin{minipage}[t]{0.49\textwidth}
    \centering
    \includegraphics[width=\linewidth]
    {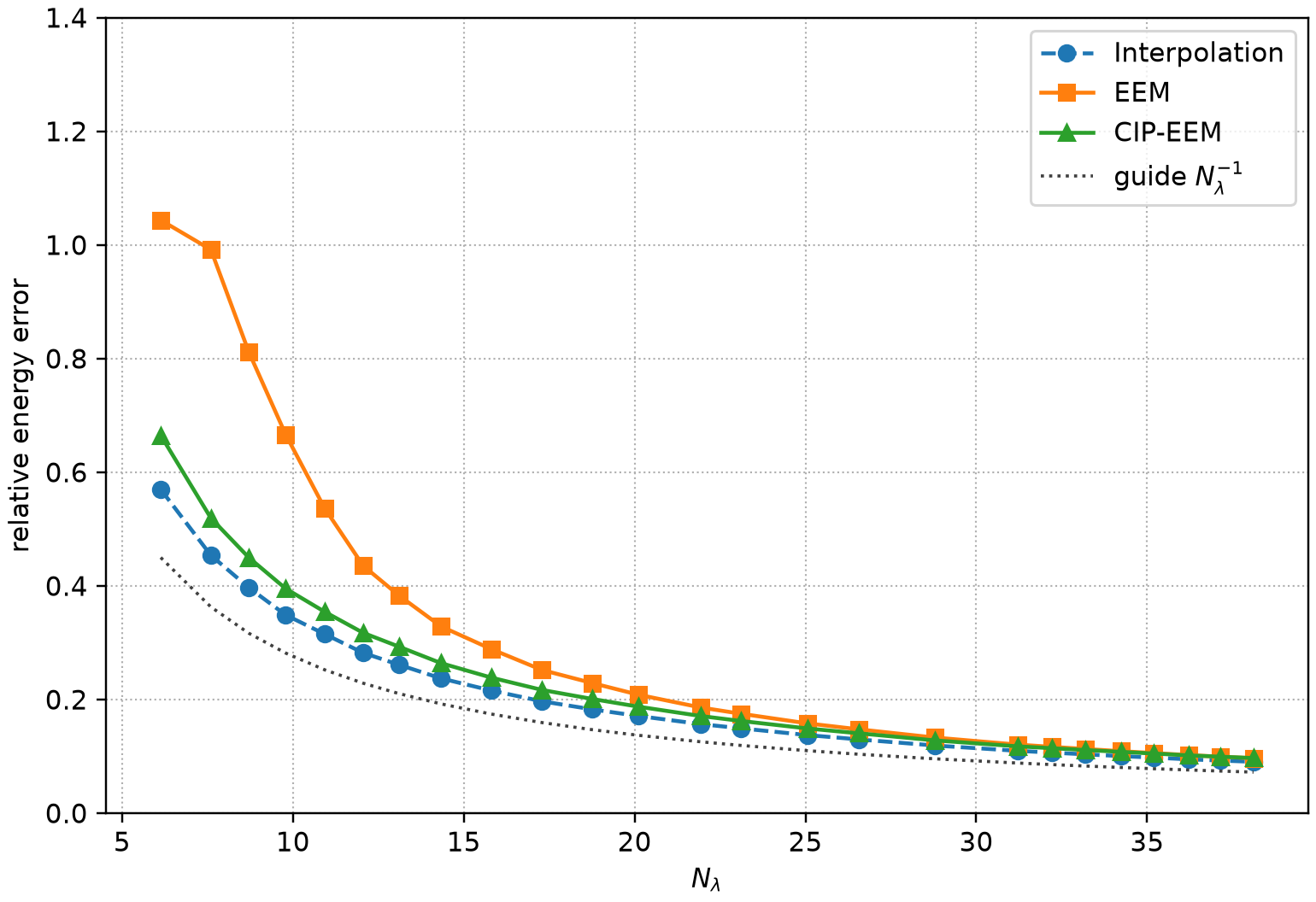}
    \par\smallskip
  (b) $\kappa=30$.
  \end{minipage}
  \caption{Relative energy errors of the EEM and CIP-EEM solutions and the
N\'ed\'elec interpolation versus $N_\lambda$ for (a) $\kappa=10$ and
(b) $\kappa=30$. The dotted line indicates the reference decay rate
$N_\lambda^{-1}$.}
  \label{fig:vsh-convergence-energy}
\end{figure}

Figure~\ref{fig:vsh-convergence-energy} plots the relative energy errors of
the EEM and CIP-EEM solutions and the N\'ed\'elec interpolation. For
$\kappa=10$, the three curves nearly coincide and decay with the expected order
$N_\lambda^{-1}$. For $\kappa=30$, the EEM error is considerably larger than
the interpolation error on coarse meshes owing to the pollution effect, but
gradually approaches the interpolation error as $N_\lambda$ increases. This behavior shows that, as the mesh is refined, the pollution effect is gradually under control and the finite element error  approaches its asymptotic regime.

We next fix $N_\lambda=12$ and increase the wave number $\kappa$. As shown in
Figure~\ref{fig:vsh-pollution-energy}, apart from the initial under-resolved
range, the interpolation error remains close to $0.25$, indicating that
$N_\lambda$ effectively measures the mesh resolution. The EEM error
increases overall with $\kappa$ because of pollution errors whereas the CIP-EEM
error is consistently smaller than the EEM error, indicating that the penalty
term reduces the pollution errors. The oscillations in the EEM and CIP-EEM error curves are mainly caused by
numerical reflections induced by the piecewise-constant PML profile. By
contrast, a smooth power-law PML profile reduces these reflections, leading
to smoother EEM and CIP-EEM error curves.
\begin{figure}[h]
  \centering
  \includegraphics[width=0.6\textwidth]
  {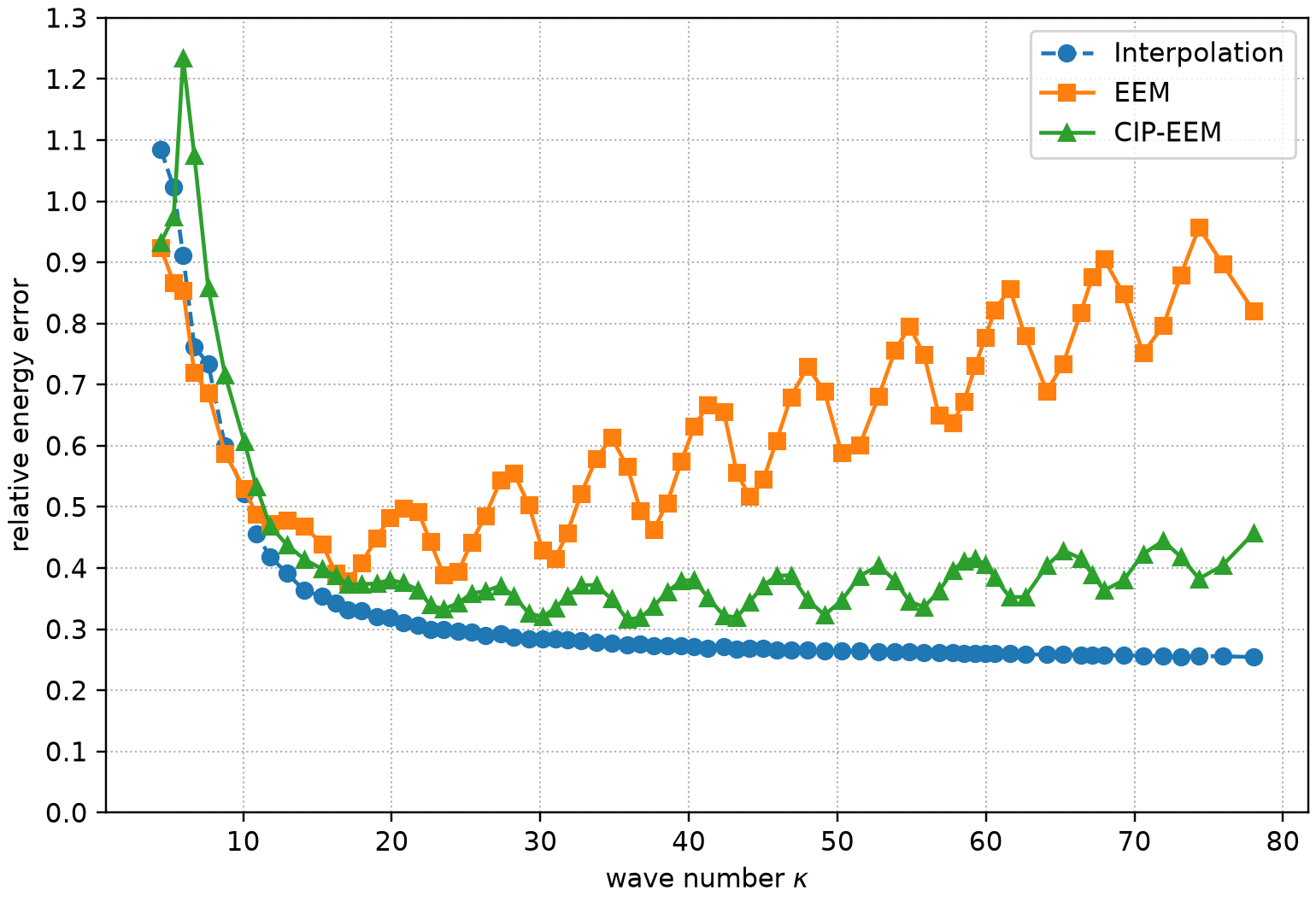}
  \caption{Relative energy errors of the  N\'ed\'elec
  interpolation, EEM, and CIP-EEM at fixed global resolution
  $N_\lambda=12$.}
  \label{fig:vsh-pollution-energy}
\end{figure}}
\endgroup

\appendix
\numberwithin{equation}{section}
\section{Auxiliary derivations in spherical coordinates}
\label{app:spherical-derivations}
\subsection{Derivation of the spherical PML formulation}\label{app:pml-derivation}
In spherical coordinates, the local orthonormal basis
$\{\bm e_r,\bm e_\theta,\bm e_\phi\}$ is defined by
\begin{equation}\label{spherical unit vectors}
\begin{aligned}
\bm e_r&=(\sin\theta\cos \phi,\sin\theta\sin\phi,\cos\theta)^T=\bm x/r,\\
\bm e_\theta&=(\cos\theta\cos\phi,\cos\theta\sin\phi,-\sin\theta)^T,\\
\bm e_\phi&=(-\sin\phi,\cos\phi,0)^T .
\end{aligned}
\end{equation}
For a vector field
\[
\bm v=v_r\bm e_r+v_\theta\bm e_\theta+v_\phi\bm e_\phi,
\]
the curl and divergence are given by
\begin{equation}\label{spherical curl div formulas}
\begin{aligned}
\mathbf{curl}\,\bm v
={}&\frac1{r\sin\theta}
\left(\frac{\partial}{\partial\theta}(\sin\theta v_\phi)
-\frac{\partial v_\theta}{\partial\phi}\right)\bm e_r\\
&+\frac1r
\left(\frac1{\sin\theta}\frac{\partial v_r}{\partial\phi}
-\frac{\partial}{\partial r}(rv_\phi)\right)\bm e_\theta
+\frac1r
\left(\frac{\partial}{\partial r}(rv_\theta)
-\frac{\partial v_r}{\partial\theta}\right)\bm e_\phi,\\
\mathbf{div}\,\bm v
={}&\frac1{r^2}\frac{\partial}{\partial r}(r^2v_r)
+\frac1{r\sin\theta}\frac{\partial}{\partial\theta}(\sin\theta v_\theta)
+\frac1{r\sin\theta}\frac{\partial v_\phi}{\partial\phi}.
\end{aligned}
\end{equation}
Let $\tilde{\bm E}(r,\theta,\phi)=\bm E(\tilde r,\theta,\phi)$ denote the
stretched field. Then $\tilde{\bm E}$ satisfies the stretched equation
\begin{equation}\label{PML equation whole space1}
\widetilde{\mathbf{curl}}\,\widetilde{\mathbf{curl}}\,\tilde{\bm E}
-\kappa^2\tilde{\bm E}
=\bm f\quad \text{in}\ \mathbb{R}^3,
\end{equation}
where $\widetilde{\mathbf{curl}}$ denotes the curl with respect to the stretched
coordinate $\tilde{\bm x}$. Write
\[
\tilde{\bm E}
=\tilde E_r\bm e_r+\tilde E_\theta\bm e_\theta+\tilde E_\phi\bm e_\phi .
\]
In term of $\alpha=\mathrm d\tilde r/\mathrm d r, 
\beta=\tilde r/r$,
a direct calculation by replacing $r$ by $\tilde{r}$ in
\eqref{spherical curl div formulas} gives
\begin{equation}\label{curl div transformation under PML stretching}
\widetilde{\mathbf{curl}}\,\tilde{\bm E}
=J^{-1}\mathbb B\mathbf{curl}(\mathbb B^T\tilde{\bm E}),
\qquad
\widetilde{\mathbf{div}}\,\tilde{\bm E}
=J^{-1}\mathbf{div}(J\mathbb B^{-1}\tilde{\bm E}).
\end{equation}
Here
\begin{equation}\label{Jacobian matrix of complex stretching}
\mathbb B=H\operatorname{diag}\{\alpha,\beta,\beta\}H^T,\qquad
J=\det(\mathbb B)=\alpha\beta^2 .
\end{equation}
are the Jacobi matrix of the radial transform
$\bm x\mapsto\tilde{\bm x}(\bm x)$ and its determinant. Set
\[
\bm E^{\mathrm{pml}}:=\mathbb B^T\tilde{\bm E}
=\alpha \tilde E_r\bm e_r
+\beta \tilde E_\theta\bm e_\theta+\beta \tilde E_\phi\bm e_\phi .
\]
Since $\operatorname{supp}\bm f\subset\Omega$, we have
$J\mathbb B^{-1}\bm f=\bm f$. Combining \eqref{PML equation whole space1} with
\eqref{curl div transformation under PML stretching}  gives the  PML equation
\begin{equation}\label{full-space PML equation}
\mathbf{curl}(\mathbb A\mathbf{curl}\,\bm E^{\mathrm{pml}})
-\kappa^2 \mathbb A^{-1}\bm E^{\mathrm{pml}}
=\bm f,\quad \text{in}\ \mathbb R^3.
\end{equation}
where the coefficient matrix 
\begin{equation}\label{PML coefficient A by Jacobian}
\mathbb A:=J^{-1}\mathbb B^T\mathbb B
=H\operatorname{diag}\left\{
\frac{\alpha}{\beta^2},\frac1\alpha,\frac1\alpha
\right\}H^T .
\end{equation}

{
\begin{lemma}\label{lemma:PML coefficient matrix properties}
The matrix $\mathbb A$ defined by
\eqref{PML coefficient A by Jacobian} is piecewise smooth and symmetric. Moreover, the matrix $\mathbb{A}$ and $\mathbb{A}^{-1}$ are uniformly coercive and bounded. Precisely, for every
$\bm\xi\in\mathbb C^3$,
\[
\Re(\mathbb A\bm\xi\cdot\overline{\bm\xi})
\geq (1+\sigma_0^2)^{-2}|\bm\xi|^2,
\qquad
\Re(\mathbb A^{-1}\bm\xi\cdot\overline{\bm\xi})
\geq (1+\sigma_0^2)^{-1}|\bm\xi|^2,
\]
and
\[
|\mathbb A\bm\xi|
\leq (1+\sigma_0^2)^{1/2}|\bm\xi|,
\qquad
|\mathbb A^{-1}\bm\xi|
\leq (1+\sigma_0^2)|\bm\xi|.
\]
\end{lemma}
\begin{proof}
The result follows directly from \eqref{sigma and delta} and
\eqref{PML coefficient A by Jacobian}. We omit the proof here.
\end{proof}
}

{Since $\alpha=\beta=1$ in $\Omega$, one has
$\bm E^{\mathrm{pml}}=\tilde{\bm E}=\bm E$ in $\Omega$. Moreover, the Helmholtz
potential representation \eqref{Maxwell solution representation by Helmholtz potential},
together with the exponential decay estimates for the corresponding stretched
Helmholtz potentials \cite{MR3899165}, gives the exponential decay of the
stretched field through the PML layer. We therefore truncate the layer at
$r=\hat R$ and impose
the perfect conductor condition on the outer boundary. This gives the truncated
PML problem \eqref{stretched Maxwell equation}, whose solution is denoted by
$\hat{\bm E}$.
}

\subsection{Auxiliary identities for vector spherical harmonics}\label{Appendix curl w}
For any vector field $\bm{w}=(w_1,w_2,w_3)^T$, 
\[
\bm{w}=w_r\bm{e}_r+w_\theta\bm{e}_\theta+w_\phi\bm{e}_\phi .
\]
{Using \eqref{spherical unit vectors}--\eqref{spherical curl div formulas}, we
derive the following identities. First, by}
$$
\nabla_{\mathbb{S}^2}Y_l^m=\frac{\partial Y_l^m}{\partial \theta}\bm e_\theta+\frac1{\sin\theta}\frac{\partial Y_l^m}{\partial \phi}\bm e_\phi,
$$
and the VSH expansions \eqref{Fourier expansions}, we have
\begin{align*}
&w_r=\sum_{l=0}^{+\infty}\sum_{m=-l}^{l}Y_l^mw_{lm}^{(1)},\\
\sqrt{l(l+1)}&w_\theta=\sum_{l=0}^{+\infty}\sum_{m=-l}^{l}\left(\frac{\partial Y_l^m}{\partial \theta}w_{lm}^{(2)}+\frac1{\sin\theta}\frac{\partial Y_l^m}{\partial \phi}w_{lm}^{(3)}\right)\\
\sqrt{l(l+1)}&w_\phi=\sum_{l=0}^{+\infty}\sum_{m=-l}^{l}\left(\frac1{\sin\theta}\frac{\partial Y_l^m}{\partial \phi}w_{lm}^{(2)}-\frac{\partial Y_l^m}{\partial \theta}w_{lm}^{(3)}\right).
\end{align*}
{Moreover,}
\begin{align*}
\frac1{r\sin\theta}\left(\frac{\partial}{\partial \theta}(\sin\theta w_\phi)-\frac{\partial w_\theta}{\partial \phi}\right)\bm e_r
&=\sum_{l=0}^{+\infty}\sum_{m=-l}^{l}\frac{w_{lm}^{(3)}}{r\sqrt{l(l+1)}}(-\Delta_{\mathbb{S}^2}Y_l^m) \bm e_r\\
&=\sum_{l=0}^{+\infty}\sum_{m=-l}^{l}\frac{\sqrt{l(l+1)}}{r} w_{lm}^{(3)} \bm{Y}_{lm},
\end{align*}

\begin{align*}
\frac1{r\sin\theta}\frac{\partial w_r}{\partial \phi}\bm e_\theta-\frac1r\frac{\partial w_r}{\partial \theta}\bm e_\phi
&=\sum_{l=0}^{+\infty}\sum_{m=-l}^{l}\frac{w_{lm}^{(1)}}{r}\left[(\nabla_{\mathbb{S}^2}Y_l^m\cdot\bm e_\phi)\bm e_\theta-(\nabla_{\mathbb{S}^2}Y_l^m\cdot\bm e_\theta)\bm e_\phi\right]\\
&=\sum_{l=0}^{+\infty}\sum_{m=-l}^{l}\frac{\sqrt{l(l+1)}}rw_{lm}^{(1)}\bm V_{lm},
\end{align*}
\begin{align*}
-\frac1r\frac{\partial}{\partial r}(rw_\phi)\bm e_\theta+\frac1r\frac{\partial}{\partial r}(rw_\theta)\bm e_\phi
&=
\sum_{l=0}^{+\infty}\sum_{m=-l}^{l}\frac1{r\sqrt{l(l+1)}}\frac{\mathrm d (rw_{lm}^{(2)})}{\mathrm d r}\left[-(\nabla_{\mathbb{S}^2}Y_l^m\cdot\bm e_\phi)\bm e_\theta+(\nabla_{\mathbb{S}^2}Y_l^m\cdot\bm e_\theta)\bm e_\phi\right]\\
&\ \ \ +\sum_{l=0}^{+\infty}\sum_{m=-l}^{l}
\frac1{r\sqrt{l(l+1)}}\frac{\mathrm d (rw_{lm}^{(3)})}{\mathrm d r}\left[(\nabla_{\mathbb{S}^2}Y_l^m\cdot\bm e_\theta)\bm e_\theta+(\nabla_{\mathbb{S}^2}Y_l^m\cdot\bm e_\phi)\bm e_\phi\right]\\
&
=\sum_{l=0}^{+\infty}\sum_{m=-l}^{l}\frac{1}r\left[-\frac{\mathrm d (rw_{lm}^{(2)})}{\mathrm d r}\bm V_{lm}+\frac{\mathrm d (rw_{lm}^{(3)})}{\mathrm d r}\bm U_{lm}\right].
\end{align*}
\begin{align*}
\frac{1}{r\sin \theta}\frac{\partial}{\partial \theta}(\sin \theta w_\theta)+\frac1{r\sin \theta}\frac{\partial w_\phi}{\partial \phi}
=\sum_{l=0}^{+\infty}\sum_{m=-l}^{l}\frac{\Delta_{\mathbb{S}^2}Y_l^m}{r\sqrt{l(l+1)}}w_{lm}^{(2)}=-\sum_{l=0}^{+\infty}\sum_{m=-l}^{l}\frac{\sqrt{l(l+1)}}rw_{lm}^{(2)}Y_l^m.
\end{align*}
{Thus,}
\begin{align}
&\mathbf{curl}\,\bm{w}=\sum_{l=0}^{+\infty}\sum_{m=-l}^{l}\left[\frac{\sqrt{l(l+1)}}r w_{lm}^{(3)} \bm{Y}_{lm}+\frac1r\frac{\mathrm d (rw_{lm}^{(3)})}{\mathrm d r}\bm U_{lm}+\left(\frac{\sqrt{l(l+1)}}rw_{lm}^{(1)}-\frac1r\frac{\mathrm d (rw_{lm}^{(2)})}{\mathrm d r}\right)\bm V_{lm}\right],\notag\\
&\mathbf{curl}\,\bm{w}\times \hat{\bm x}=\sum_{l=0}^{+\infty}\sum_{m=-l}^{l}\left[\frac1r\frac{\mathrm d (rw_{lm}^{(3)})}{\mathrm d r}{\bm V_{lm}}-\left(\frac{\sqrt{l(l+1)}}rw_{lm}^{(1)}-\frac1r\frac{\mathrm d (rw_{lm}^{(2)})}{\mathrm d r}\right)\bm U_{lm}\right],\label{VSH tangential curl expansion}\\
&\mathbf{div}\,\bm w=\sum_{l=0}^{+\infty}\sum_{m=-l}^{l}\left[\frac1{r^2}\frac{\mathrm d (r^2w_{lm}^{(1)})}{\mathrm d r}-\frac{\sqrt{l(l+1)}}rw_{lm}^{(2)}\right]Y_l^m.\label{VSH divergence expansion}
\end{align}
{Repeating the same computation gives}
\begin{equation}\label{VSH curl curl expansion}
\begin{aligned}
\mathbf{curl}\,\mathbf{curl}\,\bm{w}
={}&\sum_{l=0}^{+\infty}\sum_{m=-l}^{l}\left[\frac{\sqrt{l(l+1)}}{r^2}\left(\sqrt{l(l+1)}w_{lm}^{(1)}-\frac{\mathrm d (rw_{lm}^{(2)})}{\mathrm d r}\right) \bm{Y}_{lm}\right.\\
&+\left.\frac1r\left( \sqrt{l(l+1)}\frac{\mathrm d w_{lm}^{(1)}}{\mathrm d r}-\frac{\mathrm d^2 (rw_{lm}^{(2)})}{\mathrm d r^2}\right)\bm{U}_{lm}\right.\\
&+\left.\left(\frac{l(l+1)}{r^2}w_{lm}^{(3)}-\frac1r\frac{\mathrm d^2 (rw_{lm}^{(3)})}{\mathrm d r^2}\right)\bm{V}_{lm}\right].
\end{aligned}
\end{equation}

\section{Proof of the weighted Helmholtz estimate}
\label{appendix:weighted Helmholtz estimate}
\begin{proof}[Proof of Lemma~\ref{lemma:weighted continuous-discrete Helmholtz estimate}]
Set $\bm z:=\Pi^{\mathbf{curl}}_{\mathbb{M}}\bm w_h$ and
$\bm z_h:=\Pi_{\mathbb{M},h}^{\mathbf{curl}}\bm w_h$.
By \eqref{properties of weighted curl component} and
\eqref{discrete weighted gradient projection},
\begin{equation}\label{curl identity of weighted Helmholtz components}
\mathbf{curl}\,\bm z
=\mathbf{curl}\,\bm z_h
=\mathbf{curl}\,\bm w_h,
\qquad
\mathbf{div}(\mathbb M^H\bm z)=0.
\end{equation}
Applying \cite[Theorem~2.7]{melenk2025} with $\ell=0$,
$\bm\nu=\mathbb M^H$ and $\bm u=\bm z$, and using
\eqref{curl identity of weighted Helmholtz components}, we obtain
\begin{equation}\label{eq:regularity of weighted Helmholtz lift}
\|\bm z\|_{\bm H^1(\Omega\cup\hat\Omega)}
\lesssim
\|\mathbf{curl}\,\bm z\|_{\bm L^2(\mathcal D)}
+\|\mathbf{div}(\mathbb M^H\bm z)\|_{L^2(\mathcal D)}
+\|\bm z\times\hat{\bm x}\|_{\bm H^{1/2}(\hat\Gamma)}
=\|\mathbf{curl}\,\bm w_h\|_{\bm L^2(\mathcal D)}.
\end{equation}
Here the boundary term vanishes because
$\bm z\in\bm H_0(\mathbf{curl},\mathcal D)$.

By \eqref{curl identity of weighted Helmholtz components}--\eqref{eq:regularity of weighted Helmholtz lift}, we have
$\bm z\in\bm H(\mathbf{curl},\mathcal D)
\cap\bm H^1(\mathbf{curl},\Omega\cup\hat\Omega)$.
Hence $\pi_h\bm z$ is well defined. The commuting property of the mapped
interpolants \cite[Eq.~(B.6)]
{chaumontfrelet2025sharperrorboundsedgeelement}, together with
\eqref{curl identity of weighted Helmholtz components}, gives
\begin{align*}
\mathbf{curl}\,\pi_h\bm z
&=\Pi_{h,0}^{\mathrm{RT}}(\mathbf{curl}\,\bm z)
=\Pi_{h,0}^{\mathrm{RT}}(\mathbf{curl}\,\bm z_h)
=\mathbf{curl}\,\bm z_h,
\end{align*}
where $\Pi_{h,0}^{\mathrm{RT}}$ is the corresponding mapped lowest-order
Raviart--Thomas interpolant and the last equality follows from its
reproduction property. Hence
$\mathbf{curl}(\bm z_h-\pi_h\bm z)=\bm0$.
Since the mapped discrete de Rham sequence is exact, see
\cite[Section~7.2.2]{MR2059447} and
\cite[Appendix~B]{chaumontfrelet2025sharperrorboundsedgeelement}, there exists
$\xi_h\in S_h$ such that
\begin{equation}\label{discrete gradient representation of lift error}
\bm z_h-\pi_h\bm z=\nabla\xi_h.
\end{equation}
The continuous and discrete weighted orthogonality, together with
\eqref{discrete gradient representation of lift error}, gives
\begin{align*}
\|\nabla\xi_h\|_{\bm L^2(\mathcal D)}^2\lesssim
\Re(\mathbb M\nabla\xi_h,\nabla\xi_h)=\Re(\mathbb M\nabla\xi_h,\bm z-\pi_h\bm z)\lesssim
\|\nabla\xi_h\|_{\bm L^2(\mathcal D)}
\|\bm z-\pi_h\bm z\|_{\bm L^2(\mathcal D)}.
\end{align*}
Consequently,
\begin{equation}\label{discrete gradient estimate by interpolation error}
\|\nabla\xi_h\|_{\bm L^2(\mathcal D)}
\lesssim
\|\bm z-\pi_h\bm z\|_{\bm L^2(\mathcal D)}.
\end{equation}
By \eqref{low regularity interpolation estimate}, the inverse
estimate for $\mathbf{curl}\,\bm w_h$, and
\eqref{eq:regularity of weighted Helmholtz lift},
\begin{align}
\|\bm z-\pi_h\bm z\|_{\bm L^2(\mathcal D)}
&\lesssim h\|\bm z\|_{\bm H^1(\Omega\cup\hat\Omega)}
+h^2\|\mathbf{curl}\,\bm z\|_{\bm H^1(\Omega\cup\hat\Omega)}\notag\\
&\lesssim h\|\mathbf{curl}\,\bm w_h\|_{\bm L^2(\mathcal D)}.
\label{interpolation estimate for weighted Helmholtz lift}
\end{align}
Finally, by \eqref{discrete gradient representation of lift error}--\eqref{discrete gradient estimate by interpolation error} and  using the triangle inequality, we obtain  \eqref{curl-curlh}.
\end{proof}

\bibliographystyle{plain}
\bibliography{maxwell_pml.bib}
\end{document}